\documentclass[12pt,leqno]{article}

\usepackage{enumitem}

\usepackage[official]{eurosym}
\usepackage{fullpage,amsmath,theorem}
\usepackage{epic,eepic}
\usepackage{amssymb,amsmath}
\usepackage{cases}
\def\llvdash{{\|\hskip-2pt \raise 3pt\hbox{\vrule
			height 0.25pt width 0.4cm}}}

\def\l{{\langle}}
\def\r{{\rangle}}

\def\calK{\mathcal K}

\def\upr{\upharpoonright}

\def\oa{{\overline A^{\,\lower 7pt_{\hbox{$\scriptstyle\bet}}
			\hbox{$\scriptstyle 0\tau$}}}}

\def\bet{\beta}

\def\llvdash{{\|\hskip-2pt \raise 3pt\hbox{\vrule height 0.25pt
			width 0.4cm}}}

\newtheorem{theorem}{Theorem}[section]
\newtheorem{lemma}[theorem]{Lemma}

\newtheorem{defn}[theorem]{Definition}
\newtheorem{question}[theorem]{Question}
\newtheorem{proposition}[theorem]{Proposition}
{\theorembodyfont{\rmfamily}
	\newtheorem{definition}[theorem]{Definition}}
{\theorembodyfont{\rmfamily}
	\newtheorem{remark}[theorem]{Remark}}
{\theorembodyfont{\rmfamily}
	\newtheorem{claim}{Claim}}
{\theorembodyfont{\rmfamily}
	}
{\theorembodyfont{\rmfamily}
	}

\newcommand{\pr}{\medskip\noindent\textit{Proof}. }

\newcommand{\lusim}[1]{\smash{\underset{\raisebox{1.2pt}[0cm][0cm]{$\sim$}}
		{{#1}}}}

\def\llvdash{{\|\hskip-2pt \raise 3pt\hbox{\vrule
			height 0.25pt width 0.15cm}}}

\def\Vdashbks{\hbox{$\Vdash\!\!\!\!{\raise2pt\hbox
			{$\scriptscriptstyle\backslash$}}$}}

\begin{document}

\title{On Easton support iteration of Prikry type forcing notions}
\baselineskip=18pt
\author{Moti Gitik\footnote{The work  was partially supported by ISF grant No 882/22}  and Eyal Kaplan}

\date{\today}
\maketitle

\begin{abstract}
	We consider here  Easton support iterations of Prikry type forcing notions. New ways of constructing normal ultrafilters in extensions are presented.
	It turns out that, in contrast with other supports, seemingly unrelated measures or extenders can be involved here.
\end{abstract}

\section{Introduction}
We continue here the study of the structure of normal ultrafilters in generic extensions by iterated Prikry type forcings.

In \cite{NonStatRestElm, kaplan2022magidor, ben2014forcing}, Nonstationary support and Full support iterations were considered. When iterating Prikry forcings below a measurable limit of measurables $ \kappa $, all the normal measures it carries in the extension are characterized in terms of normal measures in the ground model; furthermore, for every normal measure on $ \kappa $ in the generic extension, the restrictions of its ultrapower to the ground model is an iteration of it by normal measures only.

Here we concentrate on Easton support iteration of arbitrary Prikry type forcings. The situation turned out to be radically different.  Namely, we show the following:
\begin{theorem} \label{Theorem: Conditions imposed on i}
Let $ \kappa $ be a measurable cardinal with $ 2^{\kappa} = \kappa^{+} $. Let $\l P_{\alpha}, \lusim{Q}_{\beta} \colon \alpha\leq \kappa, \beta<\kappa \r$ be an Easton-support iteration of Prikry type forcing notions.

Assume that $ \Delta\subseteq\kappa $ is unbounded, such that for every $ \alpha< \kappa$, $ \lusim{Q}_{\alpha} $ is forced to be trivial if and only if $ \alpha \notin \Delta $. Let $ U \in V $ is a normal measure on $ \kappa $ with $ \Delta\notin U $, and $ i \colon V\to N $ is an elementary embedding, definable in $ V $, such that the following properties hold\footnote{A typical example of such $ N $ is an ultrapower of $ V $ by its $ \kappa $-closed extender, and $ i\colon V\to N $ is its embedding.}:
	\begin{enumerate}
		\item $\mbox{crit}(i) = \kappa$.
		\item $ ^{\kappa}N \subseteq N $.
		\item $\kappa\notin i(\Delta)$.
		\item $ U = \{ X\subseteq \kappa \colon \kappa\in i(X) \} $.
		\item $ \left| i(\kappa) \right| = \kappa^{+} $.
		\item $\{ i(f)(\kappa) \colon f\in V, \ f\colon \kappa\to \kappa \}$ is unbounded in $ i(\kappa) $.
	\end{enumerate}
	Assume also that every element of $ N $ has the form $ i(f)\left( \beta_1, \ldots, \beta_l \right) $ for some $ f\in V $ and $ \beta_1<\ldots < \beta_l < i(\kappa) $.
	 
	 Then there exists a measure $ W\in V\left[  G\right] $ extending $ U $, such that, denoting $ \mbox{Ult}\left( V\left[G\right] , W \right)  \simeq M_W\left[ j_W(G) \right]$, there exists $ k\colon N\to M_W $ with $ \mbox{crit}(k)> \kappa $ such that $ j_W\restriction_{V} = k\circ i $.

\end{theorem}

Furthermore, under mild assumptions on the forcings participating in the iteration, there are  $ \left( 2^{\kappa} \right)^{+} = \kappa^{++} $ normal measures $ W $ as above extending $ U $ (see theorem \ref{Theorem: Kunen-Paris extenstion}). This generalizes the Kunen-Paris theorem on the number of normal measures \cite{kunen1971boolean}. 

In sections \ref{Section: on jwkappa above jukappa} and \ref{Section: properties of k} we analyze the properties of the ultrapower embedding $ j_W \colon V\left[G\right]\to M_W\left[ j_W(G) \right] $ for an arbitrary measure $ W\in V\left[G\right] $. 

Assume that $ V = \mathcal{K} $ is the core model. By a well known series of results in Inner Model theory, $ j_W\restriction_{V} $ is an iterated ultrapower of $ V $, provided that the variety of large cardinals in the universe is limited. For instance, by Mitchell \cite{mitchell1984core}, assuming that there is no inner model with a cardinal $ \alpha $ with $ o(\alpha) = \alpha^{++} $ and $ V = \mathcal{K} $ is the core model, $ j_W\restriction_{\mathcal{K}} $ is an iteration of $ \mathcal{K} $ by normal measures. By a result of Schindler \cite{schindler2006iterates}, assuming that there is no inner model with a Woodin cardinal, $ j_W\restriction_{\mathcal{K}} $ is an iteration of $ \mathcal{K} $ by its extenders.

Theorem \ref{Theorem: Conditions imposed on i} shows that $ j_W\restriction_{V}$ decomposes to the form $ k\circ i $. In particular, $ j_W(\kappa) \geq i(\kappa) \geq j_U(\kappa) $. In section \ref{Section: on jwkappa above jukappa}, we analyze the requirements needed to ensure strict inequality, namely $ j_W(\kappa) > j_U(\kappa) $, by concentrating on the context where $ V = \mathcal{K} $ is the core model and $ j_W\restriction_{\mathcal{K}} $  is an iteration of $ \mathcal{K} $ by measures or extenders.

In section \ref{Section: properties of k} we focus on the question what can be said about the embedding $ k\colon N\to M_W $. In particular, whether it is an iteration of $ N $ by measures or extenders (without  assuming that $ V = \mathcal{K} $ is the core model). We will prove in theorem \ref{Theorem: Restrictions of ultrapwers with simply generated measures} that this is the case where $ P = P_{\kappa} $ is an iteration of Prikry forcings (under some restrictions on the normal measures used; see subsection \ref{Subsection, structure of jW restriction V}). Furthermore, in this case, $ k $ is an iterated ultrapower with normal measures only.

\section{The General Framework} \label{Section: The general Framework}

\begin{defn}
	An iteration $ \langle P_\alpha, \lusim{Q}_\beta \colon \alpha\leq\kappa\ ,\ \beta <\kappa \rangle $ is called an Easton support iteration of Prikry-type forcings if and only if, 	
	\begin{enumerate}
		\item For every $ \alpha<\kappa $, the weakest condition in $ P_{\alpha} $ forces that $ \langle  \lusim{Q}_{\alpha}, \lusim{\leq}_{Q_{\alpha}} , \lusim{\leq}^{*}_{ Q_{\alpha} } \rangle$ is a Prikry type forcing notion.
		\item For every $ \alpha\leq\kappa $ and $p\in P_\alpha$,
	\begin{enumerate}
		\item  $ p $ is a function with domain $ \alpha $ such that for every $\beta <\alpha$, $p\restriction \beta  \in P_\beta$, and $p\restriction \beta \Vdash p(\beta) \in \lusim{Q}_\beta$. 
		\item If $ \alpha \leq \kappa $ is inaccessible, then  $\mbox{supp}(p)\cap \alpha$ is bounded in $\alpha$ ($ \mbox{supp}(p) \subseteq \alpha$ is the set of indices $ \gamma $ on which $ p(\gamma) $ is forced to be non-trivial).
	\end{enumerate}
	\end{enumerate}
	
	Suppose that $p,q \in P_\alpha$. Then $p\geq q$, which means that $p$ extends $q$, holds if and only if:
	\begin{enumerate}
		\item $\mbox{supp}(q)\subseteq \mbox{supp}(p)$.
		\item For every $\beta \in \mbox{supp}(q)$, $p\restriction \beta \Vdash p(\beta)\geq_\beta q(\beta)$ (where $ \geq_\beta$ is the order of $ Q_\beta $).
		\item There is a finite subset $b\subseteq \mbox{supp}(q)$, such that for every $\beta \in \mbox{supp}(q)\setminus b$, $p\restriction \beta \Vdash p(\beta)\geq^{*}_\beta q(\beta)$ (where $ \geq^{*}_{\beta} $ is the direct extension order of $ Q_\beta $).
	\end{enumerate} 
	If $b=\emptyset$, we say that $p$ is a direct extension of $q$, and denote it by $p\geq^* q$. 
\end{defn}

The following properties are standard (see \cite{gitik2010prikry} for example):

\begin{lemma}
	For every $ \lambda\leq \kappa $, $ P_{\lambda} $ satisfies the Prikry property.
\end{lemma}

\begin{lemma} \label{Lemma: lambda mahlo implies lambda cc}
	For every $ \lambda\leq \kappa $ which is Mahlo, $ P_{\lambda} $ has the $ \lambda-c.c. $.
\end{lemma}

Let $U$ be a normal ultrafilter over $\kappa$. Let $\l P_\alpha, \lusim{Q}_\beta\mid \alpha\leq\kappa, \beta<\kappa\r$ be
an Easton support iteration of a Prikry type forcing notions. Suppose that the following hold:

\begin{enumerate}
	\item
	There exists an unbounded subset $ \Delta\subseteq \kappa $, $ \Delta \notin U $, such that, for every $ \alpha< \kappa$,
	\begin{enumerate}
		\item 	$ \alpha\in \Delta \  \longrightarrow \  \Vdash_{ P_{\alpha} }  \lusim{Q}_\alpha \mbox{ is nontrivial.}$
		\item 	$ \alpha\notin \Delta \  \longrightarrow \  \Vdash_{ P_{\alpha} }  \lusim{Q}_\alpha \mbox{ is trivial.}$
	\end{enumerate}

	\item For every $\alpha<\kappa$, $\Vdash_{P_\alpha} \l \lusim{Q}_\alpha, \lusim{\leq}_\alpha^*\r$ is $\alpha-$closed.
	
	\item For every $\alpha\in \Delta$, $\Vdash_{P_\alpha} | \lusim{Q}_\alpha|<\min(\Delta\setminus \alpha+1)$.
	

	

\end{enumerate}

Let $G$ be a generic subset of $P=P_\kappa$. We would like to analyze the normal measures on $ \kappa $ in $ V\left[G\right] $ extending $ U $. The standard way to do so appears in \cite{gitik2010prikry}, we present it here for sake of completeness.

\begin{lemma} \label{Extending measure using * increasing sequence}
There exists a normal measure $ U^*\in V\left[G\right] $ on $ \kappa $ which extends $ U $.
\end{lemma}

\pr
Let $ \langle \lusim{A}_{\alpha} \colon \alpha<\kappa^{+} \rangle $ be an enumeration, in $ V $, of $ P = P_{\kappa} $-names, such that every $ X\in \left( \mathcal{P}(\kappa) \right)^{V\left[G\right]} $ has the form $ \left( \lusim{A}_{\alpha} \right)_G $ for some $ \alpha<\kappa^{+} $. Such list of names exists since $ P = P_{\kappa} $ is $ \kappa-c.c. $. Now, construct, in $ V\left[G\right] $, a $ \leq^* $-increasing sequence of conditions $ \langle q_{\alpha} \colon \alpha<\kappa^{+}\rangle $, such that, over $ N\left[G\right] $, $ q_{\alpha} \parallel \kappa\in j_U\left( \lusim{A}_{\alpha} \right) $. Such a sequence exists since $ V\left[G\right] \vDash \mbox{"}\langle j_U(P)\setminus \kappa, \leq^*\rangle \mbox{ is } \kappa^{+}-\mbox{closed}$."

Let $ \langle \lusim{q}_{\alpha} \colon \alpha<\kappa^{+} \rangle $ be a $ P $-name for the above sequence. Now, define $ U^*\supseteq U $ as follows: For every $ \alpha<\kappa^{+} $, $\left(\lusim{A}_{\alpha}\right)_G \in U^*$ if and only if there exists $ p\in G $ and $ \alpha<\kappa^{+} $ such that-- 
$$ p ^{\frown} \lusim{q}_{\alpha} \vDash \kappa \in i(\lusim{A}_{\alpha}) $$
 We argue that $ U^* $ defined above is a normal measure which extends $ U $. 

Assume that $ \delta<\kappa $ and $ \langle \lusim{X}_{\alpha} \colon \alpha<\delta \rangle $ is a $ P_{\kappa} $-name for a partition of $ \kappa $ in $ V\left[G\right] $. For every $ \alpha<\delta $, define--
$$ Y_{\alpha} = \{ \beta<\kappa^{+} \colon \exists p\in P_{\kappa}, \ p\Vdash \lusim{X}_{\alpha} = \lusim{A}_{\beta} \} $$
Since $ P $ is $ \kappa-c.c. $, $ \left| Y_{\alpha} \right|<\kappa $. Denote-- 
$$ Y = \bigcup_{ \alpha<\delta } Y_{\alpha} $$
Then $ Y\subseteq \kappa^{+} $ is a bounded subset. Pick $ \alpha^*<\kappa^{+} $ high enough which bounds $ Y $. Let us argue that there exists $ p\in G $ and a unique $ \beta<\delta $ such that--
$$ p^{\frown} \lusim{q}_{\alpha^*} \Vdash \kappa \in j_U\left( \lusim{A}_{\beta} \right) $$
and thus $ \left(  \lusim{A}_{\beta}\right)_G \in U^* $.

Work in $ N\left[G\right] $. Note that $ \langle A_{\beta} \colon \beta\in Y \rangle $ covers the sequence $ \langle X_{\alpha} \colon \alpha<\delta \rangle $. Since $ {q}_{\alpha^*} $ is $ \leq^* $ above any $ q_{\beta} $ for $ \beta\in Y $, 
$$\forall \xi< \alpha,   q_{\alpha^*} \parallel  \kappa\in i(\lusim{X}_{\xi}) $$
Since $ \langle i\left(X_{\xi}\right) \colon \xi<\delta \rangle $ is a partition of $ i(\kappa) $, there exists a unique $ \xi^*<\delta $ such that  $ q_{\alpha^*}\Vdash \kappa\in i\left( \lusim{A}_{\xi^*}  \right) $.  Let $ p\in G $ be a condition forcing this. Then $ p^{\frown} \lusim{q}_{\alpha^*} \Vdash \kappa\in i\left( \lusim{X}_{\xi^*} \right) $, as desired.

A similar argument shows that $ U^* $ is normal. Indeed, given a $ P_{\kappa} $-name for a regressive function $ \lusim{f}\colon \kappa\to \kappa $, define, for every $ \alpha< \kappa $,
$$ X_{\alpha}  =  \{ \xi<\kappa \colon f(\xi) = \alpha \}$$
and proceed as before to find a unique $ \alpha<\kappa $ such that $ X_{\alpha}\in U^* $.
$\square$\\

In particular, $ U $ can be extended to a normal measure $ U^*\in V\left[G\right] $, such that the ultrapower embedding $ j_{U^*} \colon V\left[G\right]\to M\left[ j_{U^*}(G) \right] $ satisfies that $ j_{U^*}\restriction_{V} = k\circ j_U $, for an embedding $ k\colon M_U \to M $ which satisfies $ \mbox{crit}(k) > \kappa $. Indeed, define $ k\left( \left[f\right]_U \right) = \left[f\right]_{U^*} $ for every $ f\colon \kappa\to V $ in $ V $.

A natural question here is whether this is the only way to generate a normal ultrafilter on $ \kappa $ in $V[G]$. In \cite{ RestElm,NonStatRestElm} it was established that this is the case when considering the Nonstationary support iteration. However, this is not true anymore once full support iterations are considered: in   \cite{ben2014forcing}  and later in \cite{kaplan2022magidor}, iterations of the standard Prikry forcing were considered. It was proved that every normal measure $ U\in V $ on $ \kappa $ with $ \Delta\notin U $ can be extended to a normal measure $ U^*\in V\left[G\right] $ similarly as above, but not every normal measure extending $ U $ is generated this way; nevertheless, all the normal measures on $ \kappa $ in $ V\left[G\right] $ were characterized, either as extensions $ U^* $ of measures $ U\in V $ with $ \Delta\notin U $, or as the projections to normal measure of extensions $ U^* $  of a normal ultrafilter $ U\in V $ with $ \Delta\in U $.

It turns out that the picture in the Easton support iteration of Prikry type forcing notions (and even of the standard Prikry forcings) is radically different. Given an elementary embedding $ i\colon V\to N $ with critical point $ \kappa $, definable in $ V $, the normal measure derived from it, $ U = \{ X\subseteq \kappa \colon \kappa\in i(X) \} $, can be extended to a normal measure $ W\in V\left[G\right] $ such that $ j_W\restriction_{V} = k\circ i $, for some $ k\colon N\to M $ with $ \mbox{crit}(k) > \kappa $. In the case of iterations of the standard Prikry forcing, $ k $ is an iterated ultrapower of $ N $ by normal measures only (see section \ref{Section: properties of k}), while $ i\colon V\to N $ can be an embedding derived from an extender (as in the formulation of theorem \ref{Theorem: Conditions imposed on i}).

Let us demonstrate that, in the Easton support iteration, there are many more possibilities to get normal measures $ W\in V\left[G\right] $. We show that an arbitrary embedding $ i\colon V\to N $ can be used to extend the normal measure $ U $ derived from it.

\begin{lemma} \label{extending U using generic} 
	Assume that $ i \colon V\to N $ is an elementary embedding definable in $ V $, with $ \mbox{crit}(i) = \kappa $, such that $ \left| i(\kappa)\right| = \kappa^{+} $, $ \kappa\notin i(\Delta) $, $ N\subseteq V $ and $^{\kappa} N \subseteq N $. Denote--
	$$ U = \{ X\subseteq \kappa \colon \kappa\in i(X) \} $$
	 Then $ G $ is $ i(P)\restriction_{\kappa} = P $-generic over $ N $, and:
	\begin{enumerate}
		\item For every $q \in i(P)\setminus \kappa$, there is $H\in V\left[G\right]$  with $ q\in H $, which is $\l i(P)\setminus \kappa, \leq^*\r$-generic over $N\left[G\right]$.
		\item Given such $ H\in V\left[G\right] $, define--
		\begin{align*}
		U_{H} = \{  \left(\lusim{A}\right)_G \colon & \lusim{A} \mbox{ is a } P-\mbox{name for a subset of } \kappa, \mbox{ and there exists }\\
		&p\in G*H \mbox{ such that } p\Vdash \kappa\in i\left( \lusim{A} \right) \} 
		\end{align*}
		Then $U_H$ is a normal, $\kappa-$complete ultrafilter on $\kappa$ which extends $ U $.
	\end{enumerate}
\end{lemma}

\pr  
\begin{enumerate}
	\item We can enumerate, in $ V\left[G\right] $, all the maximal antichains in $ \langle i(P)\setminus \kappa , \leq^* \rangle $ with order type $ \kappa^{+} $, by $ i(\kappa) $-c.c. of the forcing, and since $ V\left[G\right]\vDash \left| i(\kappa) \right| = \kappa^{+} $. Note that $\kappa\not\in i\left(\Delta\right)$, so in the sense of $ N\left[G\right] $, the forcing $ \langle i(P)\setminus \kappa , \leq^* \rangle $ is more than $ \kappa $-closed. Moreover, since $ V\vDash {}^{\kappa} N \subseteq N $, and $ P = P_{\kappa} $ is $ \kappa- $c.c., $ V\left[G\right]\vDash {}^{<\kappa}N\left[G\right]\subseteq N\left[G\right]$. Therefore, every sequence of length $ \kappa $ of conditions in $ i(P)\setminus \kappa $ which belongs to $ V\left[G\right] $ belongs to $ N\left[G\right] $ as well. Thus, in the sense of $ V\left[G\right] $, the forcing $ \langle i(P)\setminus \kappa , \leq^* \rangle $ is $ \kappa^{+} $-closed. 
	
	Starting from any condition in $ i(P)\setminus \kappa $, we can construct (in $ V\left[G\right] $) a sequence of direct extensions of it, meeting every maximal antichain. 
	This sequence generates a $ \leq^* $-generic over $ N\left[G\right] $ for $ i(P)\setminus \kappa $, which belongs to $ V\left[G\right] $.
	\item First, we prove that $ W = U_{H} $ is a normal, $ \kappa $-complete ultrafilter on $ \kappa $ which extends $ U $. It is not hard to verify that $ W $ is a filter. We prove that $ W $ is a $ \kappa $-complete ultrafilter. Assume that $ \langle X_{\alpha} \colon \alpha<\delta \rangle $ is a partition of $ \kappa $, for some $ \delta<\kappa $. Work in $ N\left[G\right] $. Let $ D\subseteq i(P)\setminus \kappa $ be the $ \leq^* $-dense open set of conditions which decide the unique $ \alpha<\delta $ for which $ \kappa\in i\left( \lusim{X}_{\alpha} \right) $. Then such a statement is forced by some $ r\in H $. Let $ p\in G $ be a condition which forces that $ r $ has this property, and also decides the value of $ \alpha $. Then $ p^{\frown} r \Vdash \kappa\in i\left( X_{\alpha} \right) $ and thus $ X_{\alpha}\in W $. Normality of $ W $ follows by a similar argument, using the dense set of conditions deciding the value of $ i(\lusim{f})(\kappa) $ for a given regressive function $ f\colon \kappa\to \kappa $. The argument works since we don't force over $ \kappa  $ in $ N $.
\end{enumerate}$\square$\\

\begin{remark}
	M. Magidor pointed out the following: Assuming that $ N\subseteq V $ and $ i\colon V\to N $ is definable in $ V\left[G\right] $, it follows that $N$ is already a class of $ V $.
	Indeed, pick a formula $ \varphi $ and a parameter $ a\in V\left[G\right] $ such that for every $ x,y $ in $ V $, $ \varphi(x,y,a) $ holds in $ V\left[G\right] $ if and only if $ i(x) = y $. For every ordinal $ \alpha $ pick a condition $ p_{\alpha}\in G $ which decides the value of the set $ \left( V_{i\left(\alpha\right)} \right)^{N} $, which is the set $ y $ for which $ \varphi\left( V_{\alpha}, y, a \right) $ holds. 
	Since $ P $ is a set forcing, there exists $ p^*\in G $ such that, for unboundedly many ordinals $ \alpha $, $ p_{\alpha} = p^* $. Then $ N $ can be defined as a class of $ V $ using $ p^* $, 
	$ N = \bigcup\{ y \colon \exists \alpha\in \mbox{ON}, \  p^* \Vdash \varphi\left( V_{\alpha}, y, \lusim{a} \right) \} $.
\end{remark}

In general, the settings of lemma \ref{extending U using generic} are not enough ensure that $ j_{U_H}\restriction_{V} = k\circ i $ for some $ k $ with $ \mbox{crit}(k) > \kappa $. For instance, given a normal measure $ U $ on $ \kappa $ in $ V $ with $ \Delta\notin U $, the embedding $ i = j_{U^2} $ satisfies the settings of lemma \ref{extending U using generic}, but cannot be used to extend $ U $ to a measure $ U_{H} $ for which $ j_{U_H} = k\circ i $ for some embedding $ k $ with $ \mbox{crit}(k)> \kappa $. This follows since $ i $ fails to satisfy clause $ 3 $ in the next claim:

\begin{proposition} \label{Proposition: Necessary conditions}
Assume that $ U\in V $ is a normal measure on $ \kappa $,  $ W\in V\left[G\right] $ is a normal measure which extends $ U $, $ i \colon V\to N $ is an elementary embedding and $ j_W\restriction_{V} = k\circ i $ for some $ k\colon N\to M $ with $ \mbox{crit}(k) > \kappa $. Then--
\begin{enumerate}
	\item $\{X\subseteq \kappa \colon  \kappa \in i(X)\}=U$.
	\item $ \left|i(\kappa)\right| = \kappa^{+} $.
	\item $\{i(f)(\kappa)  \colon f\in V, \  f:\kappa \to \kappa \}$ is unbounded in $i(\kappa)$.
\end{enumerate}
\end{proposition}

\pr 
\begin{enumerate}
	\item $\{X\subseteq \kappa \colon  \kappa \in i(X)\}=U$: Indeed, let $ X\subseteq \kappa $ in $ V $ with $ \kappa\in i(X) $. By applying $ k\colon N\to M $ it follows that $ \kappa\in j_W(X) $ and hence $ X\in W $. Since $ X\in V $ and $ U = W\cap V $, it follows that $ X\in U $.
	\item $ \left|i(\kappa)\right| = \kappa^{+} $: This holds since, in $ V\left[G\right] $, $ \left| j_W(\kappa) \right| = 2^{\kappa} = \kappa^{+} $ (since, in $ V $, $ 2^\kappa = \kappa^{+} $), and $ i(\kappa)\leq j_W(\kappa) $.
	\item $\{i(f)(\kappa)  | f:\kappa \to \kappa \}$ is unbounded in $i(\kappa)$: Given $ \beta<i(\kappa) $, let $ f\in V\left[G\right] $ be a function such that $ \left[f\right]_{W} = k\left(\beta\right) $. Since $ k\left(\beta \right)<k\left( i(\kappa) \right)= j_W(\kappa) $, we can assume that $ f\colon \kappa \to \kappa $. The Easton support ensures that there exists $ g\colon \kappa \to \kappa $ in $ V $ which dominates $ f $. Thus $ i\left(g\right)(\kappa) \geq \beta  $ (indeed, by applying $ k\colon N\to M $ on both sides, this is equivalent to $ j_W(g)(\kappa) \geq k\left( \beta  \right) = \left[f\right]_W $, which holds, since $ g $ dominates $ f $. Note that, when applying $ k $, we used the fact that $ \mbox{crit}(k)> \kappa $).
\end{enumerate}
$\square$

Theorem \ref{Theorem: Conditions imposed on i} will be proved by a sequence of lemmata, concluded in lemma \ref{Lemma: k is elementary}. The main idea in the proof of theorem \ref{Theorem: Conditions imposed on i} is to add representing functions for all the generators of $ i $ above $ \kappa $. This is needed since $ j_W\restriction_{V} $ has a single generator $ \kappa $.

\begin{definition}
	An ordinal $ \beta $ is called a generator of $ i\colon V\to N $ if there are no $ n<\omega $, ordinals $ \beta_1,\ldots, \beta_n $ below $ \beta
	$ and a function $ f\in V $ such that $ \beta = i\left( f \right)\left( \beta_1, \ldots, \beta_n \right) $.
\end{definition}

In the next lemma we construct a function $ \alpha\mapsto \theta_{\alpha} $ in $ V\left[G\right] $, which will be utilized, alongside functions in $ V $, to represent the generators of $ i $ in $ \mbox{Ult}\left( V\left[G\right], W \right) $.

\begin{lemma} \label{Existence of the sequence  theta alpha}
There exists a $ P_{\kappa} $-name for a sequence of ordinals, $ \langle \theta_{\alpha} \colon \alpha< \kappa \rangle $, such that the following property holds:\\
\begin{enumerate}
	\item For every $ \beta<\kappa $ and $ p\in P_{\kappa} $, there is $ \alpha_0 < \kappa $ such that for every $ \alpha\geq \alpha_0 $ there exists $ p^*\geq^* p $ such that $ p^*\Vdash  \ \theta_{\alpha} = \beta $.
	\item For every $ \alpha<\kappa $ and condition $ p\in P_{\kappa} $, there exists a condition $ p^*\geq^* p $ which decides the value of $ \lusim{\theta}_{\alpha} $.
\end{enumerate}
\end{lemma}

\begin{remark}
When iterating Prikry forcings, the natural candidate for the function $ \alpha\mapsto \theta_{\alpha} $ is the function which maps every $ \alpha\in \Delta $ to the first element in its Prikry sequence (this function does not have domain $ \kappa $, but this can be fixed by defining the function of elements outside of $ \Delta $ as the value of the first element of $ \Delta $ above them). The main problem with such a function is that it fails to satisfy clause $ 2 $ of the lemma (from density, every $ \leq^* $-generic set has some $ \alpha\in \Delta $ for which it does not decide the first element in its Prikry sequence).  However, we can base the function $\alpha \mapsto \theta_\alpha$ on $\alpha \mapsto \mbox{the first element of the Prikry sequence of } \alpha$. There will be still many places of disagreement between them,  but for every given $\beta<\kappa$ there will be $\alpha<\kappa$ such that $\beta$ is the first element of the Prikry sequence of $\alpha$.
\end{remark}

\pr 
For every $ \alpha<\kappa $, let $ \tau_{\alpha} <\kappa $ be the least ordinal such that $ P\restriction_{ \left( \alpha, \tau_{\alpha} \right)} $ is not $ \alpha-c.c. $. We will argue below that such $ \tau_{\alpha}<\kappa$ exists, but first, let us show that this suffices: Pick an unbounded subset $ X\subseteq \kappa $, such that, for every $ \alpha, \alpha'\in X $,
$$ \alpha<\alpha' \implies \tau_{\alpha}< \tau_{\alpha'} $$
(for instance, let $ X $ be the club of closure points of the function $ \alpha\mapsto \tau_{\alpha} $).
Enumerate $ X = \langle x_{\alpha} \colon \alpha\in \Delta \rangle $. For every $ \alpha\in \Delta $, let $ \langle q_{x_\alpha, \xi} \colon \xi<x_{\alpha} \rangle $ be an antichain in $ P_{ \left( x_{\alpha}, \tau_{x_{\alpha}} \right) } $ of cardinality $ \alpha $. Define $ \theta_{\alpha} $ to be the unique ordinal $ \xi<x_{\alpha} $ for which $ q_{x_\alpha, \xi}\in G\restriction_{ \left( x_{\alpha}, \tau_{x_{\alpha}} \right) }  $ (if there is no such $ \xi $, which is possible since the antichain is not necessarily maximal, set $ \theta_\alpha =0 $). 

Now, given $ \beta<\kappa $ and a condition $ p\in P_{\kappa} $, pick first $ \alpha\in \Delta $ for which $ x_{\alpha} $ bounds the support of $ p $. Direct extend $ p $ to $ p^* $ such that $ p^*\restriction_{ \left( x_{\alpha}, \tau_{x_{\alpha}} \right) } = q_{x_{\alpha}, \beta}  $. Then\footnote{ It is crucial here that the Easton support is used.}
 by our definition, $ p^* $ forces that $ \theta_{\alpha} = \beta $. 

Let us prove now that for every $ \alpha<\kappa $ and condition $ p\in P_{\kappa} $, there exists $ p^*\geq^* p $ which decides the value of $ \theta_{\alpha} $. 

We will direct extend $ p $ is the interval $ \left( x_{\alpha}, \tau_{ x_{\alpha} } \right) $, $ x_{\alpha} $-many times, to decide whether $ q_{ x_{\alpha}, \xi }\in G\restriction_{ \left( x_{\alpha}, \tau_{x_{\alpha}} \right)} $, for every $ \xi<x_{\alpha} $. Note that this is possible since $ \l P\restriction_{ \left( x_{\alpha}, \tau_{ x_{\alpha} } \right) } , \leq^* \r $ is more than $ x_{\alpha} $-closed. Let $ p^*\geq^* p $ be the obtained condition. Then either there exists $ \xi< x_{\alpha} $ such that $ p^*  $ forces that $ q_{ x_{\alpha}, \xi } $ is in the generic, and then $ p^*\Vdash \lusim{\theta}_{\alpha} = \xi $; or, there is no such $ \xi $, and then $ p^*\Vdash \lusim{\theta}_{\alpha} = 0 $.

Let us argue now that indeed, for every $ \alpha<\kappa $ there exists  $ \tau_{\alpha}<\kappa $ such that $ P\restriction_{ \left( \alpha, \tau_{\alpha} \right) } $ is not $ \alpha-c.c. $: Pick $ \tau_{\alpha} $ such that there are $ \alpha $-many elements of $ \Delta $ in the interval $ \left( \alpha, \tau_{\alpha} \right) $. Let $ \langle \tau_{\alpha,\xi} \colon \xi<\alpha \rangle $ be an enumeration of the first $ \alpha$-many elements in $ \left(  \alpha,\tau_{\alpha} \right)\cap \Delta $. For every $ \xi<\alpha $, let $ \lusim{x}_{\xi}, \lusim{y}_{\xi} $ be $ P_{\tau_{\alpha,\xi}} $-names, forced by $ 0_{P_{ \tau_{\alpha,\xi}  }} $ to be pair of incompatible elements of $ \lusim{Q}_{\tau_{\alpha,\xi} } $. Such a pair exists since $ \lusim{Q}_{ \tau_{\alpha, \xi} } $ is nontrivial. 

Now, for every $ \sigma\in 2^{\alpha} $, let $ p_{\sigma}\in P\restriction_{ \left( \alpha, \tau_{\alpha} \right) } $ be the condition which satisfies, for every $ \xi<\alpha $, that-- 
$$ p_{\sigma}\restriction_{\xi} \Vdash p_{\sigma}\left( \xi \right) =  \begin{cases}
	\lusim{x}_{\xi} &  \mbox{ If } \sigma\left( \xi \right) =0 \\
	\lusim{y}_{\xi} &  \mbox{ If }\sigma \left( \xi \right) = 1
\end{cases}$$
Note that  $ \tau_{\alpha} $ is the limit of the first $ \alpha $ many elements above $ \alpha $ in $ \Delta $, and thus $ \tau_{\alpha} $ is singular, so the support of a condition in $ P = P_{\kappa} $ may be unbounded in $ \tau_{\alpha} $.

Then $ \langle p_{\sigma} \colon \sigma \in 2^{<\alpha} \rangle $ is an antichain in $ P\restriction_{\left(  \alpha, \tau_{\alpha} \right)} $ of cardinality at least $ \alpha $.
$\square$\\

\begin{remark} \label{Remark: Abuse of notation with theta}
Given a function $ \alpha\mapsto \lusim{\theta}_{\alpha} $ as in lemma \ref{Existence of the sequence  theta alpha}, we slightly abuse the notation and denote $ i\left(  \alpha\mapsto \lusim{\theta}_{\alpha} \right) $ by $ \l \lusim{\theta}_{\alpha} \colon \alpha<i(\kappa) \r  $. 
\end{remark}

\begin{lemma} \label{Lemma: Genereic satisfying *}
Under the assumptions of theorem \ref{Theorem: Conditions imposed on i},  there exists $ H\in V\left[G\right] $ which is $ \langle i(P)\setminus \kappa, \leq^* \rangle $-generic over $ N\left[G\right] $, with the following property:
\begin{align*}
	(*)  \ \  & \mbox{For every generator } \ \beta\in i(\kappa)\setminus \left(\kappa+1\right)  \mbox{ of } i, \mbox{ there exists a function } f = f_{\beta}\in V, \\
	&f\colon \kappa \to \kappa \ \mbox{and a condition } q\in H \mbox{ such that }  q\Vdash   \beta =  i\left(  \alpha\mapsto \lusim{\theta}_{ f(\alpha)} \right)\left(\kappa  \right) .
\end{align*}
where $ \langle \lusim{\theta}_{\alpha} \colon \alpha < i(\kappa) \rangle $ is as in remark \ref{Remark: Abuse of notation with theta}.
\end{lemma}

\pr In $ V\left[G\right] $, let $\l A_\xi \mid \xi<\kappa^+\r$ be an enumeration of maximal antichains in $i(P)$. Let $\l \beta_\xi \mid\xi<\kappa^+\r$ be an enumeration of all the generators of $ i $ below  $i(\kappa)$. Define in $V[G]$ a $\leq^*-$increasing sequence $\l r_\xi\mid \xi<\kappa^+\r$. Assume that $ \l r_{\xi} \colon \xi < \xi^* \r $ has been constructed for some $ \xi^*< \kappa^{+} $. Pick a condition $ r $ which $ \leq^* $ extends all the conditions $ \langle r_{\xi} \colon \xi < \xi^* \rangle $ constructed so far, and, by extending it, assume that $ r $ extends a condition in $ A_{\xi^*} $. Finally, let $ \alpha_0 < i(\kappa)  $ be such that for every $ \alpha\geq \alpha_0 $ there exists $ r^* \geq^* r $ which forces that $ i\left( \xi\mapsto \theta_{\xi} \right)\left(  \alpha  \right) = \beta_{\xi^*} $. Pick any $ \alpha \geq \alpha_0 $ below $ i(\kappa) $ which has the form $ i\left( f \right)\left( \kappa \right) $ for some $ f = f_{ \beta_{\xi^*} } \in V $, and let $ r_{\xi^*} \geq^* r $ be a condition which forces that $ i\left( \xi\mapsto \theta_{\xi} \right)\left(  \alpha  \right) = \beta_{\xi^*} $.

Finally, let $ H $ be the $ \leq^* $-generic generated from $ \langle r_{\xi} \colon \xi <\xi^* \rangle $.
$\square$\\

\begin{remark}  \label{Remark: Kunen Paris}
Repeating the above argument, we can construct $ 2^{\kappa^{+}} $-many distinct generic sets $ H $ satisfying property $ (*) $, by constructing a binary tree $ \langle r_{ \sigma }  \colon \sigma\in 2^{<\kappa^{+}}\rangle $ of conditions, which are $ \leq^* $-increasing in each branch, and for each $ \sigma\in 2^{<\kappa^{+}} $, $ r_{ \sigma^{\frown} \langle 0 \rangle } $ and $ r_{ \sigma^{\frown} \l1\r  } $ are $ \leq^* $-incompatible. 
Assuming $ 2^{\kappa^{+}} = \kappa^{++} $, this provides the maximal number of generic sets $ H $ in $ V\left[G\right] $ for $ \l  i(P)\setminus \kappa, \leq^*\r $ over $ N\left[G\right] $. 

Below we will define for every such $ H $ a measure $ U_{H}\in V\left[G\right] $ on $ \kappa $ which extends $ U $; under mild assumptions on the forcing notions $ \lusim{Q}_{\alpha} $, we will prove that for $ H\neq H' $ satisfying property $ \left( * \right) $, $ U_{H} \neq U_{  H'} $ (see theorem \ref{Theorem: Kunen-Paris extenstion}). Assuming $ \mbox{GCH} $, this produces the maximal number $ \kappa^{++} $ of normal measures on $ \kappa $, generalizing the well known result of Kunen and Paris \cite{kunen1971boolean}.
\end{remark}

\begin{remark}
Not every generic set $ H\in V\left[G\right] $ for $ \l i(P)\setminus \kappa, \leq^* \r $ satisfies property  $ (*) $.

Indeed, assume that $ \Delta $ consists only of inaccessibles and $ i\colon V\to N $ has a nonempty set of generators in $ \left( \kappa, i(\kappa) \right)$ which is bounded by some ordinal $ \eta =i(f)(\kappa) $ below $ \min\left( i(\Delta)\setminus \kappa \right) $, for some $f\in V$. This holds in the typical case where $ \Delta $ consists of measurables below $ \kappa$ and $ i $ is a $ \left( \kappa, \kappa^{+} \right) $-extender (the length of $ i $ is $ \kappa^{+} $ since $ i $ has to satisfy the requirement $ \left| i(\kappa) \right| = \kappa^{+} $ of theorem \ref{Theorem: Conditions imposed on i}).  Let $ \sigma \colon M_U \to N $ be the embedding which maps each element $ \left[g\right]_U $ of $ M_U $ to $ i(g)(\kappa) $ (here $ g\in V $ is any function with domain $ \kappa $). $ \sigma $ has critical point strictly above $ \kappa^{+} $, since $ \left( \kappa^{+} \right)^N = \kappa^{+}$.

In $ V\left[G\right] $, let $ H_U \subseteq j_U(P)\setminus \kappa $ be $ \leq^* $-generic over $ M_U\left[G\right] $. Let $ H\subseteq i(P)\setminus \kappa  $ be the generic set generated from $ \sigma'' H_U $. We argue that $ H $ is indeed $ \leq^* $-generic over $ N $. Let $ D\in N\left[G\right] $ be a $ \leq^* $-dense open subset of $ i(P)\setminus \kappa $. Write $ \lusim{D} = i(F)\left( \kappa, \beta_1, \ldots, \beta_l \right) $ for some function $ F\in V $, $ l<\omega $ and generators $ \beta_1, \ldots, \beta_l < i(f)(\kappa) $ of $ i $. We can assume that for every $ \xi,\eta_1, \ldots, \eta_l < f(\xi) $, $ F\left( \xi, \eta_1, \ldots, \eta_l \right)\subseteq P\setminus \xi $ is forced to be $ \leq^* $-dense open subset of $ P\setminus \xi $. Define, in $ M_U $,
$$ D_U = \bigcap_{ \gamma_1 , \ldots , \gamma_l < j_U(f)(\kappa)} j_U(F)\left( \kappa, \gamma_1, \ldots,\gamma_l \right) $$
and note that, since the amount of sequences $ \gamma_1, \ldots, \gamma_l < j_U(f)(\kappa)  $ in $ M_U $ is below $ \min\left( \Delta\setminus \kappa \right) $, and $ \l j_U(P)\setminus \kappa , \leq^* \r $ is more than $ \min\left( j_U(\Delta)\setminus \kappa \right) $-closed, $ D_U $ is $ \leq^* $-dense open subset of $ j_U(P)\setminus \kappa $. Pick any $ q\in H_U \cap D_U $. Then $ \sigma(q)\in D\cap H $, since $   \sigma(D_U)  \subseteq D$. 

Since $ \sigma'' G*H_U \subseteq G*H $, the embedding $ \sigma \colon M_U\to N $ can be lifted to an embedding $ \sigma^* \colon M_U\left[ G*H_U \right]\to N\left[ G*H \right] $.

Pick now any generator $ \beta $ of $ i $ in the interval $ \left( \kappa, i(\kappa) \right) $. We argue that there is no $ f\in V $ such that $ H\Vdash \beta = i\left( \alpha\mapsto \lusim{\theta}_{ f\left( \alpha \right) } \right)(\kappa) $. Indeed, otherwise, by elementarity of $ \sigma^* $, there exists $ \beta^*< j_U(\kappa) $ such that--
$$ H_U \Vdash \beta^* = j_U\left( \alpha\mapsto \lusim{\theta}_{f\left( \alpha \right) } \right)\left( \kappa \right) $$
Let $ g\in V $ be a function such that $ \beta^* = j_U(g)(\kappa) $. Then--
$$ \beta = \sigma^*\left( \beta^* \right) = i\left( g \right)\left( \kappa \right) $$
contradicting the fact that $ \beta $ is a generator of $ i $.
\end{remark}

Given $ i, N , U $ as in theorem \ref{Theorem: Conditions imposed on i} and a generic set $ H\in V\left[G\right] $ for $ \langle i(P)\setminus \kappa , \leq^*\rangle $ over $ N\left[G\right] $, define--
\begin{align*}
	U_{H} = \{  \left(\lusim{A}\right)_G \colon & \lusim{A} \mbox{ is a } P-\mbox{name for a subset of } \kappa, \mbox{ and there exists }\\
	&p\in G*H \mbox{ such that } p\Vdash \kappa\in i\left( \lusim{A} \right) \} 
\end{align*}
Then $ U_{H} $ is a normal, $ \kappa $-complete ultrafilter which extends $ U $. This follows by repeating the argument of lemma \ref{extending U using generic}. 

The model $M_{U_H} \simeq \mbox{Ult}\left( V\left[G\right], U_{H} \right)$ is of the form $M[G^*]$, where $M$ is the image of $V$ and $G^*=j_{U_H}(G)$ is $j_{U_H}(P)-$generic over $M$ in sense of $M_{U_H}$. We conclude the proof of theorem \ref{Theorem: Conditions imposed on i} by defining an elementary embedding $ k\colon N\to M $ and proving that $ \mbox{crit}(k)> \kappa $.

In the next lemma we continue the abuse of notation as in remark \ref{Remark: Abuse of notation with theta}, and denote--
$$ j_{U_{H}}\left(  \l \theta_{\xi} \colon \xi\in \Delta \r \right) = \langle \theta_{\xi} \colon \xi\in j_{U_H}\left( \Delta \right) \rangle $$

\begin{lemma} \label{Lemma: k is elementary}
Assume the settings of theorem \ref{Theorem: Conditions imposed on i}. 	Suppose that $H$ is a generic set for $ \langle i(P)\setminus \kappa, \leq^* \rangle $ over $ N\left[G\right] $ with the property $ (*) $. 	Define then $k:N\to M$ as follows:
$$k\left(i(f)(\kappa, \beta_1,...,\beta_l)\right)=j_{U_H}(f)\left(\kappa, \theta_{ \left[ f_{ \beta_1 }  \right]_{ U_H }  },..., \theta_{ \left[ f_{ \beta_l }  \right]_{U_H} }\right)$$
For every $ l<\omega $, $ \beta_1, \ldots, \beta_l <i(\kappa) $ generators of $ i $ and $ f\in V $ (the functions $ f_{\beta_i} $, $ 1\leq i\leq l $, are as in lemma \ref{Lemma: Genereic satisfying *}).

Then $ k\colon N\to M $ is elementary, $\mbox{crit}(k) > \kappa$ and $ j_{U_H}\restriction_{V} = k\circ i $.
\end{lemma}

\pr
Denote $ W  = U_{H} $. Let us prove that the embedding $ k $ defined above is elementary. Assume that $ x,y\in N $.  There are functions $ f,g $ in $ V $, generators $ \beta_1<\ldots< \beta_l < i(\kappa) $ such that--
$$ x= i(f)\left( \kappa,\beta_1, \ldots, \beta_l \right) , \ y = i(g)\left(\kappa, \beta_1, \ldots, \beta_l \right)$$
Assume now that $ k(x) = k(y) $, namely--
$$ j_W(f)\left( \kappa,  \theta_{ j_W\left( f_{\beta_1} \right)(\kappa)  },   \ldots, \theta_{  j_W\left( f_{\beta_l} \right)(\kappa) }  \right)    \in  j_W(g)\left( \kappa,  \theta_{  j_W\left( f_{\beta_1} \right)(\kappa)  },   \ldots, \theta_{  j_W\left( f_{\beta_l} \right)(\kappa)  } \right)    $$
Then--
$$ \{ \xi<\kappa \colon   f\left( \xi,  \theta_{  f_{ \beta_1 }(\xi) }  , \ldots,  \theta_{  f_{ \beta_l }(\xi) } \right) \in g\left( \xi,  \theta_{  f_{ \beta_1 }(\xi) }  , \ldots,  \theta_{  f_{ \beta_l }(\xi) } \right) \}\in W $$
and by the definition of $ W $, there exists $ p\in G $ and $ r\in H $ such that-- 
$$ p^{\frown}r \Vdash \kappa \in i\Big( \{ \xi<\kappa \colon  f\left( \xi,   \lusim{\theta}_{  f_{ \beta_1 }(\xi) }  , \ldots,  \lusim{\theta}_{  f_{ \beta_l }(\xi) } \right) \in g\left(\xi,  \lusim{\theta}_{  f_{ \beta_1 }(\xi) }  , \ldots,  \lusim{\theta}_{  f_{ \beta_l }(\xi) } \right)  \} \Big) $$
By extending $ r\in H$ finitely many times, $ p^{\frown}r \Vdash \lusim{\theta}_{\left( i(f_{\beta_m } )\left( \kappa \right) \right)} =  \beta_m  $ holds for every $ 1\leq m\leq k $. Thus, the last equation can be replaced with--
$$ p^{\frown} r \Vdash i(f)\left( \kappa, \beta_1, \ldots, \beta_l \right) \in i(g)\left( \kappa, \beta_1 , \ldots, \beta_l\right) $$
but the forced statement above is entirely in $ N $, and since a condition forces it, it is true in $ N $. Thus--
$$ i(f)\left( \kappa, \beta_1, \ldots, \beta_l \right)\in i(g)\left( \kappa, \beta_1, \ldots ,\beta_l \right) $$
as desired. The implication in the other direction is proved similarly.

Clearly $ \mbox{crit}(k)> \kappa $. We finish the proof by showing that $ j_W\restriction_{V} = k\circ i $. Let $ x\in V $ and let $ c_x \colon \kappa \to V $ be the constant function with value $ x $. Then--
$$ k\left(i(x)\right) = k\left( i\left( c_x \right)(\kappa) \right) = j_W\left( c_x \right)(\kappa)= j_W(x) $$
as desired.
$\square$

Let us now study the properties of the embedding $ k\colon N\to M $. We assume the settings of theorem \ref{Theorem: Conditions imposed on i}.

\begin{lemma}

If $\leq=\leq^*$, or at least $\leq_\alpha=\leq_\alpha^*$, for a final segment of $\alpha<\kappa$, then $k$ is the identity and $M=N$.

\end{lemma}
\pr
Fix an ordinal $ \eta $, and let $ f\in V\left[G\right] $  be a function such that $ \eta = \left[f\right]_W $. We will prove that $ \eta\in \mbox{Im}(k) $. Indeed, consider the set--
$$\{p\in i(P)/G \mid \exists \tau (p\Vdash i(\lusim{f})(\kappa)=\tau)\}$$
It is $\leq-$dense in $N[G]$.
So, if $\leq=\leq^*$, then $H$ meets it. Thus, there exists a condition $ q\in H $, a function $ g\in V $ and generators $ \beta_1, \ldots, \beta_l $ of $ i $, such that $ q\Vdash i\left( \lusim{f} \right)(\kappa) = i(g)\left( \kappa, \beta_1, \ldots, \beta_l \right) $. Thus, by the definition of $ W $,
$$  \{  \xi<\kappa \colon f(\xi) = g\left( \xi,  \theta_{  f_{ \beta_1 }(\xi)  }, \ldots, \theta_{  f_{ \beta_l }(\xi) }  \right) \} \in W $$
and thus $ \eta = \left[  f \right]_W = k\left(   i(g)\left( \kappa, \beta_1, \ldots, \beta_l \right) \right) $. 
\\
$\square$

In general, $M$ should not be equal to $N$. Thus, for example, they will differ if the Prikry forcing was used unboundedly often below $\kappa$. 

However, we do not know whether the assumption of 2.15 is necessary.

\begin{question}
Suppose that for unboundedly many $\alpha<\kappa$,  $\leq_\alpha\not=\leq_\alpha^*$. Is then $M\neq N$?	
\end{question}

We do not know what are the requirements on the forcings $ \lusim{Q}_{\alpha} $ for $\alpha\in \Delta $ which imply $ M=N $. We conjecture that the requirement should be that there is $\delta<\kappa$, such that every set of ordinals $x$ of $V[G]$ can be covered by a set $y \in V$ of cardinality $\leq  |x|+\delta$.

\begin{lemma}\label{lem3-1}

$k''H\subseteq G^*\setminus \kappa$.

\end{lemma}
\pr
Let $q $ be in $H$, and let $ p\in G $ be a condition such that $ p\Vdash \lusim{q}\in \lusim{H} $ (recall that $ H\in V\left[G\right] $). Clearly,
$$p^{\frown} \lusim{q} \Vdash \lusim{q} \in \Gamma\setminus \kappa$$
where $ \Gamma $ is the canonical $ i(P) $-name for the generic set for $ i(P) $ over $ V $.

Pick $f:[\kappa]^n\to \kappa, f\in V$ and $\beta_1<...<\beta_n<i(\kappa)$ such that
$\lusim{q}=i(f)(\beta_1,...,\beta_n)$.
For every $m,1\leq m\leq n$, there are  $f_m:\kappa\to \kappa,f_m\in V$ such that $ q_{ i\left( f_m \right)(\kappa), \beta_m }\in H $, namely, $\beta_m =\theta_{i(f_m)(\kappa)}$.\\
\\Let us argue that the set--
$$A_q=\{\nu<\kappa\mid f(\theta_{f_1(\nu)},...,\theta_{f_n(\nu)})\in \lusim{G}\setminus \nu \}$$ is in W. Pick any $q\leq^* q^*\in H$ which $ \leq^* $ which forces that $ \beta_{m} = \theta_{ i\left( f_{ \beta_m } \right)(\kappa) } $, for every  $ 1\leq m\leq n $. Recall that--
$$q=i(f)(\beta_1,...,\beta_n)=i(f)(\theta_{i(f_1)(\kappa)},...,\theta_{i(f_n)(\kappa)})$$
and thus $p^{\frown}\lusim{q}^* \Vdash \kappa\in i(\lusim{A}_q)$.
\\
$\square$

The next lemma generalizes a Kunen-Paris result (see remark \ref{Remark: Kunen Paris}).

\begin{theorem}\label{Theorem: Kunen-Paris extenstion}
Let $ H,H'\in V\left[G\right] $ be generic sets for $ \l i(P)\setminus \kappa, \leq^* \r $ over $ N\left[G\right] $.
Suppose that $H$ and $H'$ satisfy  $ (*) $. 
Assume that for every $\beta<\kappa$, if $q,q' \in Q_\beta$ are incompatible according to the order $\leq^*$, then--
$$D_{\beta}(q) = \{r\in Q_\beta\mid  r \text{ is } \leq-\text{incompatible with } q\}$$ is $\leq^*-$dense   above $q'$, or
$$D_{\beta}(q') = \{r\in Q_\beta\mid  r \text{ is } \leq-\text{incompatible with } q'\}$$ is $\leq^*-$dense   above $q$.\footnote{ This type of condition usually holds. For example, if we iterate Prikry forcings, then just shrinking sets of measure one will produce such type of incomparability.}
\\Suppose that $H\not =H'$, then $U_H\not = U_{H'}$.
\end{theorem}

\begin{remark}
Note that if the $ Q_{\beta} $-s are taken to be Prikry forcings, then the above property holds. Indeed, assume that $ q = \langle t,A \rangle  $ and $ q' = \langle t', A' \rangle $ are $ \leq^* $-incompatible. Then $ t\neq t' $. Assume without loss of generality that $ t $ is an end extension of $ t' $. Then $ D(q) = \{ r \colon r,q \mbox{ are }\leq^* \mbox{-incompatible}  \} $ is $ \leq^* $-dense open above $ q' $. Indeed, pick a condition $ \langle t', B \rangle\geq^* \langle t', A \rangle $. Shrink $ B $ to the set $ B^* = B\setminus \left( \max(t)+1 \right) $. Then $ \langle t', B^* \rangle\geq^* \l t', B \r $ and is incompatible with $q = \langle t, A \rangle $.
\end{remark}

\pr   Suppose otherwise, i.e. $H\not =H'$, but $U_H = U_{H'}:=W$.
\\Let $k:N \to M$ be the elementary embedding defined from $H$ and $k':N \to M$ from $H'$.

\begin{claim} \label{Claim: Technical First claim in Kunen Paris}
$k\neq k'$.
\end{claim}

\pr
Assume for contradiction that $ k= k' $. Thus, by Lemma \ref{lem3-1}, every pair of elements from $ H, H' $ are $ \leq$-compatible. We will argue that this implies that $ H = H' $. It suffices to prove that every pair of conditions $ q\in H , q'\in H'$ are $ \leq^* $-compatible.

Assume otherwise. Let $ \alpha<\kappa $ be the least ordinal such that
there are pair of conditions $ q\in H, q'\in H' $ for which $ q\restriction_{\alpha} , q'\restriction_{\alpha}$ are $ \leq^* $-incompatible. $ \alpha $ cannot be limit, since $ \leq^* $-compatibility of all the initial segments of $ q,q' $ below $ \alpha $ implies that $ q\restriction_{\alpha} $ and $ q'\restriction_{\alpha} $ are $ \leq^* $-compatible themselves (if $ \alpha $ is inaccessible, this is clear since the support of $ q,q' $ is bounded in $ \alpha $; if the supports of $ q,q' $ are unbounded in $ \alpha $, just intersect sets of measure one to find a common direct extension). Thus $ \alpha = \beta+1 $ is successor, and $ q\left( \beta \right) $, $ q'\left( \beta \right) $ are $ \leq^* $-incompatible. By the property of the forcing $ Q_\beta $, without loss of generality, $ D_\beta\left( q \right) $ is $ \leq^* $-dense open above $ q' $. Since $ q'\in H'\left( \beta \right) $, $ q' $ can be extended to a condition $r\in  H' $, such that $ r(\beta)\in D_{\beta}(q) $. In particular, $ q\in H $, $ r\in H' $ are $ \leq $-incompatible, which is a contradiction. $\square \mbox{ of claim}$ \ref{Claim: Technical First claim in Kunen Paris}.\\

Since $ k\neq k' $, there exists a generator $ \beta $ of $ i $ such that $ k(\beta) \neq k'(\beta) $. Pick the least such generator $ \beta $. 

\begin{claim} \label{Claim: Technical claim in Kunen Paris}
For every generator $ \beta'< \beta $ of $ i $, there exists a function $ f_{\beta'}\in V $ such that each generic $ H,H' $ has a condition which forces that $ \beta' =  \lusim{\theta}_{i\left(f_{\beta'}\right)(\kappa)}  $.
\end{claim}

\pr Let $ f,f' $ be functions such that some condition in $ H $ forces that $ \beta = \lusim{\theta}_{i(f)(\kappa)} $, and some condition in $ H' $ forces that $ \beta' = \lusim{\theta}_{ i(f')(\kappa) } $. Let $ q\in H $ be a condition which decides the statement $ \lusim{\theta}_{ i(f)(\kappa) }  =  \lusim{\theta}_{ i(f')(\kappa) }  $ and assume for contradiction that it is decided negatively. By applying $ k $, $ k(q)\in j_W(G) $ forces that--
$$ \lusim{\theta}_{ \left[ f \right]_W }  \neq \lusim{\theta}_{ \left[f'\right]_W } $$
namely $ k(\beta') \neq k'\left( \beta' \right) $, contradicting the minimality of $ \beta $. $\square \mbox{ of claim}$ \ref{Claim: Technical claim in Kunen Paris}.

Recall now that $ k(\beta) \neq k'(\beta) $. Thus, there are two distinct functions $ f,f' $ in $ V $ such that--
\begin{enumerate}
	\item Some condition in $ H $ forces that $ \beta = \lusim{\theta}_{i(f)(\kappa)} $.
	\item Some condition in $ H' $ forces that $ \beta = \lusim{\theta}_{i(f')(\kappa)} $.
	\item Without loss of generality, $ \{\xi< \kappa \colon \theta_{f'(\xi)} < \theta_{f(\xi)} \}\in W $.
\end{enumerate}

By property $ (2) $ of the names $ \langle \lusim{\theta}_{\alpha} \colon \alpha<\kappa \rangle $, presented in lemma \ref{Existence of the sequence  theta alpha}, there exists an ordinal $ \beta' $ such that some condition in $ H $ forces that $ \lusim{\theta}_{ i\left( f' \right)(\kappa) } = \beta' $. By the above assumptions, $ \beta'< \beta $.

We argue that $ \beta' $ is a generator of $ i $ as well. This will finish the proof: once we prove that $ \beta' $ is a generator of $ i $, it follows from claim \ref{Claim: Technical claim in Kunen Paris} that $ \theta_{ i(f')(\kappa) } $ represents $ \beta' $ in the sense of both generics, $ H,H' $. However, in the sense of $ H' $, it represents $ \beta $, which is a contradiction.

Assume for contradiction that $ \beta' $ is not a generator of  $ i $. Then there is a function $ g\in V $ and $ \beta_1, \ldots, \beta_l $ below $ \beta' $, such that $\beta' = i(h)\left(  \kappa, \beta_1, \ldots, \beta_l \right)  $. Since $ H $ forces that $ \beta' = \lusim{\theta}_{ i\left( f' \right)(\kappa) } $, it follows that--
$$ \{ \xi <\kappa \colon  g\left( \xi, \theta_{f_{\beta_1}(\xi) }, \ldots, \theta_{  f_{ \beta_l }(\xi) } \right)  = \theta_{  f'(\xi) }  \}\in U_H = W $$
Thus the same set belongs to $ U_{H'} $. Therefore, $ H' $ forces that--
$$ \beta = \theta_{ i(f')(\kappa) } = i(g)\left( \kappa, \beta_1, \ldots, \beta_l \right) $$
contradicting the fact that $ \beta $ is a generator of $ i $ (note that we used claim \ref{Claim: Technical claim in Kunen Paris} when arguing that the generators $ \beta_{i} $, $ 1\leq i \leq l $, are represented the same way in the sense of $ H, H' $). $\square$

\begin{definition}
A measure $ W\in V\left[G\right] $ is called simply generated if $ W = U_H $ for some $ U\in V $, where $ H $ is generic for $ \langle j_U(P)\setminus \kappa, \leq^* \rangle $ over $ M_U\left[G\right] $. 
\end{definition}

\begin{remark}
Given a simply generated normal measure $ W\in V\left[G\right] $ as above, with $ \Delta\notin W $, the parameters $ U $ and $ H $ are uniquely defined from it\footnote{Note that when iterating Prikry forcings, $ \Delta\notin W $ holds for every normal measure $ W\in V\left[G\right] $ on $ \kappa $,  since every such $ W $ concentrates on regulars.}. Indeed, we will prove in the next lemma that $ U = W\cap V $ belongs to $ V $, and is a normal measure there with $ \Delta\notin U $,  Now, assume that there are $ H, H' $, generic over $ M_U\left[G\right] $ for $ \langle j_U(P)\setminus \kappa , \leq^* \rangle $, with $ W = U_{H} = U_{H'} $.  Then $ H,H' $ satisfy the conditions of lemma \ref{Theorem: Kunen-Paris extenstion} (since $ j_U $ has no generators other than $ \kappa $). Thus, by the theorem, $ H = H' $.
\end{remark}

Given $ W\in V\left[G\right] $ normal  on $ \kappa $ (which is not necessarily simply generated),  we can say the following:

\begin{lemma} \label{Lemma: Every ultrafilter extends one from V}
Every normal measure $ W\in V\left[G\right] $ on $ \kappa $ extends a measure $ U = W\cap V \in V $.  
\end{lemma}

\pr
First, let us argue that $ U = W\cap V $ belongs to $ V $. By \cite{RestElm}, it suffices to prove that there are no new fresh unbounded subsets of cardinals in the interval $ \left[\kappa, \left(2^{\kappa}\right)^V\right] = \left[ \kappa, \kappa^{+} \right] $. Thus, it suffices to prove the following pair of claims:

\begin{claim}
	$ P = P_{\kappa} $ does not add fresh unbounded subsets to $ \kappa $.
\end{claim}

\pr
The fact that there are no fresh unbounded subsets of $\kappa$ follows essentially from the facts that $ 2^{\kappa} = \kappa^{+} $, and that there exists a normal measure on $ \kappa $ in $ V\left[G\right] $: Given a normal measure $ U\in V $ with $ \Delta \notin U $, take any $ U^*\in V\left[G\right] $ which extends it. Given a fresh unbounded $ A\subseteq \kappa $, $ A = j_{U^*}\left( A \right)\cap \kappa $ and thus, by elementarity, $ A $ belongs to the ground model $ M $ of $ \mbox{Ult}\left( V\left[G\right], U^* \right) $. Now set $ k_U\colon M_{U}\to M $ to be the function which maps $ \left[f\right]_{U} $ to $ \left[f\right]_{U^*} $. Then $ k_U $ is a well defined elementary embedding since $ U\subseteq U^* $, and $ \mbox{crit}\left( k_U \right) > \kappa $ by normality of $ U^* $. Since $ 2^{\kappa} = \kappa^{+} $ holds in $ M $, $ k_U $ maps the sequence of subsets of $ \kappa $ to itself, and thus every subset of $ \kappa $ which belongs to $ M $, already belongs to $ V $. So the above set $ A $ belongs to $ V $, which is a contradiction.
$\square$

\begin{claim}
	For every measurable (in $ V $) $ \lambda \leq \kappa $, $ P_{\lambda} $ doesn't add fresh unbounded subsets of $ \lambda^{+} $. In particular, $ P_{\kappa} $ does not add fresh subsets to $ \lambda^{+} $.
\end{claim}

\pr
Let $ f\in V\left[G\right] $ be the characteristic function of a fresh unbounded subset of $ \lambda^{+} $. Let $ \lusim{f} $ be a $ P_{\lambda} $-name and assume that $ p\in P $ forces that $ \lusim{f} $ is fresh.

Let $ G\subseteq P_{\lambda} $ be generic over $ V $. For every $ \xi<\lambda^{+} $, let $ p_{\xi} \in G $ be a condition which decides $ \lusim{f}\restriction_{ \xi} $.  For every $ \xi<\lambda^{+} $ there exists $ \alpha_{\xi} <\lambda$ such that the support of $ p_{\xi} $ is bounded by $ \alpha_{\xi} $. Let $ A\subseteq \lambda^{+} $ and $ \alpha^*<\lambda $ be such that $ \left|A\right|=\lambda^{+} $  and  $ \alpha_{\xi} =\alpha^* $ for every $ \xi\in A $. 

By shrinking $ A\subseteq \lambda^{+} $ even further (to a set of cardinality $ \lambda^{+} $), we can assume that there exists $ q^* \in P_{\lambda} $ such that, for every $ \xi\in A $, $ p_{\xi} \restriction_{ \alpha^*} = {q^*}\restriction_{ \alpha^*} $, and $ q^* \restriction_{\left[ \alpha^* , \lambda \right)} $ is trivial. 

Let $ h = \bigcup\{ g\colon  \exists \xi<\lambda  \ q^* \Vdash \lusim{f}\restriction_{\xi} = g  \} $. Clearly, $ h\colon \lambda^{+} \to 2 $ is a function and $ q^* \Vdash \lusim{f} = \check{h} $.
$\square$ \\
$\square$ \mbox{ of lemma } \ref{Lemma: Every ultrafilter extends one from V}.

\section{On $j_W(\kappa)>j_U(\kappa)$ and the existence of $ N $} \label{Section: on jwkappa above jukappa}

Assume in this section that $ P = P_{\kappa} $ is an iteration of Prikry forcings. Let $ W,U = W\cap V, i\colon V\to N $ be as in the previous section. 

Clearly, $j_W(\kappa)\geq j_U(\kappa)$. Our interest here will be in situations where a strict inequality holds. 

Note such phenomenon is impossible with the non-stationary support, where, for every normal measure $ W\in V\left[G\right] $ on $ \kappa $, $ j_W(\kappa) = j_U(\kappa) $ (see \cite{NonStatRestElm}). 

On the other hand, in the full support iteration, it is possible that $ j_W(\kappa)> j_U(\kappa) $ starting with $ o(\kappa) \geq 2 $. Indeed, under the assumption that there are normal measures $ U_0 \vartriangleleft U_1 $ on $ \kappa $ in $ V $, take $ W = \left( U_1 \right)^{\times} $ (in the notations of \cite{kaplan2022magidor}). Assume that $ \xi \mapsto U_0(\xi) $ is a function in $ V $ which represents $ U_0 $ in $ \mbox{Ult}\left( V, U_1 \right) $, and, for every $ \xi\in \Delta $, $ U_0(\xi) $ is the normal measure used to singularize $ \xi $ at stage $ \xi $ in the iteration $ P = P_{\kappa} $.  Then $ W $ extends $ U_0 $, but $ j_W(\kappa) = j_{U_1}(\kappa) > j_{U_{0}}(\kappa) $. The other direction is also true: if $ j_W(\kappa) > j_{U}(\kappa) $, let $ U' $ be a measure on $ \kappa $ in $ V $ such that $ W = \left( U' \right)^{\times} $. Then $ U\vartriangleleft U' $ and thus $ o(\kappa) \geq 2 $.

Let us discuss this situation in the context of the Easton support iteration.

\subsection{On $j_W(\kappa)>j_U(\kappa)$} \label{Siubsection of jw above ju}

Start with the following simple observation:

\begin{proposition} \label{proposition: set unbounded in jWkappa}
	The set--
	$$\{j_W(f)(\kappa)\mid f:\kappa\to \kappa, f\in V\}$$
	is unbounded in $j_W(\kappa)$.
	\\Hence, $k''j_U(\kappa)$ is unbounded in $j_W(\kappa)$, where $k([f]_U)=[f]_W$ is the embedding defined in lemma \ref{Lemma: k is elementary}.
	
\end{proposition}
\pr
$P$ satisfies $\kappa-$c.c. Hence for every $g:\kappa\to \kappa$ in $V[G]$ there is $f:\kappa\to \kappa$ in $V$ which dominates it, i.e.,
for every $\nu<\kappa$, $g(\nu)<f(\nu)$.
\\
$\square$

Let us present a first example of a situation where $j_W(\kappa)>j_U(\kappa)$.

\begin{definition}(W. Mitchell)
	A cardinal $\kappa$ is called \emph{$\mu-$measurable} iff there exists an extender $E$ over $\kappa$ such that $E_\kappa\in M_E$,
	where $E_\kappa=\{A\subseteq \kappa\mid \kappa\in j_E(A)\}$.

\end{definition}

Note that we can use a witnessing extender $E$ with two generators only - $\kappa$ and the ordinal $\eta<2^{2^\kappa}$
which codes $E_\kappa$.
The ultrapower by such extender is closed under $\kappa-$sequences.

The next lemma is obvious:

\begin{lemma}
	Suppose that $\kappa$ is $\mu-$measurable and $E$ is an extender witnessing this.
	Then $j_{E_\kappa}(\kappa)<j_E(\kappa)$.

\end{lemma}

\begin{proposition}\label{prop1}
	Suppose that $\kappa$ is $\mu-$measurable and $E$ is an extender witnessing this which ultrapower is closed under $\kappa-$sequences.
	Let $U=E_\kappa$ and $\Delta\subseteq \kappa$ be a set of measurable cardinals which is not in $U$.
	Force with an Easton support iteration $P$ of the Prikry forcings over $\Delta$.
	Let $G\subseteq P$ be a generic.
	\\Then, in $V[G]$, there is a normal ultrafilter $W$ which extends $U$ such that  $j_W(\kappa)>j_U(\kappa)$.
	
\end{proposition}
\pr
Construct $W$ as in theorem \ref{Theorem: Conditions imposed on i}  using $E$, i.e., $i=j_E$ and $N=M_E$.
\\Then $j_U(\kappa)<j_E(\kappa)=i(\kappa)=j_W(\kappa)$.
\\
$\square$

Let us observe now that we need a $\mu-$measurable in order to have $j_W(\kappa)>j_U(\kappa)$, provided $V=\calK$, where $\calK$ denotes the core model.

\begin{proposition}\label{prop2}
	Assume $\neg 0^{\P}$.
	Suppose that $V=\calK$.
	Let $U$ be a normal ultrafilter over $\kappa$ and $\Delta\subseteq \kappa$ be a set of measurable cardinals which is not in $U$.
	Force with an Easton support iteration $P$ of the Prikry forcings over $\Delta$.
	Let $G\subseteq P$ be a generic.
	\\Suppose that, in $V[G]$, there is a normal ultrafilter $W$ which extends $U$ such that  $j_W(\kappa)>j_U(\kappa)$.
	\\Then $\kappa$ is a $\mu-$measurable in $V$. Moreover, $U$ is a normal measure of a witnessing extender.

\end{proposition}
\pr
Suppose otherwise.
\\Consider $j_W\upr V$.
\\By Mitchell \cite{mitchell2009beginning}, it is a normal\footnote{ Extenders with smaller indexes are used first.} iterated ultrapower of $\calK=V$ by its measures and extenders.
Recall that $W\cap V=U$, and so,
$U=\{A\subseteq \kappa \mid A\in V, \kappa\in j_W(A)\}$.
The assumption that   $U$  is not a normal measure of an extender which witnesses a $\mu-$measurability of $\kappa$ implies then that
$U$
must be used first in this iterated ultrapower.

Apply now the arguments of \cite{gitik1993measurable, G-M1996indiscernible} in $\calK_U$ the core model of $M_U$.
For every measurable $\alpha, \kappa\leq \alpha<j_U(\kappa)$, there will be a bound $\eta_\alpha$ (which depends on $o(\alpha)$) on number
of possible applications of measures and extenders over $\alpha$ with their images, and, by the assumption that there are no strong cardinals,
$\eta_\alpha<j_U(\kappa)$. Therefore, for every such $ \alpha $, there exists an upper bound $ \mu^*_{\alpha}< j_U(\kappa) $ on the image of $ \alpha $ in the iterated ultrapower by the measures or extenders taken on $ \alpha $ or its images. $ \mu^*_{\alpha} < j_U(\kappa) $ since none of the extenders participating has length $ j_U(\kappa) $ or above (by $\neg 0^{\P}$), and $ \eta_{\alpha}< j_U(\kappa) $. 

This implies that the rest of the iteration from $ j_U $ to $ j_W\restriction_{\mathcal{K}} $ cannot move $ j_U(\kappa) $: otherwise, $ j_U(\kappa) $ participates as a critical point in the iteration (it cannot be moved by an extender with a critical point below $ j_U(\kappa) $, as explained above; it surely cannot moved by an extender with a critical point above $ j_U(\kappa) $; thus, it moves since an extender on it participates in the iteration). Decompose $ j_W\restriction_{V} = k_1\circ i_1$, where $ i_1 $ is the iteration with all the extenders below $ j_U(\kappa) $, and $ k_1 $ is an iteration with critical point $ j_U(\kappa) $. Then--
$$ \{ j_W(f)(\kappa) \colon f\colon \kappa\to \kappa, \ f\in V \}  $$
is bounded in $ j_W(\kappa) $, since, for every $ f\colon \kappa\to \kappa $ in $ V $, 
$$ j_W(f)(\kappa) = k_1\circ \left( i_1\left( f \right)(\kappa) \right) = i_1\left( f \right)(\kappa)< j_U(\kappa) $$
where we used the fact that $  i_1\left( f \right)(\kappa) $ does not move in $ k_1 $ since it is strictly below $ j_U(\kappa) $. This, contradicts proposition \ref{proposition: set unbounded in jWkappa}.
%
%
$\square$\\

The situation changes if we do not assume $V=\calK$.
Let us argue now that the consistency strength of $j_W(\kappa)>j_U(\kappa)$ is just a measurable which is a limit of measurable cardinals.

\begin{proposition}\label{prop3}
	Let $V_0$ be a model of GCH with a measurable cardinal $\kappa$ which is a limit of measurable cardinals.
	\\
	Then there is a cardinal preserving generic extension $V$ of $V_0$ which satisfies the following:
	\\
	Let $\Delta$ be an unbounded subset of $ \kappa$ consisting of measurable cardinals.
	Force with an Easton support iteration $P$ of the Prikry forcings over $\Delta$.
	Let $G\subseteq P$ be a generic.
	\\There  exists  a normal ultrafilter $U$ over $\kappa$ in $V$ and
	a normal ultrafilter $W$ in $V[G]$ which extends $U$ such that  $j_W(\kappa)>j_U(\kappa)$.

\end{proposition}
\pr
The idea is as follows. Let $W$ be a normal ultrafilter over $\kappa$ in $V_0$ which concentrates on non-measurable cardinals.
Consider $W^2=W\times W$ and $W^3=W\times W\times W$.
\\Let $j_1=j_W, j_2=j_{W^2}, j_3=j_{W^3}, M_1=M_W, M_2=M_{W^2},M_3=M_{W^3},\kappa_1=j(\kappa), \kappa_2=j_2(\kappa),\kappa_3=j_3(\kappa)$.
We have natural commuting embeddings  $j_{12}:M_1\to M_2, j_{23}:M_2\to M_3$ and $j_{13}:M_1 \to M_3$.
Namely, $j_{12}(j_1(f)(\kappa))=j_2(f)(\kappa), j_{23}(j_2(g)(\kappa,\kappa_1))=j_3(g)(\kappa,\kappa_1)$, etc.
Note that the critical point of $j_{12}, j_{13}$ is $\kappa_1$ and of $j_{23}$ is $\kappa_2$.
However there is an additional way to embed $M_2 $ into $M_3$. Define $\sigma:M_2 \to M_3$ by setting
$\sigma(j_1(f)(\kappa,\kappa_1))=j_3(f)(\kappa,\kappa_2)$. Clearly, $\sigma$ is elementary and its critical point is $\kappa_1$ and it is moved to $\kappa_2$.
\\The idea will be to force in order to extend $W$ to a normal ultrafilter $U$ such that
\begin{enumerate}
	\item $M_U$ is a generic extension of $M_2$,
	\item $W^3$ extends to a $\kappa-$complete ultrafilter $E$ with $M_E$ a generic extension of $M_3$,
	\item $U$ is the normal ultrafilter which is strictly below $E$ with the corresponding embedding\\  extending $\sigma$.
\end{enumerate}

Now, $\kappa_2<\kappa_3$ will imply $j_U(\kappa)<j_E(\kappa)$, since  $j_U(\kappa)=\kappa_2$ and $j_E(\kappa)=\kappa_3$.

Such construction was used in \cite{gitik1981non}. We refer to this paper for details. Let us only sketch the argument.

We force a Cohen function $f_\alpha:\alpha\to \alpha$ for every inaccessible $\alpha\leq \kappa$ using the iteration with an Easton support.
\\Denote a generic object which produces such $\l f_\alpha\mid \alpha\leq \kappa, \alpha \text{ is an inaccessible }\r$ by $G_0$.
\\Let $V=V_0[G_0]$.\\
It is possible to extend all the embeddings, $j_1, j_2, j_3, j_{12}, j_{13}, j_{23}, \sigma$.
We change one value of $f_{\kappa_3}$ at $\kappa$ by setting it to $\kappa_2$. Let $G_3$ be such generic over $M_3$
Then, $j_3:V_0\to M_3$ extends to $j_3^*:V_0[G_0]\to M_3[G_3]$.
Derive now $U$ and $E$ from $j_3^*$ ,in $V=V_0[G_0]$, by setting
$U=\{A\subseteq \kappa \mid \kappa\in j_3^*(A)\}$ and
$E=\{B\subseteq \kappa^3 \mid \l\kappa,\kappa_1,\kappa_2\r\in j_3^*(B)\}$.

Finally we apply the construction of Section \ref{Section: The general Framework} to $U$ and $E$ to produce an extension $W$ of $U$ in $V[G]$.
\\
$\square$

Note that $U$ produced in \ref{prop3} can be picked to be the minimal in the Mitchell order, which is not true about one of \ref{prop1}, where $V=\calK$.
Let us argue that under rather strong assumptions it is possible to find such $U$ in $\calK$.

\begin{proposition}\label{prop4}
	Let $U$ be a normal ultrafilter over $\kappa$.
	Suppose that the set
	$$\{\alpha<\kappa \mid \alpha \text{ is } \kappa-\text{strong }\}$$ is unbounded in $\kappa$.
	Force with $P$ as above. Let $G\subseteq P$ be a generic.
	Then, in $V[G]$, there is a normal ultrafilter $W$ over $\kappa$ such that
	
	\begin{enumerate}
		\item $U\subseteq W$,
		\item $j_U(\kappa)<j_W(\kappa)$,
		\\Moreover, $j_W\upr V=k\circ i$, where
		\begin{itemize}
			\item $i:V \to N$,
			\item $j_U(\kappa)<i(\kappa)$,
			\item $i, N$ satisfy the conditions of theorem \ref{Theorem: Conditions imposed on i}.
		\end{itemize}
		
	\end{enumerate}

\end{proposition} \label{Proposition: kappa limit of lappa strings suffices for jW above jU}
\pr
Work in $M_U$. Pick some $\alpha, \kappa<\alpha<j_U(\kappa)$ which is $j_U(\kappa)-$strong.
Let $E\in M_U$ be an $(\alpha, j_U(\kappa))-$extender witnessing this.
Set $N$ to be the ultrapower of $M_U$ by $E$ and let $i=j_E\circ j_U$.
We have $$j_U(\kappa)\leq j_E(\alpha)<j_E(j_U(\kappa))=i(\kappa).$$
Note that the embedding $ i $ satisfies the assumptions of theorem \ref{Theorem: Conditions imposed on i}. Indeed, the only nontrivial properties of $ i $ that require verification are:
\begin{enumerate}
	\item $ \{ i(f)(\kappa) \colon f\colon \kappa \to \kappa \} $ is unbounded in $ i(\kappa) $: In $ M_{U} $, denote $ \lambda = j_U(\kappa) $. Then $ \lambda $ is regular, and thus the $ \left( \alpha, \lambda \right) $-extender embedding $ j_E $ is continuous at $ \lambda $. Thus, for every $ \beta< i(\kappa) = j_E\left( \lambda \right) $, there exists $ \beta'< \lambda $ such that $ j_E(\beta')>\beta $. Now, find $ f\in V $, $ f\colon \kappa \to \kappa $, such that $\beta' = j_U(f)(\kappa) $. Then $ i(f)(\kappa) = j_E(\beta') > \beta $, as desired.
	$ \left| j_E(\kappa) \right| = \left(j_{U}(\kappa)\right)^{+} $. 
	\item $ \left| i(\kappa) \right| = \kappa^{+} $: Using the above notations, $ M_U\Vdash \left| j_{E}(\lambda) \right| = \lambda^{+} $. But $ V\Vdash \left| \lambda^{+} \right| = \kappa^{+} $  since $ 2^{\kappa} = \kappa^{+} $.
\end{enumerate}
Now apply theorem \ref{Theorem: Conditions imposed on i} to construct the desired measure $ W $.
$\square$

We do not know whether the assumption of \ref{prop4} is really necessary.
However it is possible to show the following.

\begin{proposition}\label{prop5}
	Suppose  $\neg 0^{\P}$.
	\\Assume $V=\calK$.
	\\Let $U$ be a normal ultrafilter over $\kappa$ which is minimal in the Mitchell order.
	\\Let $P$ be an Easton support iteration of Prikry type forcing notions up to $\kappa$ and $G\subseteq P$ be a generic.
	\\Suppose that $W$ is a normal ultrafilter in $V[G]$ which extends $U$.
	\\Then $j_U(\kappa)=j_W(\kappa)$.

\end{proposition}
\pr
By W. Mitchell \cite{mitchell2009beginning},  $j_W\upr \calK$ is a normal iterated ultrapower of $\calK$ by its measures and extenders.
\\The minimality of $U$ implies that it must be used first in this iteration.
\\
Apply now the arguments of \cite{gitik1993measurable,G-M1996indiscernible} in $\calK_U$ the core model of $M_U$.
For every measurable $\alpha, \kappa\leq \alpha<j_U(\kappa)$, there will be a bound $\eta_\alpha$ (which depends on $o(\alpha)$) on number
of possible applications of measures and extenders over $\alpha$ with their images, and, by the assumption that there is no strong cardinals,
$\eta_\alpha<j_U(\kappa)$. Now complete the argument as in proposition \ref{prop2} by showing that, if $ j_U(\kappa) < j_W(\kappa)$ then an extender with critical point $ j_U(\kappa) $ participates in the iteration, and thus the set $ \{ j_W(f)(\kappa) \colon f\colon \kappa\to \kappa, \ f\in V \} $ is bounded by $ j_U(\kappa) $; this contradicts proposition \ref{proposition: set unbounded in jWkappa}.
$\square$\\

We conjecture that the needed strength (for \ref{prop5}) is exactly $$\{\alpha<\kappa \mid \alpha \text{ is } \kappa-\text{strong }\} \text{ is unbounded in } \kappa.$$
Thus, R. Schindler \cite{schindler2006iterates} extension of the Mitchell result can be used to argue that  $j_W\upr \calK$ is a normal iterated ultrapower of $\calK$ by its measures and extenders. A missing part is an extension of \cite{gitik1993measurable}
beyond strongs which is likely to hold.

\subsection{On existence of $N$} \label{Subsection: on existence of N}

As before, let $ W\in V\left[G\right] $ be a normal measure on $ \kappa $, and $ U = W\cap V\in V $. In section \ref{Section: properties of k} we will prove that if $ W $ is constructed as in theorem \ref{Theorem: Conditions imposed on i} then $ j_W\restriction_{V} = k\circ i $, where $ k $ is an iteration of $ N $ be normal measures only. A natural question in view of this result is whether for every $ W\in V\left[G\right] $  there  exists $N, {}^\kappa N\subseteq N$
such that $M$ is obtained from it by iterating normal measures only.
We do not know the answer in general. However, it turns out to be an affirmative provided some anti large cardinal assumptions and $V=\calK$.

\begin{proposition}
	Assume $\neg 0^{\P}$ and $V=\calK$.
	\\Let $U$ be a normal ultrafilter over $\kappa$ and $\Delta\subseteq \kappa$ be a set of measurable cardinals which is not in $U$.
	Force with an Easton support iteration $P$ of the Prikry forcings over $\Delta$.
	Let $G\subseteq P$ be a generic.
	\\Suppose that, in $V[G]$, there is a normal ultrafilter $W$ which extends $U$.
	\\Then there are $N, i:V\to N$ which satisfy the conditions of theorem \ref{Theorem: Conditions imposed on i} such that $j_W\upr V=k\circ i$ and $k$ is formed by iterating normal measures only,
	starting from $N$.

\end{proposition}
\pr
As in subsection \ref{Siubsection of jw above ju}, we analyze $j:=j_W\upr \calK$.
\\ By elementarity, $j:\calK \to (\calK)^{M_W}$ and $M_W$ is a generic extension of $(\calK)^{M_W}$ by an Easton support iteration of Prikry forcings
with normal measures in $j(\Delta)$.

By Mitchell \cite{mitchell2009beginning}, $j$ is an iterated ultrapower of $\calK$ by its measures and extenders.
Recall that $W\cap \mathcal{K}=U$, and so,
$U=\{A\subseteq \kappa \mid A\in V, \kappa\in j_W(A)\}$.
So, this iterated ultrapower starts with $U$ or with an extender $F$ which normal measure is $U$.
\\Note that $M_F$ must be closed under $\kappa-$sequences.
Otherwise, there will be a set of ordinals $a, |a|<\kappa$ which consists of generators and which is not in $M_F$.
The further Easton support iteration of Prikry forcings will not be able to add such $a$.
Thus, by our assumption, the length of $F$ must be below first  measurable cardinal above $\kappa$ in $M_F$.
The iteration of Prikry forcings above $\kappa$ does not add new bounded subsets below the first measurable $>\kappa$.

By the same reason, extenders used to continue the iteration must be $\kappa-$closed.

None of them can be used infinitely many times (or infinitely many extenders cannot be used), since otherwise, $\omega-$sequences which cannot be added by
an Easton support iteration of Prikry forcings,  will be produced. It follows from the strong Prikry condition of the forcing, which can be shown for the relevant parts as in Ben Neria \cite{ben2021mathias}.
\\
This leaves us with a finite iteration by $\kappa-$closed extenders (measures).
\\It is the first part of the iteration.
\\The rest consisting of iteration of normal measures, each of them is applied $\omega-$many times.
\\
Take $N$ to be the first part of the iteration and $i:\calK \to N$ be the corresponding embedding.
\\
$\square$

\section{Properties of k} \label{Section: properties of k}

We continue and use the notations of theorem \ref{Theorem: Conditions imposed on i}. We first state the following lemma.

\begin{lemma}
Let $ P = P_{\kappa} $ be an Easton support iteration of Prikry type forcings, and $ i\colon V\to N $, $ \Delta\subseteq \kappa $, $ U\in V  $,  $ W\in V\left[G\right] $ and $ k\colon N\to M $ be as in section \ref{Section: The general Framework}.

Assume that there are no elements in $ \left( \kappa, \mbox{crit}(k) \right)\cap i(\Delta) $. Then $ \mbox{crit}(k)\in i(\Delta) $, namely, it  is the least element above $ \kappa $ in $ i\left( \Delta \right) $.
\end{lemma}

\begin{remark} \label{Remark: Crit k in Iterations of Prikry forcings}
	The assumption $ \left( \kappa, \mbox{crit}(k) \right)\cap i(\Delta) = \emptyset $ holds in the typical case where $ P = P_{\kappa} $ is an iteration of Prikry forcings. Indeed, assume, by contradiction, that there exists  $ \mu\in \left( \kappa, \lambda \right)\cap i(\Delta) $. Then $ \mu = k\left( \mu \right) $, and thus in $ M\left[j_W(G)\right] $, $ \mu $ changes cofinality to $ \omega $. Therefore, in $ V\left[G\right] $, $ \mbox{cf}\left( \mu \right) = \omega $, and, in $ V $, $ \mbox{cf}\left( \mu \right) \leq \kappa $. The sequence witnessing this belongs to $ V\cap \left(  {^\kappa}N \right) $ and thus, by our assumption on $ N $, belongs already to $ N $. This contradicts the measurability of $ \mu $ in $ N $.
\end{remark}

\pr
Denote $ \lambda = \mbox{crit}\left( k \right) $. Then for some $ h\in V $ and $ \kappa = \beta_0 < \beta_1 < \ldots < \beta_k $, 
$$ \lambda = i(h)\left( \kappa, \beta_1, \ldots, \beta_k \right) $$
By the definition of $ k $, $ \lambda >\kappa $. 

We first prove that $ \lambda \in i\left( \Delta \right) $. Assume otherwise. We can assume without loss of generality that for every $ \xi, \nu_1, \ldots, \nu_k $ below $ \kappa $, $ h\left( \xi, \nu_1, \ldots, \nu_k \right) > \xi $ does not belongs to $ \Delta $: this can be assumed by replacing the function $ h $ with the function $ h' \colon \left[ \kappa \right]^{n+1} \to \kappa $ defined as follows: For every $ \xi, \eta_1, \ldots, \eta_k $, $ h'\left( \xi, \eta_1, \ldots, \eta_k \right) $ equals $ h\left( \xi, \eta_1, \ldots, \eta_k \right) $ if $ h\left( \xi, \eta_1, \ldots, \eta_k \right) > \xi $ is not measurable in $ V $; and else,  $ h'\left( \xi, \eta_1, \ldots, \eta_k \right)  $ is an arbitrary non-measurable above $ \xi $. By our assumption, 
$$ i\left( h \right)\left(  \kappa, \beta_1, \ldots, \beta_k  \right) =  i\left( h' \right)\left(  \kappa, \beta_1, \ldots, \beta_k  \right) $$
so we can replace $ h $ with $ h' $. Since $ \lambda $ is regular (as a critical point of an elementary embedding), we can assume, using a similar argument, that each $  h\left( \xi, \nu_1, \ldots, \nu_k \right) $ is regular.

We can assume that for every $ \xi, \mu_1,\ldots ,\nu_k $, there are no elements of $ \Delta $ in the interval $ \left( \xi, h\left( \xi, \nu_1, \ldots, \nu_k \right) \right) $.

Let $ f\in V\left[G\right] $ be a function such that $ \left[f\right]_W = \lambda $. Then--
$$  \left[f\right]_W=\lambda < k(\lambda) = j_W\left( h \right)\left( \kappa,  \theta_{ \left[ f_{\beta_1} \right]_W }, \ldots, \theta_{  \left[ f_{ \beta_k  } \right]_W }  \right)  $$
By the definition of $ W $, there exists $ p\in G $ and $ r\in H $ such that--
$$ p^{\frown} r \Vdash  i(\lusim{f})(\kappa) < i(h)\big(\kappa, \theta_{ i\left( f_{ \beta_1 } \right)(\kappa)  } , \ldots,  \theta_{ i\left( f_{ \beta_k } \right)(\kappa)  } \big) $$
Recall that, for every $ 1\leq i \leq k $, there exists a condition in $ H $ forcing that $ \theta_{ i\left( f_{ \beta_i } \right)(\kappa)  } = \beta_i $. Thus by extending $ r $ inside $ H $, 
$$ p^{\frown} r \Vdash  i(\lusim{f})(\kappa) < i(h)\left( \kappa, \beta_1, \ldots, \beta_k \right) $$
Since there are no measurables of $ N $ in the interval $ \left( \kappa, i(h)\left( \kappa, \beta_0, \ldots, \beta_k \right) \right] $,  we can find $ r'\geq^* r $ inside $ H $ such that--
$$ p\Vdash \exists \alpha < i(h)\left( \kappa, \beta_1, \ldots, \beta_k \right), \ r'\Vdash i(f)(\kappa)< \alpha $$
and since $ P = P_{\kappa} $ is $ \kappa $-c.c. and $  i(h)\left( \kappa, \beta_1, \ldots, \beta_k \right) $ is regular, there exists $ \alpha < i(h)\left( \kappa, \beta_1, \ldots, \beta_k \right) $ such that--
$$ p^{\frown}r' \Vdash i(\lusim{f})(\kappa) < \alpha $$
Now apply $ k $ on both sides. By lemma \ref{lem3-1},  
$$ M\left[j_W(G)\right] \vDash \lambda = \left[ f \right]_W <  k\left( \alpha \right) $$
but $ \alpha < i(h)\left( \kappa, \beta_1, \ldots, \beta_k \right) = \lambda $ and thus $\lambda < k(\alpha) = \alpha < \lambda $, which is a contradiction.
$\square$

\begin{remark}
	Assume that $  P = P_{\kappa} $ is an iteration of the one point Prikry forcings. A one point Prikry forcing on a measurable $ \alpha $ is a forcing, which depends on a normal measure $ U $ on $ \alpha $, and is defined as follows: Conditions are of the form $ A $ where $ A\in U $ or $ \xi $ for some ordinal $ \xi<\alpha $. The latter kind of condition cannot be extended. A condition of the form $ A $ for $ A\in U $ can be extended in two ways: A direct extension is a condition $ B $ where $ B\in U $ and $ B\subseteq A $; a non-direct extension is of the form $ \xi $ where $ \xi\in A $ is an ordinal. 
	
	We argue that in this case, the question whether $ \left( \kappa, \mbox{crit}(k) \right)\cap i(\Delta) \neq \emptyset $, and, as a result, the value of $ \mbox{crit}(k) $,  depend of the choice of $ H $:
	
	\begin{enumerate}
		\item Denote by $ \mu $ the first element above $ \kappa $ in $ i(\Delta) $. Assume first that $ H $ is chosen such that the condition on coordinate $ \mu $ is a measure one set. In this case, $ \mu = \mbox{crit}(k) $. Indeed, $\mbox{crit}(k) < \mu$ cannot hold, since then $ \left( \kappa, \mbox{crit}(k) \right)\cap i(\Delta) =\emptyset $ which implies, by the last lemma, that $ \mu = \mbox{crit}(k) $. And $ \mu < \mbox{crit}(k) $ cannot hold since then $ k(\mu) = \mu $. Denote by $ \mu_0< \mu $ the one point added below $ \mu $ in $ j_W(G) $. Then $ H $ at coordinate $ \mu $ has a condition which is incompatible with $ \mu_0 $ (by shrinking the large set and applying a density argument). Thus $ \mu= \mbox{crit}(k) $.
		\item Denote now by $ \mu $ the least element in $ i(\Delta) $, for which $ H $ does not specify the one-point element added to it. We argue that $ \mbox{crit}(k) = \mu $, even though $ \mu $ doesn't have to be the least element above $ \kappa $ in $ i(\Delta) $.
		
		Repeat the proof of the last lemma, and note that the $ \leq^* $  forcing in the interval $ \left( \kappa, \mu \right) $ is trivial, since no condition in this interval can be non-trivially extended. This replaces the assumption that there are no elements of $ i\left(\Delta\right) $ in the interval $ \left( \kappa, i(h)\left( \kappa, \beta_1, \ldots, \beta_k   \right) \right) $. Therefore, $ \mu = \mbox{crit}(k) $. 
	\end{enumerate}
\end{remark}

Let us deal here with an Easton support iteration $P$ of the Prikry forcings over a set $\Delta$ of a measurable length $\kappa$.
Let $U$ be a normal ultrafilter over $\kappa$ in $V$ with $\Delta\not \in U$.
Let $G\subseteq P$ be a generic and $W$ be a normal ultrafilter in $V[G]$ which extends $U$.

Let $ i\colon V\to N $ be an elementary embedding as in theorem \ref{Theorem: Conditions imposed on i}, and assume that $ W = U_H $ and $ k\colon N\to M $ are as in lemma \ref{Lemma: k is elementary}. 

In the setting of iteration of Prikry forcings, much more can be said about the embedding $ k\colon N\to M $. From remark \ref{Remark: Crit k in Iterations of Prikry forcings}, it follows that $ \mbox{crit}(k) $ is the least element in $ i(\Delta) $ above $ \kappa $. In particular, by elementarity, $ k(\mu)\in j_W\left( \Delta \right) $ in $ M $, and thus a Prikry sequence is added to $ k(\mu) $ in $ j_W(G) $.

\begin{lemma}
	Denote $ \mu = \mbox{crit}(k) $. Then $\mu$ appears in the Prikry sequence of $k(\mu)$.
\end{lemma}

\begin{remark}
	$ \mu $ is not necessarily the first element in the Prikry sequence of $ k(\mu) $. The initial segment of this Prikry sequence below $ \mu $ depends on the choice of $ H $. For every finite sequence $ t\in \left[\mu\right]^{<\omega} $, we can choose $ H \subseteq i(P)\setminus \kappa$ such that $ t $ is an initial segment of the Prikry sequence of $ \mu $. This way, in $ M\left[ j_W(G) \right] $, $ t $ will be an initial segment of the Prikry sequence of $ k(\mu) $ below $ \mu $.
\end{remark}

\pr
	Let $ t $ be the finite initial segment of the Prikry sequence of $ k(\mu) $ below $ \mu $, and assume that $ \langle \xi, \eta_1, \ldots, \eta_l \rangle\mapsto t(\xi, \eta_1 , \ldots, \eta_l) $ is a function in $ V $, such that--
	$$ t=   i\left(  \langle \xi, \eta_1, \ldots, \eta_l \rangle\mapsto t(\xi, \eta_1 , \ldots, \eta_l)  \right)\left( \kappa, \beta_1, \ldots, \beta_l \right)$$ 
	for some generators $ \beta_1, \ldots, \beta_l $ of $ i $.
	For every $ \xi<\kappa $, let $ s(\xi) =\min\{\Delta\setminus \left(\xi+1\right)  \}  $, so $ \left[ \xi\mapsto s(\xi) \right]_W = \mu $.  In $ V\left[G\right] $, define, for every $ \xi<\kappa $,
	$$ \mu(\xi) = \mbox{the first element above }  t\left(\xi, \theta_{ f_{\beta_1 }(\xi) }, \ldots, \theta_{ f_{ \beta_k }(\xi) } \right)  \mbox{ in the Prikry sequence of }  s(\xi)$$
	and, if $ t\left(\xi, \theta_{ f_{\beta_1 }(\xi) }, \ldots, \theta_{ f_{ \beta_k }(\xi) } \right) $  is not an initial segment of the Prikry sequence of $ s(\xi) $, set $ \mu\left( \xi \right) =0 $.
	
	It suffices to prove that $ \left[ \xi\mapsto \mu(\xi) \right]_W = \mu $. 
	
	Assume first that $ \eta < \mu $. Work in $ N\left[G\right] $. Since $ H $ is $ \leq^* $-generic, it meets an element $ q\in i(P)\setminus \kappa $, for which $A^{q}_{\mu} \subseteq \mu\setminus \left(\eta+1\right) $. Since $ q\in H $, we can assume that $ t^{q}_{\mu} $ is an initial segment of $ t $: Indeed, $ t, t^{q}_{\mu} $ are compatible sequences, since, for any $ p\in G $ which forces that $ q\in H $ and decides the value of $ t^{q}_{\mu} $, the condition $ k\left( p^{\frown} q \right) = p^{\frown} k(q) $ belongs to $ j_W(G) $, and decides an initial segment, below $ \mu$, of the Prikry sequence of $ k(\mu) $. By our assumption, this initial segment is contained in $ t $, and $ p^{\frown} k(q) $ forces that every possible extension of it is above $ \eta $. Thus, in $ M\left[j_W(G)\right] $, each element in the Prikry sequence of $ k(\mu) $ after $ t $ is strictly above $ \eta $.
	
	The argument given in the previous paragraph also shows that for every $ q\in H $, $ t^{q}_{\mu}  $ is either empty or equals to $ t $: As mentioned, it must be an initial segment of $ t $. Let us argue that if it is proper, then it is empty. Apply the above paragraph for $ \eta = \max(t) $. Then by direct extending   $ q $ inside $ H $, it forces that the element after $ t^{q}_{\mu} $ in the Prikry sequence of $ \mu $ is strictly above $ \eta $. By applying $ k\colon N\to M $, there exists a condition in $ j_W(G) $ which forces that the Prikry sequence of $ k(\mu) $ has an initial segment $ t^{q}_{\mu} $, followed only by elements above $ \eta $. So $ t^{q}_{\mu} $ cannot be a proper initial segment of $ t $.

	Assume now that $ \eta < \left[ \xi \mapsto \mu(\xi) \right]_W $. Write $ \eta = \left[f\right]_W $ and assume that for every $ \xi<\kappa $, 
	$$  f(\xi) < \mu(\xi) < s(\xi) $$
	Let $ p\in G $ be a condition which forces this. Work in $ N\left[G\right] $. Take $ q\in H $ such that $ t^{q}_{\mu} = t $. Then $ i(p)^{\frown} q = p^{\frown} q$ forces that $i(\lusim{f})(\kappa)  $ is below the first element above $ t $ in the Prikry sequence of $ \mu $. Thus, its value can be decided by taking a direct extension. So, by direct extending $ q $ inside $ H $ we can assume that--
	$$ p\Vdash  \exists \alpha<\mu , \ q\Vdash i\left(\lusim{f}\right)(\kappa) < \alpha $$
	and thus there exists $ \alpha<\mu $ in $ V $, such that--
	$$ p^{\frown} q \Vdash i(f)(\kappa)< \alpha $$
	Thus, in $ M\left[j_W(G)\right] $, $ \eta =j_W(f)(\kappa) < k(\alpha) = \alpha < \mu $, as desired.
$\square$

In the next subsection we will decompose the embedding $ k $ to an iterated ultrapower of $ N $. We now demonstrate the first step in the iteration:

\begin{lemma}
	Let $ \mu = \mbox{crit}(k) $ and let $ U_{\mu} = \{ X\subseteq \mu \colon \mu \in k(X) \}\cap N $. Then $ U_{\mu}\in N $.
\end{lemma}

\pr
	For every $ \xi <\kappa $, denote by $ W_\xi $ the measure in $ V\left[G_{\xi}\right] $ used to singularize $ \xi $ in the Prikry forcing at stage $ \xi $ in the iteration. Let $ U_\xi = W_\xi\cap V $. We first argue that there exists a set $ \mathcal{F}\in N $ of measures on $ \mu $, with $ \left| \mathcal{F} \right| < \mu $, such that, for some $ p\in G $ and $ q\in H $, 
	\begin{equation} \label{Equation: lemma on Umu}
		 p^{\frown} q \Vdash i\left( \xi \mapsto \lusim{U}_{\xi} \right)(\mu) \in \mathcal{F} 
	\end{equation}
	Indeed, let $ \lusim{\alpha} $ be a $ j_U(P) $-name for the index of $ i\left( \xi \mapsto \lusim{U}_{\xi} \right)(\mu) $ in a prescribed well order of the normal measures $ \mu $ carries in $ N $. Work in $ N\left[G\right] $. For some $ q\in H $, there exists an ordinal $ \beta$ such that $ q\Vdash \lusim{\alpha} = \beta $. Thus, by $ \kappa-c.c. $ of the forcing $ i(P)_{\mu} = P_{\kappa} $, there exist $ p\in G $ and a set $S\subseteq 2^{2^{\mu}}$ of ordinals with $ \left|S\right| < \mu $, such that $ {p}^{\frown} q \Vdash \lusim{\alpha} \in S $. In particular, $  {p}^{\frown} q $ forces that $ i\left( \xi \mapsto \lusim{U}_{\xi} \right)(\mu)  $ belongs to $ \mathcal{F} $, where $ \mathcal{F} $ is the set of measures on $ \mu $ indexed in $ S $.
	
	Now apply $ k $ on equation (\ref{Equation: lemma on Umu}), and work in $ M\left[ j_W(G) \right] $.  Since $ \left|  \mathcal{F}\right|<\mu $,  it follows that there exists a measure $ F\in \mathcal{F} $ such that--
	$$ j_W\left( \xi \mapsto U_{\xi} \right)(k\left(\mu\right)) = k\left(F\right) $$
	so it suffices to argue that $ F =\{ X\subseteq \mu \colon \mu\in k(X) \}\cap N $. Fix $ X\in F $. Write $ X = i(g)\left(\kappa, \beta_0 ,\ldots, \beta_k \right) $. Then--
	$$ j_W(g)\left( \kappa,  \theta_{ \left[ f_{\beta_1 } \right]_W }, \ldots,  \theta_{ \left[ f_{\beta_k } \right]_W }   \right) \in   j_W\left( \xi \mapsto U_{\xi} \right)\left(  k(\mu) \right) $$
	Recall the function $ \xi\mapsto s(\xi) = \min\left(  \Delta \setminus \left(\xi+1\right) \right) $, for which $ \left[ \xi\mapsto s(\xi) \right]_W = k(\mu) $. We can assume that for every $ \xi<\kappa $,
	$$  g\left(  \xi, \theta_{f_{\beta_1}(\xi)}, \ldots, \theta_{f_{ \beta_k }(\xi)} \right)   \in U_{  s(\xi) } $$
	and let $ p\in G $ be a condition which forces this. Then for strong enough $ q\in H $, 
	$$  p^{\frown} q \Vdash i(g)\left( \kappa, \beta_1, \ldots, \beta_k \right) \in i\left( \xi \mapsto \lusim{U}_{\xi} \right)\left( \mu \right) $$
	and thus by direct extending $ q $ further, we can assume that $ q $ forces that the first element after $ t $ in the Prikry sequence of $ \mu$ belongs to $ i(g)\left( \kappa, \beta_1, \ldots, \beta_k \right) = X $. Thus $ k(q)\in j_W(G) $ forces that the first element after $ t $ in  the Prikry sequence of $ k(\mu) $ belongs to $ k(X) $. By the previous lemma, it follows that $ \mu\in k(X) $, as desired.
$\square$\\

\subsection{Description of $ j_W\restriction_{V} $} \label{Subsection, structure of jW restriction V}

We now generalize the previous subsection, in order to completely decompose $ j_W\restriction_{V} $. We continue to assume that $ P = P_{\kappa} $ is an iterations of Prikry forcings. For technical reasons, we will assume that the measures used in the iteration $ P = P_{\kappa} $ to singularize  the measurables in $ \Delta $ are all simply generated; this is needed only in the proof of claim \ref{Claim: Inductive step For multivar Fusion} which will be presented in the next subsection.

At each stage $ \alpha \in \Delta$, let $ \lusim{Q}_{\alpha} $ be the $ P_{\alpha} $-name for the Prikry forcing on $ \alpha $, using a simply generated normal measure $ \lusim{W}_{\alpha} $ on $ \alpha $. Denote $ \lusim{U}_{\alpha} = \lusim{W}_{\alpha}\cap V \in V $. Let $ \lusim{H}_{\alpha} \subseteq \left(j_{ \lusim{U}_{\alpha} }\left( P_{\alpha} \right) \setminus \alpha ,  \leq^*   \right)$, $ \lusim{H}_{\alpha}\in  V\left[G_{\alpha}\right] $, be $\leq^*$-generic over $ M_{\lusim{U}_\alpha}\left[G_{\alpha}\right] $, such that $ \lusim{W}_{\alpha} = \left( \lusim{U}_{\alpha} \right)_{ \lusim{H}_{\alpha} } $.

Let $ G\subseteq P_{\kappa} $ be generic over $ V $. 

Our goal is to prove the following theorem:

\begin{theorem} \label{Theorem: Measure generated from i restrics to iteration of N}
	Let $ H\in V\left[G\right] $ be a generic set for $ \langle i(P)\setminus \kappa, \leq^* \rangle $ which satisfies $ (*) $. Let $ W = U_{H} $  be the corresponding normal measure on $ \kappa $ extending $ U $, and denote its ultrapower embedding $ j_W\colon V\left[G\right]\to M\left[ j_W(G) \right] \simeq \mbox{Ult}\left( V\left[G\right], W \right) $ for some model $ M $.
	Then $ j_W\restriction_{V} $ factors to the form $ j_W\restriction_{V} = k\circ i $ for some elementary $ k\colon N\to M $.
	
	Moreover, if $ P $ is an Easton support iteration, where at each step $ \beta\in \Delta $, $ \lusim{Q}_{\beta} $ is forced to be Prikry forcing with a simply generated normal measure on $\beta$,  then $ k $ is an iterated ultrapower of $ N $ by normal measures and $ j_W(\kappa) = i(\kappa) $. 
\end{theorem}

This, in contrast to Full-Support and Nonstationary-Support iterations of Prikry forcings, where, assuming $ \mbox{GCH}_{\leq \kappa} $, $ j_W\restriction_{V} $ is an iteration of $ V $ by normal measures only. 

If all the measures considered, including $ W $, are simply generated, $ j_W\restriction_{V} $ is an iterated ultrapower by normal measures only:

\begin{theorem} \label{Theorem: Restrictions of ultrapwers with simply generated measures}
	Assume that $ P $ is an Easton support iteration, where at each step $ \beta\in \Delta $, $ \lusim{Q}_{\beta} $ is forced to be Prikry forcing with a simply generated normal measure on $\beta$. Then for every simply generated measure $ W\in V\left[G\right] $ on $ \kappa $, $ j_W\restriction_{V} $ is an iteration of $ V $ by normal measures. Moreover, if $ U = W\cap V $ then $j_W(\kappa) =  j_U(\kappa) $.
\end{theorem}

We will prove theorems \ref{Theorem: Restrictions of ultrapwers with simply generated measures} and \ref{Theorem: Measure generated from i restrics to iteration of N} simultaneously. Assume that $ H\in V\left[G\right] $ is a  generic for $ \langle i(P)\setminus \kappa, \leq^* \rangle $ over $ N\left[G\right] $ with the property $ (*) $. In the case where $ i = j_U $ and $ N = M_U $, any generic for $ \langle i(P)\setminus \kappa \rangle , \leq^*\rangle $ is such. Let $ W = U_{H}\in V\left[G\right] $  be the corresponding normal measure on $ \kappa $. Let $ j_W\colon V\left[G\right]\to M\left[ j_W(G) \right] $ be the corresponding ultrapower embedding. 

Denote by $ B\subseteq \left( \kappa, i(\kappa) \right) $ the set of generators of $ i $. By property $ (*) $ of $ H $, for every $ \beta\in B $, there exists a function $ f_{\beta} $ in $ V $ such that $ H $ forces that $ \beta = \theta_{i\left( f \right)(\kappa) }  $. The mapping $ \beta\mapsto f_{\beta} $ is available in $ V\left[G\right] $.

Recall the embedding  $ k\colon N\to M $ defined in lemma \ref{Lemma: k is elementary}:
$$  k\left( i(f)\left( \kappa, \beta_1 ,\ldots, \beta_k \right) \right) = j_W(f)\left(  \kappa, \theta_{\left[f_{\beta_1}\right]_W} , \ldots, \theta_{\left[ f_{ \beta_k }  \right]_W}  \right)  $$
for every $ f\in V $ and $ \beta_1, \ldots, \beta_k \in B $. Then $ k $ is elementary, $\mbox{crit}(k) > \kappa$ and $ j_W\restriction_{V} = k\circ i $.

Denote $ \kappa^* = i(\kappa) $. Define by induction a linear directed system $ \langle \langle  M_{\alpha} \colon \alpha\leq \kappa^*  \rangle, \langle j_{\alpha, \beta} \colon \alpha< \beta \leq \kappa^* \rangle \rangle $ such that:
\begin{enumerate}
	\item  $ M_0 = N $, $ j_0 = i $.
	\item \textbf{Successor Step:} Assume that $ \alpha<\kappa^* $ and $ M_{\alpha} $ has been defined. We will define an elementary embedding $ k_{\alpha} \colon M_{\alpha} \to M $, such that $ j_W\restriction_{V} = k_{\alpha} \circ j_{\alpha} $. We denote $ \mu_{\alpha} = \mbox{crit}\left(  k_{\alpha} \right) $ and define--
	$$ U_{\mu_{\alpha} } = \{ X\subseteq \mu_{\alpha} \colon \mu_{\alpha} \in k_{\alpha}(X) \}\cap M_{\alpha} $$
	We will prove that $ U_{ \mu_{\alpha} }\in M_{\alpha} $ and take $ M_{\alpha+1} \simeq \mbox{Ult}\left(  M_{\alpha}, U_{ \mu_{\alpha}  } \right) $. We also take $ j_{\alpha, \alpha+1} \colon M_{\alpha}\to M_{\alpha+1} $ to be the ultrapower embedding $ j^{ M_{\alpha} }_{ U_{ \mu_{\alpha} } } $, and $ j_{\alpha+1} = j_{\alpha, \alpha+1} \circ j_{\alpha} $.
	\item \textbf{Limit Step:} For every limit $ \alpha\leq\kappa^* $, the system $ \langle M_{\beta} \colon \beta<\alpha \rangle, \langle j_{\beta, \gamma} \colon \beta< \gamma < \alpha \rangle $ is linearly directed, and we take direct limit to form the model $ M_{\alpha} $ and the embedding $ j_{\alpha} \colon V\to M_{\alpha}  $.
	
\end{enumerate}

For every $ \alpha<\kappa^* $, define $ k_{\alpha} \colon M_{\alpha}\to M $ as follows:
$$ k_{\alpha}\left(   j_{\alpha}\left( f \right)\left(  \kappa,  j_{0,\alpha}(\beta_1), \ldots, j_{0,\alpha}\left( \beta_l \right), \mu_{\alpha_1}, \ldots, \mu_{\alpha_k}   \right)  \right)  =  j_W\left( f \right)\left(  \kappa, \theta_{ \left[  f_{\beta_1} \right]_W} , \ldots,   \theta_{\left[ f_{ \beta_l }\right]_W}  , \mu_{\alpha_1}, \ldots, \mu_{\alpha_k}  \right) $$
for every $ f\in V $, $ \beta_1, \ldots, \beta_l $ generators of $ i $  and $ \alpha_1< \ldots < \alpha_k < \alpha $.

Our goal is to prove by induction on $ \alpha<\kappa^* $ the following properties:
\begin{enumerate}[label=(\Alph*)]
	\item \label{Property:   k alpha is elementary }
	 $ k_{\alpha} \colon M_{\alpha}\to M $ is an elementary embedding, and $ j_W\restriction_{V} = k_{\alpha} \circ j_{\alpha} $.
	\item \label{Property:   crit k alpha} $\mu_{\alpha}$ is measurable in $ M_{\alpha} $. Moreover, it is the least measurable in $ j_\alpha\left( \Delta \right) $, which is greater or equal to $ \mbox{sup}\{ \mu_{\beta} \colon \beta<\alpha \} $, and whose cofinality is above $ \kappa $ in $ V $.
	\item \label{Property:   mu alpha appears in the Prikry sequence} $ \mu_{ \mu_{\alpha} } $ appears in the Prikry sequence of $ k_{\alpha}\left(  \mu_{\alpha} \right) $.
	\item \label{Property:   U mu alpha}
	 Let  $ U_{\mu_{\alpha}} $ be defined in $ V\left[G\right] $ as above. Then $ U_{ \mu_{\alpha} }\in M_{\alpha} $ is a normal measure which concentrates on $ \mu_{\alpha} \setminus j_{\alpha}\left( \Delta \right) $. Moreover,
	$$  k_{\alpha}\left( U_{ \mu_{\alpha} } \right) = j_{W}\left(  \delta\mapsto U_{\delta} \right)\left( k_{\alpha}\left( \mu_{\alpha} \right) \right) $$
	where, for every $ \delta\in \Delta $, $ U_{\delta} = W_{\delta}\cap V $, for $ W_\delta $ which is the measure used in the Prikry forcing at stage $ \delta $ in the iteration $ P $.
\end{enumerate}

After that, we will prove in lemma \ref{Lemma:  k_kappa* is the identity},  that $ k_{\kappa^*} \colon M_{\kappa^*}\to M $ is the identity, and thus $ j_W\restriction_{V} = j_{\kappa^*} $. This will conclude the proof of theorems \ref{Theorem: Restrictions of ultrapwers with simply generated measures} and \ref{Theorem: Measure generated from i restrics to iteration of N}.

\begin{remark}
	We remark that $ k_{\alpha} $ is well defined is the sense that there is no $ \alpha'<\alpha $ and generator $ \beta $ of $ i $, for which $ j_{ 0, \alpha }(\beta) = \mu_{\alpha'} $. Indeed, assume otherwise. Note that $ \mu_{\alpha'} = j_{0,\alpha}\left( \beta \right) \geq j_{ 0,\alpha' }(\beta) $. Strict inequality is not possible here, since if $ j_{ 0,\alpha' }(\beta) <\mu_{\alpha'} $ then $ j_{0,\alpha'}(\beta) = j_{0,\alpha}(\beta) = \mu_{\alpha'} $, which is a contradiction. Thus, $ j_{0,\alpha'}\left( \beta \right) = \mu_{\alpha'} $ (which is, by itself, possible for $ \alpha'<\alpha $ - see remark \ref{Remark - a generator can be a critical point}), but then, applying $ j_{\alpha',\alpha} $ on both sides, we get--
	$$ j_{0,\alpha}(\beta) = j_{\alpha',\alpha}\left(  \mu_{\alpha'} \right) > \mu_{\alpha'} $$
	where the last inequality follows since $ \mu_{\alpha'} = \mbox{crit}\left( j_{\alpha',\alpha} \right) $.
\end{remark}

\begin{remark} \label{Remark - a generator can be a critical point}
	It is possible that a generator $ \beta $ of $ i $ is measurable in $ N $ and belongs to $ i(\Delta) $. In this case, there exists $ \alpha<\kappa^* $ such that $ \mu_{\alpha} = \beta = j_{0,\alpha}(\beta) $. Such $ \beta $ will appear as an element in the Prikry sequence of $ k_{\alpha}\left( \beta \right)\in j_W(\Delta) $, which also has the form $ \theta_{ \left[ f_{ \beta } \right]_W } $.
\end{remark}

Properties $ (A)-(D) $ of $ k_{\alpha} $, presented above, will be proved by induction on $ \alpha<\kappa^* $. The proof of the inductive step at stage $ \alpha<\kappa^* $ will be carried out in subsection \ref{Subsection: Properties of kalpha}, using the tools presented in  \cite{NonStatRestElm} and  \cite{kaplan2022magidor}. Fixing $ \alpha<\kappa^* $, we can assume by induction that $k_{\alpha'} \colon M_{\alpha'}\to M$ and  $\mu_{\alpha'}, U_{ \mu_{\alpha'} } $, for $ \alpha'<\alpha $, satisfy properties $ (A)-(D) $. Denote by $ t_{\alpha'}$  the initial segment of the Prikry sequence of $ k_{\alpha'}\left( \mu_{\alpha'} \right) $ below $ \mu_{\alpha'} $. 

\begin{defn}
	Fix $ \alpha < \kappa^* $ and a sequence of generators $ \langle \beta_1, \ldots, \beta_l \rangle $ for $ i $. An increasing sequence $ \langle \alpha_1, \ldots, \alpha_k \rangle $ below $ \alpha $ is called a $ \langle \beta_1,\ldots, \beta_l \rangle $-nice sequence if there are functions $ g_1, \ldots, g_k, t_1, \ldots, t_k$ in V, such that--
	$$ \mu_{\alpha_1} = j_{\alpha_1}\left( g_1 \right)\left(  \kappa, j_{ 0,\alpha_1 }\left(\beta_1\right), \ldots, j_{ 0, \alpha_1 }\left(\beta_l\right) \right) $$
	$$ t_{\alpha_1} = j_{\alpha_1}\left( t_{\alpha_1} \right)\left(  \kappa, j_{ 0,\alpha_1 }\left(\beta_1\right), \ldots, j_{ 0, \alpha_1 }\left(\beta_l\right) \right) $$
	$$ U_{\mu_{\alpha_1}} = j_{\alpha_1}\left( F_1 \right)\left(  \kappa, j_{ 0,\alpha_1 }\left(\beta_1\right), \ldots, j_{ 0, \alpha_1 }\left(\beta_l\right) \right) $$
	and, for every $ 1\leq i < k $,
	$$ \mu_{\alpha_{i+1}} = j_{\alpha_{i+1}}\left( g_{i+1} \right)\left(  \kappa, j_{ 0,\alpha_1 }\left(\beta_1\right), \ldots, j_{ 0, \alpha_1 }\left(\beta_l\right), \mu_{\alpha_1}, \ldots, \mu_{\alpha_i} \right) $$
	$$ t_{\alpha_{i+1}} = j_{\alpha_{i+1}}\left( t_{i+1} \right)\left(  \kappa,j_{ 0,\alpha_1 }\left(\beta_1\right), \ldots, j_{ 0, \alpha_1 }\left(\beta_l\right), \mu_{\alpha_1}, \ldots, \mu_{\alpha_i} \right) $$
	$$ U_{\mu_{\alpha_{i+1}}} = j_{\alpha_{i+1}}\left( F_{i+1} \right)\left(  \kappa, j_{ 0,\alpha_1 }\left(\beta_1\right), \ldots, j_{ 0, \alpha_1 }\left(\beta_l\right), \mu_{\alpha_1}, \ldots, \mu_{\alpha_i} \right) $$
\end{defn}

Fix now $ \alpha<\kappa^* $. Assume by induction that properties $(A)-(D)$  above hold for every $ \alpha'<\alpha $. Fix also a sequence of generators $ \langle \beta_1, \ldots, \beta_l \rangle $ for $ i $, and a $ \langle \beta_1, \ldots, \beta_l \rangle $-nice sequence $ \langle \alpha_1, \ldots, \alpha_k \rangle  $ below $ \alpha $. We define, in $ V\left[G\right] $, functions which can be used to represent $ \mu_{\alpha_i}, \  t_{\alpha_i} , U_{\alpha_i}$. Assume that $ \mu_{\alpha_i} $ is the $ n_{i} $-th element in the Prikry sequence of $ k_{ \alpha_{i} }\left( \mu_{\alpha_i} \right) $.

First, set--
$$ \mu_{\alpha_1}(\xi) = \mbox{the } n_1\mbox{-th element in the Prikry sequence of } g_1(\xi, \theta_{f_{ \beta_1 }(\xi)} , \ldots, \theta_{ f_{ \beta_l }(\xi) }   ) $$
By induction, define, for every $ i<k $,
\begin{align*}
	 \mu_{\alpha_{i+1}}(\xi) = &\mbox{the } n_{i+1}\mbox{-th element in the Prikry sequence of }\\
	 &g_{i+1}(\xi, \theta_{f_{ \beta_1 }(\xi)} , \ldots, \theta_{ f_{ \beta_l }(\xi) } , \mu_{\alpha_1}(\xi), \ldots, \mu_{\alpha_i}(\xi)    ) 
\end{align*}
and $ U_{ \mu_{\alpha_i} }(\xi) = W_{ \mu_{\alpha_i}(\xi) }\cap V $. Here, given  $ \delta\in \Delta $, $ W_{ \delta } $ is the measure on $ \delta $ used in the Prikry forcing which was applied at stage $ \delta $  in the  iteration.

\begin{claim}
	$ \left[  \xi \mapsto \mu_{\alpha_i}(\xi) \right]_W = \mu_{\alpha_i} $ and $ \left[ \xi \mapsto U_{ \mu_{\alpha_i}(\xi) } \right]_W = k_{\alpha_i}\left(U_{ \mu_{\alpha_i} } \right)$.
\end{claim}

\pr
	We begin by proving that $ \left[  \xi \mapsto \mu_{\alpha_i}(\xi) \right]_W = \mu_{\alpha_i} $. We present the argument for $ i=1 $. Higher values of $ i\leq k $ are proved similarly, using induction. Recall that--
	$$  \mu_{\alpha_1} = j_{\alpha_1}\left(  g_1 \right)(\kappa, j_{0,\alpha_1}\left(\beta_1\right), \ldots, j_{ 0,\alpha_1 }\left(\beta_l\right) ) $$
	and by applying $ k_{\alpha_1} $ on both sides,
	$$  k_{\alpha_1}\left(\mu_{\alpha_1}\right) = j_{W}\left( g_1   \right)(\kappa,  \theta_{\left[  f_{ \beta_1 }(\xi)  \right]_W }   ,   \ldots, \theta_{ \left[  f_{ \beta_l }(\xi)  \right]_W }  )  $$
	By induction, $ \mu_{\alpha_1} $ is the $ n_1 $-th element in the Prikry sequence of $ k_{\alpha_1}\left( \mu_{\alpha_1} \right) $, and thus it is represented as the $ n_1 $-th element in the Prikry sequence of $ g_1\left(   \xi, \theta_{f_{ \beta_1 }(\xi)} , \ldots, \theta_{ f_{ \beta_l }(\xi) }  \right) $.
	
	As for  $ \left[ \xi \mapsto U_{ \mu_{\alpha_i}(\xi) } \right]_W = k_{\alpha_i}\left(U_{ \mu_{\alpha_i} } \right)$, this follows since, by induction, $$ k_{\alpha_i}\left( U_{ \mu_{\alpha_i} } \right) = j_W\left( \delta\mapsto U_\delta \right)\left(  k_{ \alpha_i }\left( { \mu_{\alpha_i} } \right) \right) $$
$\square$

Let us argue that $ k_{\alpha} \colon M_{\alpha}\to M $ is elementary. 

\begin{lemma}
	$ k_{\alpha} \colon M_{\alpha}\to M $ is elementary.
\end{lemma}

\pr
Assume that $ x,y\in M_{\alpha} $, and let us prove, for example, that $ x\in y $ if and only if $ k(x) \in  k(y) $. Let $ f,g\in V $, $\beta_1, \ldots, \beta_l $ and $ \alpha_1 < \ldots < \alpha_k< \alpha $ be such that--
$$  x = j_{\alpha}(f)\left( \kappa, j_{0,\alpha}\left(\beta_1\right), \ldots, j_{0,\alpha}\left(\beta_l\right), \mu_{\alpha_1}, \ldots, \mu_{\alpha_k} \right)  \ , \ y = j_{\alpha}(g)\left( \kappa, j_{0,\alpha}\left(\beta_1\right), \ldots, j_{0,\alpha}\left(\beta_l\right), \mu_{\alpha_1}, \ldots, \mu_{\alpha_k} \right)$$
Assume that $ \alpha = \alpha'+1 $ is successor (the limit case is simpler). For simplicity, we assume also that $ \alpha_k = \alpha' $. Then $ x\in y $ if and only if--

\begin{align*}
\mu_{\alpha'}\in j_{\alpha',\alpha}  & \left( \{     \xi< \mu_{\alpha'} \colon  j_{\alpha'}(f)\left( \kappa, j_{0,\alpha'}(\beta_1), \ldots, j_{0,\alpha'}(  \beta_l), \mu_{\alpha_1}, \ldots, \mu_{\alpha_{k-1}}, \xi \right)\in \right. \\
& \left. j_{\alpha'}(g)\left( \kappa, j_{0,\alpha'}(\beta_1), \ldots, j_{0,\alpha'}(  \beta_l), \mu_{\alpha_1}, \ldots, \mu_{\alpha_{k-1}}, \xi \right) \} \right)
\end{align*}

which is equivalent to--
\begin{align*}
\{  \xi< \mu_{\alpha'} \colon & j_{\alpha'}(f)\left( \kappa, j_{0,\alpha'}(\beta_1), \ldots, j_{0,\alpha'}(  \beta_l), \mu_{\alpha_1}, \ldots, \mu_{\alpha_{k-1}}, \xi \right)\in  \\
& j_{\alpha'}(g)\left( \kappa, j_{0,\alpha'}(\beta_1), \ldots, j_{0,\alpha'}(  \beta_l), \mu_{\alpha_1}, \ldots, \mu_{\alpha_{k-1}}, \xi \right)   \} \in U_{  \mu_{\alpha'} } 
\end{align*}
which, by the definition of $ U_{ \mu_{\alpha'} } $, is equivalent to--
\begin{align*}
\mu_{\alpha'} \in k_{\alpha'}&	\left(\{  \xi< \mu_{\alpha'} \colon  j_{\alpha'}(f)\left( \kappa, j_{0,\alpha'}(\beta_1), \ldots, j_{0,\alpha'}(  \beta_l), \mu_{\alpha_1}, \ldots, \mu_{\alpha_{k-1}}, \xi \right)\in \right. \\
&\left.	 j_{\alpha'}(g)\left( \kappa, j_{0,\alpha'}(\beta_1), \ldots, j_{0,\alpha'}(  \beta_l), \mu_{\alpha_1}, \ldots, \mu_{\alpha_{k-1}}, \xi \right)   \}\right)
\end{align*}
namely $ k_{\alpha}(x) \in  k_{\alpha}(y) $.\\
$\square$\\

Let us describe now the main ideas behind the proof that $ \mu_{\alpha} = \mbox{crit}\left( k_{\alpha} \right) $ is measurable in $ M_{\alpha} $.  Note that this is not trivial since $ k_{\alpha} \colon M_{\alpha}\to M $ is not definable in $ M_{\alpha} $.  
The full argument will be presented in lemma \ref{Lemma: mu alpha is measurable in M alpha}, but will require a technical theorem (theorem \ref{Theorem: Multivariable Fusion Replacement}). Mainly we would like to follow the methods developed in \cite{NonStatRestElm} and  \cite{kaplan2022magidor}, which deal with nonstationary and full support iterations of Prikry forcings, respectively. 

We  consider the function $ f\in V\left[G\right] $, for which $ \mu_{\alpha} = \left[ f \right]_W $. We will prove that if $ \mu_{\alpha} $ is not measurable in $ M_{\alpha} $, then $\mu_{\alpha} = \left[f\right]_W \in \mbox{Im}\left( k_{\alpha} \right) $, contradicting the fact that $ \mu_{\alpha} = \mbox{crit}\left( k_{\alpha} \right) $. For that, we first fix a function $ h\in V $ such that, for some sequence $ \beta_1, \ldots, \beta_l$ of generators of $ i $, and for some nice sequence $ \langle \alpha_1, \ldots, \alpha_k \rangle $ below $ \alpha $,
$$  \mu_{\alpha} = j_{\alpha}\left( h \right)\left( \kappa, j_{0,\alpha}\left( \beta_1\right),\ldots, j_{0,\alpha}\left(\beta_l\right), \mu_{\alpha_1}, \ldots, \mu_{\alpha_k} \right)  $$
since $ \mu_{\alpha} = \mbox{crit}\left( k_{\alpha} \right) $, we can assume that for every $ \xi<\kappa $,
$$  f(\xi) < h\left(  \xi, \theta_{f_{ \beta_1 }(\xi)} , \ldots, \theta_{ f_{ \beta_l }(\xi) } ,  \mu_{\alpha_1}(\xi), \ldots, \mu_{\alpha_k}(\xi)  \right)  $$
Pick a condition $ p\in G $ which forces this. For every $ \xi<\kappa, \ \vec{\eta} = \langle  \eta_1, \ldots, \eta_l \rangle $ and $\vec{\nu} = \langle  \nu_1,\ldots, \nu_k \rangle $, denote--
$$  e\left( \xi, \vec{\eta}, \vec{\nu} \right) = \{  r\in P\setminus \nu_k \colon \mbox{ there exists a bounded subset } A\subseteq h\left( \xi, \vec{\eta}, \vec{\nu} \right)  \mbox{ such that } r\Vdash \lusim{f}(\xi) \in A   \}$$
This set is $ \leq^* $-dense open above conditions which extend $ p $ and force that--
\begin{equation} \label{Equation: The correct values}
	\langle \theta_{f_{ \beta_1 }(\xi)} , \ldots, \theta_{ f_{ \beta_l }(\xi) } , \mu_{\alpha_1}(\xi), \ldots, \mu_{\alpha_k}(\xi)  \rangle = \langle \vec{\eta}, \vec{\nu} \rangle 
\end{equation}
We would like to follow \cite{NonStatRestElm} and \cite{kaplan2022magidor},  and construct a condition $ p^*\in G $ above $ p $, such that, very roughly\footnote{We omitted some of the details in the version described here, for sake of simplicity.}, for every $ \xi, \vec{\eta}, \vec{\nu} $ as above, and for every extension $ r $ of $ p^* $ which forces (\ref{Equation: The correct values}), 
$$ r\restriction_{ \nu_k } \Vdash r\setminus_{\nu_k} \in e\left( \xi, \vec{\eta}, \vec{\nu} \right) $$
Essentially, such $ p^* $ will have the following property: every extension $ r $ of it which forces that equation (\ref{Equation: The correct values}) holds, forces also that $ f\left( \xi \right) $ belongs to a bounded subset $ A\left( \xi, \vec{\eta}, \vec{\nu} \right)\subseteq h\left( \xi, \vec{\eta}, \vec{\nu} \right) $ (which depends only on $ p^* $ and $ \langle \xi, \vec{\eta}, \vec{\nu} \rangle $, and not on the choice of the extension of $ p^* $ which forces (\ref{Equation: The correct values})). In \cite{NonStatRestElm} and \cite{kaplan2022magidor} the construction of such $ p^* $ was done by a Fusion argument which allows, in a sense, to absorb a lot of data into a single direct extension  $ p^* $ of $ p $. Such a method is not available in the Easton support iteration. We bypass this problem by constructing, for every sequence $ \langle \xi, \eta_1, \ldots, \eta_l \rangle $, a system of non-direct extensions of $ p $,
$$  \langle p\left(   \xi, \eta_1, \ldots, \eta_l, \nu_1, \ldots, \nu_k  \right) \colon \nu_1 <\ldots < \nu_k <\kappa \rangle $$
and sets--
$$ \langle A\left(   \xi, \eta_1, \ldots, \eta_l, \nu_1, \ldots, \nu_k  \right) \colon \nu_1 <\ldots < \nu_k <\kappa  \rangle $$
such that the following properties hold:
\begin{enumerate}
	\item If $   p\left(   \xi, \eta_1, \ldots, \eta_l, \nu_1, \ldots, \nu_k  \right) $ forces (\ref{Equation: The correct values}), then it also forces that $ \lusim{f}(\xi) \in A\left( \xi, \vec{\eta}, \vec{\nu}  \right) $, which is a bounded subset of $ h\left( \xi, \vec{\eta}, \vec{\nu} \right) $.
	\item For a set of $ \xi $-s in $ W $, $  p\left(   \xi, \theta_{f_{ \beta_1 }(\xi)} , \ldots, \theta_{ f_{ \beta_l }(\xi) } , \mu_{\alpha_1}(\xi), \ldots, \mu_{\alpha_k}(\xi)  \right)$ belongs to $G$.
\end{enumerate}
This suffices, since, by combining the above properties, 
$$ V\left[G\right] \vDash  \{  \xi<\kappa \colon f\left( \xi \right)\in A\left(  \xi, \theta_{f_{ \beta_1 }(\xi)} , \ldots, \theta_{ f_{ \beta_l }(\xi) } , \mu_{\alpha_1}(\xi), \ldots, \mu_{\alpha_k}(\xi) \right) \}\in W$$
and thus, in $ M\left[j_W(G)\right] $,
\begin{align*}
	\mu_{\alpha} = \left[  f\right]_W \in & \left[  \xi \mapsto  A\left(  \xi, \theta_{f_{ \beta_1 }(\xi)} , \ldots, \theta_{ f_{ \beta_l }(\xi) } , \mu_{\alpha_1}(\xi), \ldots, \mu_{\alpha_k}(\xi) \right)
	\right]_W\\
	=& k_{\alpha}\left( \  j_{\alpha}\left(  \langle \xi, \vec{\eta} ,\vec{\nu} \rangle\mapsto A\left( \xi, \vec{\eta}, \vec{\nu} \right)  \right)\left( \kappa, \beta_1, \ldots, \beta_l, \mu_{1}, \ldots, \mu_k \right)   \  \right)\subseteq  \mbox{Im}\left( k_{\alpha}   \right) 
\end{align*} 
where the last inclusion follows since $  j_{\alpha}\left(  \langle \xi, \vec{\eta} ,\vec{\nu} \rangle\mapsto A\left( \xi, \vec{\eta}, \vec{\nu} \right)  \right)\left( \kappa, \beta_1, \ldots, \beta_l, \mu_{1}, \ldots, \mu_k \right)   $ is a bounded subset of $ \mu_{\alpha} = j_{\alpha} \left(h\right)\left( \kappa, \beta_1, \ldots, \beta_l, \mu_{\alpha_1}, \ldots, \mu_{\alpha_k} \right) $.

We will complete the missing details in the proof in lemma \ref{Lemma: mu alpha is measurable in M alpha}. Before that, we present the proof of theorem  \ref{Theorem: Multivariable Fusion Replacement}.

\subsection{Theorem \ref{Theorem: Multivariable Fusion Replacement} and its proof}

We devote this subsection to the proof of the following theorem:

\begin{theorem} \label{Theorem: Multivariable Fusion Replacement}
	Let $ p\in G $ be a condition.  Assume that for every increasing sequence $ \langle \xi, \nu_1, \ldots, \nu_k \rangle $, and for every $ \vec{\eta} = \langle \eta_1, \ldots, \eta_l \rangle $ above $ \xi $, the set--
	$$ e\left(  \xi , \eta_1, \ldots, \eta_l, \nu_1, \ldots, \nu_k \right) \subseteq P\setminus \nu_k$$
	is $ \leq^* $ dense open above conditions in $ P\setminus \nu_k $ which force that--
	$$  \langle \eta_1, \ldots, \eta_l, \nu_1, \ldots, \nu_k \rangle = \langle  \theta_{f_{ \beta_1 }(\xi)} , \ldots, \theta_{ f_{ \beta_l }(\xi) } , \mu_{\alpha_1}(\xi), \ldots, \mu_{\alpha_k}(\xi) \rangle $$
	Then there are $ s<\omega $, a new sequence of generators $ \beta'_{l},\ldots, \beta'_{s}  $ of $ i $ which contains $ \beta_1, \ldots, \beta_l $, and a system of extensions of $ p $, 
	$$ \langle p\left( \xi, \eta_1, \ldots , \eta_{s}, \nu_1, \ldots, \nu_k\right) \colon \eta_1,\ldots, \eta_{s}< \kappa, \nu_1< \ldots < \nu_{k}<\kappa \rangle $$
	with the following properties: 
	\begin{enumerate}
		\item There exists a set of $ \xi $-s in $ W $ for which--
		\begin{align*}
			&p\left( \xi, \theta_{f_{ \beta'_1 }(\xi)} , \ldots, \theta_{ f_{ \beta'_s }(\xi) } , \mu_{\alpha_{1}}(\xi), \ldots, \mu_{\alpha_k}(\xi) \right)\restriction_{ \mu_{ \alpha_k }(\xi) } \Vdash \\
			&p\left( \xi,\theta_{f_{ \beta'_1 }(\xi)} , \ldots, \theta_{ f_{ \beta'_s }(\xi) } , \mu_{\alpha_{1}}(\xi), \ldots, \mu_{\alpha_k}(\xi) \right)\setminus \mu_{\alpha_k}(\xi) \in \\
			&e\left( \xi, \theta_{f_{ \beta'_1 }(\xi)} , \ldots, \theta_{ f_{ \beta'_s }(\xi) } , \mu_{\alpha_{1}}(\xi), \ldots, \mu_{\alpha_k}(\xi) \right) 
		\end{align*}
		\item There exists a set of $ \xi $-s in $ W $ for which-- 
		$$p\left( \xi, \theta_{f_{ \beta'_1 }(\xi)} , \ldots, \theta_{ f_{ \beta'_s }(\xi) } , \mu_{\alpha_{1}}(\xi), \ldots, \mu_{\alpha_k}(\xi) \right) \in G $$
	\end{enumerate}
\end{theorem}

(Intuitively, for the majority of values of $ \langle \xi, \eta_1, \ldots, \eta_s, \nu_1, \ldots, \nu_k \rangle $,  the condition\\ $p\left( \xi, \eta_1, \ldots , \eta_{s}, \nu_1, \ldots, \nu_k\right)$ which we will construct, forces that--
$$\langle \theta_{f_{ \beta_1 }(\xi)} , \ldots, \theta_{ f_{ \beta_s }(\xi) } , \mu_{\alpha_1}(\xi), \ldots, \mu_{\alpha_k}(\xi)  \rangle = \langle \eta_1, \ldots, \eta_s, \nu_1, \ldots, \nu_k \rangle $$
and its final segment belongs to $ e\left( \xi, \eta_1, \ldots, \eta_s, \nu_1, \ldots, \nu_k \right) $).

\begin{remark}
	When we extend a sequence of generators $ \langle \beta_1, \ldots, \beta_l \rangle $ to a sequence $ \langle \beta'_{1}, \ldots, \beta'_{s} \rangle $
	we will naturally identify the set  $ e\left( \xi, \eta_1, \ldots, \eta_l \right) $, with--
	$$ e'\left( \xi, \eta_1, \ldots, \eta_s \right) = e\left(  \xi, \eta_{i_1}, \ldots, \eta_{i_l} \right) $$
	where $ i_j$ is the index for which $ \beta'_{i_j} =\beta_{j} $, for every $ 1\leq j \leq l $.
	
	Similarly, whenever a function $ g\in V$ is given, whose variables are $ \xi, \eta_1, \ldots, \eta_l, \nu_1, \ldots, \nu_k $, we abuse the notation and denote $ g\left( \xi, \eta_1,\ldots, \eta_{s}, \nu_1, \ldots, \nu_k \right) $ to mean $g\left( \xi, \eta_{i_1},\ldots, \eta_{ i_l }, \nu_1, \ldots, \nu_k  \right)$.
\end{remark}

The proof of theorem \ref{Theorem: Multivariable Fusion Replacement} goes by generalizing the given sets $ e\left( \xi, \eta_1, \ldots, \eta_l, \nu_1, \ldots, \nu_k \right) $:

\begin{definition} \label{Def: the sets e}
For every $ \eta_1, \ldots, \eta_l< \kappa $, $ 1\leq i\leq k $ and an increasing sequence $ \langle \xi, \nu_1, \ldots, \nu_{i} \rangle $, we  define a set $ e\left( \xi, \eta_1, \ldots, \eta_l, \nu_1, \ldots, \nu_{i} \right)\subseteq P\setminus \nu_{i} $. \\
For $ i=k $ this is the set $ e\left( \xi, \eta_1, \ldots, \eta_l, \nu_1, \ldots, \nu_k \right) $ given in the formulation of the theorem.\\
Assume that $ 1\leq i<k $. Work by recursion. Assume that for every $ \nu< g_{i+1}\left( \xi, \eta_1, \ldots, \eta_l, \nu_1, \ldots, \nu_i \right) $, the set $ e\left( \xi, \eta_1, \ldots, \eta_l, \nu_1,\ldots, \nu_{i}, \nu  \right) $ is defined. Denote $ g_{i+1} = g_{i+1}\left( \xi, \eta_1, \ldots, \eta_l, \nu_1, \ldots,\nu_i \right) $. Let us define the set $ e\left( \xi, \eta_1, \ldots, \eta_l, \nu_1, \ldots, \nu_i \right) $, as follows:  A condition $ q\in P\setminus \nu_{i} $ belongs to $ e\left( \xi, \eta_1, \ldots, \eta_l, \nu_1, \ldots, \nu_{i} \right)$ if and only if the following properties hold:
\begin{enumerate}
\item (A technical requirement) $q\restriction_{ g_{i+1} }$ decides the statements-- 
$$ F_{i+1}\left(  \xi, \eta_1, \ldots, \eta_l, \nu_1, \ldots, \nu_j \right) = \lusim{W}_{ g_{i+1}  }\cap V  \ , \ t^{q}_{ g_{i+1} } = t_{i+1}\left(  \xi, \eta_1, \ldots, \eta_l, \nu_1, \ldots, \nu_i \right)$$ Also, if $q\restriction_{ g_{i+1} }$ decides that $  t^{q}_{ g_{i+1} }\neq t_{i+1}\left(  \xi, \eta_1, \ldots, \eta_l, \nu_1, \ldots, \nu_i \right) $, it also decides whether one of the sequences is an initial segment of the other, and if so, which one it is. Finally, if it forces that $ t^{q}_{g_{i+1}} $ is a strict initial segment of $t_{i+1}\left(  \xi, \eta_1, \ldots, \eta_l, \nu_1, \ldots, \nu_i \right)$, it also forces that $ A^{q}_{g_{i+1}} \subseteq g_{i+1}\setminus \max\left( t_{i+1}\left(  \xi, \eta_1, \ldots, \eta_l, \nu_1, \ldots, \nu_i \right) \right) $.
\item (The essential requirement) If both statements in the technical requirement are decided positively, there exists a sequence--
	$$ \langle \  q\left(  \nu \right) \colon \nu< g_{i+1}\left(  \xi, \eta_1, \ldots, \eta_l, \nu_1, \ldots, \nu_i  \right)  \   \rangle $$
	such that, for every $ \nu< g_{i+1}\left(  \xi, \eta_1, \ldots, \eta_l, \nu_1, \ldots, \nu_i  \right)$ above $ \nu_{i} $, $ q(\nu)\in P\setminus \nu $ extends $ q\setminus \nu $, and--
	$$ q\Vdash \mbox{if } \lusim{\mu}_{ \alpha_{i+1} }(\xi) = \nu, \mbox{ then } q(\nu)\in G\setminus \nu  \mbox{ and } q(\nu) \in e\left( \xi, \eta_1, \ldots, \eta_l, \nu_1, \ldots, \nu_{i}, \nu  \right) $$

\end{enumerate}

Similarly, given $ \langle \xi, \eta_1, \ldots, \eta_l \rangle $, define $e\left( \xi, \eta_1, \ldots, \eta_l \right) $ to be the set of conditions $ q\in P\setminus \xi $ which decide whether $ F_1\left(  \xi, \eta_1, \ldots, \eta_l \right) = W_{g_1\left( \xi,\eta_1, \ldots,\eta_l \right)}\cap V $, $ t_1\left( \xi, \eta_1, \ldots, \eta_l \right) = t^{q}_{g_1\left( \xi,\eta_1,\ldots, \eta_l \right)} $, and, assuming that it is decided positively, have a system of extensions--
$$ \langle q\left(  \nu\right) \colon  \nu < g_1\left( \xi, \eta_1, \ldots, \eta_l \right) \rangle $$
such that, for every $ \nu< g_1\left( \xi, \eta_1, \ldots, \eta_l \right)$, $ q(\nu) \in P\setminus \nu $, and--
$$ q\Vdash \mbox{if } \lusim{\mu}_{ \alpha_1 }(\xi) = \nu \mbox{ then } q(\nu)\in G\setminus \nu \mbox{ and } q(\nu)\in e\left( \xi, \eta_1, \ldots, \eta_l, \nu \right) $$
If it is decided negatively, then $ q\restriction_{g_1} $ knows how to compare $ t^{q}_{g_1} $ and $ t_{1}\left( \xi, \eta_1, \ldots, \eta_l \right) $ as in the second point above.

\end{definition}

By induction, we will argue that for every $ i\leq k $ and $ \xi, \eta_1, \ldots, \eta_l, \nu_1, \ldots, \nu_i $, the set $ e\left( \xi, \eta_1, \ldots, \eta_l, \nu_1, \ldots, \nu_i \right)\subseteq P\setminus \nu_i $ is $ \leq^* $-dense open above conditions $ q\in P\setminus \nu_{i} $ for which--
\begin{align*}
	q \Vdash\ & \langle  \theta_{f_{ \beta_1 }(\xi)} , \ldots, \theta_{ f_{ \beta_l }(\xi) } , \lusim{\mu}_{\alpha_1}(\xi), \ldots, \lusim{\mu}_{\alpha_i}(\xi) \rangle = \langle \eta_1, \ldots, \eta_l, \nu_1, \ldots, \nu_{i} \rangle, \mbox{ and for}\\
	&\mbox{every } 1\leq j\leq i, \ F_{j+1}\left(  \xi, \eta_1, \ldots, \eta_l, \nu_1, \ldots, \nu_j \right) = \lusim{W}_{ g_{j+1}\left(   \xi, \eta_1, \ldots, \eta_l, \nu_1, \ldots, \nu_j  \right)  } \mbox{ and }\\
	& t_{j+1}\left(  \xi, \eta_1, \ldots, \eta_l, \nu_1, \ldots, \nu_j \right) = t^{q}_{ g_{j+1}\left( \xi, \eta_1, \ldots, \eta_l, \nu_1, \ldots, \nu_i \right) }
\end{align*}

The induction will be inverse: The basis, for $ i=k $, is true, as it is known that the set $ e\left( \xi, \eta_1, \ldots, \eta_l, \nu_1, \ldots, \nu_k\right)\subseteq P\setminus \nu_k $ is $ \leq^* $ dense--open above conditions $ q\in P\setminus \nu_k $ which force that--
\begin{align*}
	\langle  \theta_{f_{ \beta_1 }(\xi)} , \ldots, \theta_{ f_{ \beta_l }(\xi) } , \lusim{\mu}_{\alpha_1}(\xi), \ldots, \lusim{\mu}_{\alpha_k}(\xi) \rangle = \langle \eta_1, \ldots, \eta_l, \nu_1, \ldots, \nu_{k} \rangle
\end{align*}
The inductive step is given in the following lemma:

\begin{lemma} \label{Lemma: Inductive step for the sets e}
	Fix $ \eta_1, \ldots, \eta_l< \kappa $, $ 1\leq i< k $ and an increasing sequence $ \langle \xi, \nu_1, \ldots, \nu_{i} \rangle $. Denote $ g_{i+1} = g_{i+1}\left( \xi, \eta_1, \ldots, \eta_l, \nu_1, \ldots, \nu_i \right) $. Assume that for every $ \nu_{i+1} \in \left( \nu_i, g_{i+1} \right) $, the set--
	$$ e\left( \xi, \eta_1, \ldots, \eta_l, \nu_1, \ldots, \nu_{i}, \nu_{i+1} \right)\subseteq P\setminus \nu_{i+1} $$
	is $ \leq^* $-dense open above conditions $ q\in P\setminus \nu_{i+1} $ for which--
	\begin{align*}
		q \Vdash\ & \langle  \theta_{f_{ \beta_1 }(\xi)} , \ldots, \theta_{ f_{ \beta_l }(\xi) } , \lusim{\mu}_{\alpha_1}(\xi), \ldots, \lusim{\mu}_{\alpha_i}(\xi), \lusim{\mu}_{\alpha_{i+1}}(\xi) \rangle = \langle \eta_1, \ldots, \eta_l, \nu_1, \ldots, \nu_{i}, \nu \rangle, \mbox{ and for}\\
		&\mbox{every } 1\leq j\leq i+1, \ F_{j+1}\left(  \xi, \eta_1, \ldots, \eta_l, \nu_1, \ldots, \nu_j \right) = \lusim{W}_{ g_{j+1}\left(   \xi, \eta_1, \ldots, \eta_l, \nu_1, \ldots, \nu_j  \right)  } \mbox{ and }\\
		& t_{j+1}\left(  \xi, \eta_1, \ldots, \eta_l, \nu_1, \ldots, \nu_j \right) = t^{q}_{ g_{j+1}\left( \xi, \eta_1, \ldots, \eta_l, \nu_1, \ldots, \nu_i \right) }
	\end{align*}
	then $ e\left( \xi, \eta_1, \ldots, \eta_l, \nu_1, \ldots, \nu_{i} \right) $
	is $ \leq^* $-dense open above conditions $ q\in P\setminus \nu_{i} $ for which--
	\begin{align*}
		q \Vdash\ & \langle  \theta_{f_{ \beta_1 }(\xi)} , \ldots, \theta_{ f_{ \beta_l }(\xi) } , \lusim{\mu}_{\alpha_1}(\xi), \ldots, \lusim{\mu}_{\alpha_i}(\xi) \rangle = \langle \eta_1, \ldots, \eta_l, \nu_1, \ldots, \nu_{i} \rangle, \mbox{ and for}\\
		&\mbox{every } 1\leq j\leq i, \ F_{j+1}\left(  \xi, \eta_1, \ldots, \eta_l, \nu_1, \ldots, \nu_j \right) = \lusim{W}_{ g_{j+1}\left(   \xi, \eta_1, \ldots, \eta_l, \nu_1, \ldots, \nu_j  \right)  } \mbox{ and }\\
		& t_{j+1}\left(  \xi, \eta_1, \ldots, \eta_l, \nu_1, \ldots, \nu_j \right) = t^{q}_{ g_{j+1}\left( \xi, \eta_1, \ldots, \eta_l, \nu_1, \ldots, \nu_i \right) }
	\end{align*}
\end{lemma}

\pr
Let $ q\in P\setminus \nu_{i} $ be a condition which forces that--
\begin{align*}
	&\theta_{f_{ \beta_1 }(\xi)} , \ldots, \theta_{ f_{ \beta_l }(\xi) } , \lusim{\mu}_{1}(\xi), \ldots, \lusim{\mu}_{i}(\xi) \rangle = \langle \eta_1, \ldots, \eta_l, \nu_1, \ldots, \nu_{i} \rangle\\
	&\mbox{and for  every } 1\leq j\leq i, \ F_{j+1}\left(  \xi, \eta_1, \ldots, \eta_l, \nu_1, \ldots, \nu_j \right) = \lusim{W}_{ g_{j+1}\left(   \xi, \eta_1, \ldots, \eta_l, \nu_1, \ldots, \nu_j  \right)  } \\
	&\mbox{and } t_{j+1}\left(  \xi, \eta_1, \ldots, \eta_l, \nu_1, \ldots, \nu_j \right) = t^{q}_{ g_{j+1}\left( \xi, \eta_1, \ldots, \eta_l, \nu_1, \ldots, \nu_i \right) }
\end{align*}
Denote--
$$  g = g_{i+1}\left( \xi, \eta_1, \ldots, \eta_l, \nu_1, \ldots, \nu_{i}   \right) $$
$$ U_g = F_{i+1}\left( \xi, \eta_1, \ldots, \eta_l, \nu_1, \ldots, \nu_{i} \right) $$
$$ t = t_{i+1}\left( \xi, \eta_1, \ldots, \eta_l, \nu_1, \ldots, \nu_{i} \right) $$
Assume that $ q\restriction_{g} $ forces that--
$$ \lusim{W}_{g}\cap V = U_g , \  t = \lusim{t}^{q}_{g} $$
(if not, we are done since $ q \in e\left( \xi, \eta_1, \ldots, \eta_l, \nu_1, \ldots, \nu_i \right) $).
Denote $ n=\mbox{lh}(t) $. We will now apply the following claim:

\begin{claim} \label{Claim: Inductive step For multivar Fusion}
	Assume that $ p\in G $ is a condition, $ n<\omega $ and $ g\in \Delta$ is measurable in $ V $. Assume that $ U_g $ is a normal measure on $ g $ in $ V $, $ t $ is a finite sequence below $ g $ of length $ n $, and--
	$$ p\Vdash  \lusim{t}^{q}_{g}  = t, \ \lusim{W}_{g}\cap V = U_g $$
	For every $ \nu<g $, assume that $ e\left( \nu \right)\subseteq P\setminus \nu $ is a $ P_{\nu} $-name for a subset of $ P\setminus \nu $, which is $ \leq^* $ dense-open above conditions which force that $ \nu $  is the $ \left(n+1\right) $-th element in the Prikry sequence of $ g $. Then there exists a direct extension $ p^* \geq^* p $ and a sequence $ \langle p\left(  \nu \right) \colon \nu<g \rangle $, such that, for every $ \nu< g $, 
	\begin{align*}
		p^*\Vdash &\mbox{if } \nu \mbox{ appears after } t \mbox{ in the Prikry sequence of }g, \mbox{ then } p(\nu)\in \left(G\setminus \nu\right) \cap e(\nu) 	\\
		& \mbox{and }  p^*\restriction_{\nu} \Vdash   p\left( \nu \right) \geq^* {{p^*}\restriction_{\left[ \nu,g \right)}}^\frown {\langle t^\frown \langle \nu \rangle, \lusim{A}^{p^*}_{g}\setminus \nu \rangle}^\frown p^*\setminus \left(g+1\right)  .
	\end{align*}
\end{claim}

\pr
For every $ \nu< g $, consider the set--
\begin{align*}
	d(\nu) = \{  r\in P\restriction_{\left[ \nu, g \right)} \colon&  r\parallel \nu\in \lusim{A}^{p}_{g} \mbox{, and if } r\Vdash \nu\in \lusim{A}^{p}_{ g } \mbox{ then } \\
	&r\Vdash \exists s\geq^* \langle t^{\frown} \langle \nu \rangle, \lusim{A}^{p}_{ g }\setminus \nu \rangle^{\frown} p\setminus \left(g+1\right) , \ r^{\frown} s \in e(\nu) \} 
\end{align*}
Then $ d(\nu)\subseteq P\restriction_{\left[ \nu, g \right)} $ is $ \leq^* $-dense open above $ p\restriction_{\left[ \nu, g \right)} $. Let $ H_g $ be the $ P_{g} $-name, forced by $ p\restriction_{g} $, to be the $ \leq^* $-generic subset of $ j_{U_{g}}\left( P_{g} \right)\setminus g $, for which-- 
$$ \lusim{W}_g =   \left(U_{ g }\right)_{  \lusim{H}_{g} } $$
(such a generic exists since $ W_g $ is simply generated). 
Let $ \lusim{q}\in \mbox{Ult}\left( V, U_g \right) $ be a $ P_{g} $-name, forced by $ p $ to be a condition in $ \left[   \nu\mapsto d(\nu) \right]_{U_{g}} \cap \lusim{H}_{g} $. Let $ \nu\mapsto \lusim{q}(\nu) \in P\restriction_{ \left[  \nu, g \right) } $ be a function in $ V $ such that $ \left[   \nu\mapsto \lusim{q}\left( \nu \right)  \right]_{U_g} = \lusim{q} $. Then we can assume that for a set of $ \nu $-s in $ U_{g} $, 
\begin{equation} \label{Equation large set in lemma}
	p\restriction_{\nu} \Vdash \lusim{q}(\nu) \in d\left( \nu \right) 
\end{equation}
and, by lemma \ref{lem3-1}, $ p\restriction_{g} $ forces that there exists a set $ \lusim{C}\in W_{g} $,  such that for every $ \nu\in C $,
$$ {p\restriction_{\nu}}^{\frown} \lusim{q}(\nu) \in \lusim{G}\restriction_{g} $$
By shrinking $ C $ if necessary, we can assume that every $ \nu\in C $ also satisfies equation (\ref{Equation large set in lemma}).
Now let us define the extension $ p^*\geq^* p $, and, for every $ \nu<g $, the condition $ p(\nu)\in P\setminus \nu $. First, set--
$$ p^*\restriction_{g} = p\restriction_{g} $$
and, in $ V^{P\restriction_{\nu} }$, set--
$$ p(\nu)\restriction_{g} = \lusim{q}(\nu) $$
Work in an arbitrary generic extension for $ {P\restriction_{g}} $, where $ p^*\restriction_{g} $ belongs. For every $ \nu\in C\cap A^{p}_{g} $ (which thus satisfies ${ p\restriction_{\nu}}^{\frown}\lusim{q}(\nu)\in G\restriction_{g} $), there exists $ s(\nu) \in P\setminus g $, $ s\left( \nu \right) \geq^* \langle t^{\frown} \langle \nu \rangle, \lusim{A}^{p}_{ g }\setminus \nu \rangle^{\frown} q\setminus \left(g+1\right)  $, such that  $ {p(\nu)\restriction_{g}}^{\frown} s(\nu) \in e(\nu) $. Set--
$$  p^*\left( g \right) = \langle \lusim{t}^{p}_{g} ,  \lusim{A}^{p}_{g} \cap C \cap  \left(  \triangle_{ \nu<g, \   \nu\in C\cap A^{p}_{g} }  \lusim{A}^{ s(\nu) }_{g}  \right) \rangle  $$
(the definition above is carried in $ V\left[G\restriction_{g}\right] $, so  $ \lusim{C} $ is available there).

Let $ p^*\setminus \left( g+1 \right) =  s\left( \lusim{\nu} \right)  $, where $ \lusim{\nu} $ is the $ \left(n+1\right) $-th element in the Prikry sequence of $ g $. Finally, let--
$$ p\left( \nu \right)\setminus g = \langle t^{\frown} \langle \nu \rangle, A^{p^*}_{g}\setminus \nu \rangle^{\frown} p^*\setminus \left( g+1 \right) $$
where the above definition is possible if $ {p\restriction_{\nu}}^{\frown}p(\nu)\restriction_{g}\Vdash \nu \in 
\lusim{A}^{p^*}_{g} $; if not, let  $p\left( \nu \right)\setminus g$ be arbitrary.

This completes the definition of $ q^*\geq^* q $ and $ \langle p\left( \nu \right) \colon \nu<g \rangle $. Let us prove that for every $ \nu<g $,
\begin{align*}
	p^*\Vdash &\mbox{if } \nu \mbox{ appears after } t \mbox{ in the Prikry sequence of }g, \mbox{ then } p(\nu)\in \left(G\setminus \nu\right) \cap e(\nu) 	\\
	& \mbox{and }  p^*\restriction_{\nu} \Vdash   p\left( \nu \right) \geq^* {{p^*}\restriction_{\left[ \nu,g \right)}}^\frown {\langle t^\frown \langle \nu \rangle, \lusim{A}^{p^*}_{g}\setminus \nu \rangle}^\frown p^*\setminus \left(g+1\right)  .
\end{align*}
Fix $ \nu<g $ and let $ G $ be a generic set for $ P $ which includes $ p^* $, such that, in $ V\left[G\right] $, $ \nu $ appears after $ t $ in the Prikry sequence of $ g $. In particular, $ \nu\in C $ and thus $ q(\nu)\in G\restriction_{\left[ \nu, g \right)} $. By the definition of $ p(\nu) $, and since $ p^*\in G, q(\nu)\in G\restriction_{ \left[ \nu, g \right) } $, it follows that $ p(\nu)\in G\setminus \nu $, as desired. 

$\square \mbox{ of claim }\ref{Claim: Inductive step For multivar Fusion}.$

Apply claim \ref{Claim: Inductive step For multivar Fusion} with respect to the set $ e\left( \xi, \eta_1, \ldots, \eta_l, \nu_1, \ldots, \nu_i, \nu \right) \subseteq P\setminus \nu$ (recall that $ \xi, \eta_1, \ldots, \eta_l, \nu_1, \ldots, \nu_i $ are fixed), and direct extend $ q $ further, to a condition $ q^*\geq^* q $, which has a system of extensions--
$$\langle  q\left( \nu \right) \colon \nu< g \rangle $$
as in the statement if the lemma.

 It follows that,  for every $ \nu<g $,
$$ q^* \Vdash \mbox{if } \lusim{\mu}_{ \alpha_{i+1} }(\xi)= \nu \mbox{ then } q(\nu)\in G\setminus \nu_{i} \mbox{ and } q(\nu)\setminus \nu \in e\left( \xi, \eta_1, \ldots, \eta_l, \nu_1, \ldots, \nu_i, \nu \right)$$
Therefore $ \langle q(\nu) \colon \nu<g \rangle $ witnesses the fact that $ q^*\in e\left( \xi, \eta_1, \ldots, \eta_l, \nu_1, \ldots, \nu_k \right) $.\\
$\square \mbox{ of lemma } \ref{Lemma: Inductive step for the sets e}$.\\

We now proceed towards the proof of theorem \ref{Theorem: Multivariable Fusion Replacement}. We use the same notations as in the formulation of the theorem.

By induction, the following holds: For every $ \xi, \eta_1, \ldots, \eta_l $, the set $ e\left( \xi, \eta_1, \ldots, \eta_l \right)\subseteq P\setminus \xi
$ is $ \leq^* $ dense open above conditions $ q\in P\setminus \xi $ which force that--
$$ \langle  \theta_{f_{ \beta_1 }(\xi)} , \ldots, \theta_{ f_{ \beta_l }(\xi) }   \rangle = \langle \eta_1, \ldots, \eta_l \rangle $$
and that--
$$ F_{1}\left(  \xi, \eta_1, \ldots, \eta_l \right) = \lusim{W}_{ g_{1}\left(   \xi, \eta_1, \ldots, \eta_l  \right)  } \mbox{ and }  t_{1}\left(  \xi, \eta_1, \ldots, \eta_l \right) = t^{q}_{ g_{1}\left( \xi, \eta_1, \ldots, \eta_l \right) }$$

We would like to perform another step, and move from conditions in $ P\setminus \xi $ to conditions in $ P $. This might require extending the sequence generators $ \beta_1, \ldots, \beta_l $. We do this in the following lemma, which concludes the proof of theorem \ref{Theorem: Multivariable Fusion Replacement}.

\begin{lemma} \label{Lemma: Final lemma in the proof of the theorem - adding more generators}
	There exists $ s<\omega $, a sequence of generators $\langle  \beta'_1,\ldots, \beta'_s\rangle $ of $ i $ which extends $ \langle \beta_1, \ldots, \beta_l \rangle $, and a system of conditions--
	$$ \langle p\left( \xi, \eta'_1, \ldots, \eta'_s, \nu_1, \ldots, \nu_k \right) \colon \eta'_1, \ldots, \eta'_s <\kappa, \ \ \xi< \nu_1< \ldots < \nu_k \rangle $$
	(all of them extend the condition $ p\in G $ given in the statement of theorem \ref{Theorem: Multivariable Fusion Replacement}), such that,
	\begin{align*}
		\{ \xi<\kappa \colon & p\left( \xi, \theta_{f_{ \beta'_1 }(\xi)} , \ldots, \theta_{ f_{ \beta'_s }(\xi) } , \mu_{\alpha_1}(\xi), \ldots, \mu_{\alpha_{k}}(\xi) \right) \restriction_{ \mu_{\alpha_k}(\xi) }  \Vdash \\
		& p\left( \xi, \theta_{f_{ \beta'_1 }(\xi)} , \ldots, \theta_{ f_{ \beta'_s }(\xi) } , \mu_{\alpha_1}(\xi), \ldots, \mu_{\alpha_{k}}(\xi) \right)\setminus \mu_{\alpha_k}(\xi)\in  \\
		&e\left( \xi, \theta_{f_{ \beta'_1 }(\xi)} , \ldots, \theta_{ f_{ \beta'_s }(\xi) } , \mu_{\alpha_1}(\xi), \ldots, \mu_{\alpha_k}(\xi) \right)   \mbox{ and- }\\
		&p\left( \xi, \theta_{f_{ \beta'_1 }(\xi)} , \ldots, \theta_{ f_{ \beta'_s }(\xi) } ,  \mu_{\alpha_1}(\xi)  , \ldots, \mu_{\alpha_k}(\xi)\right)   \in G
		\}
	\end{align*}
\end{lemma}

\pr
	Recall that $ W = U_{H} $ is generated from the elementary embedding $ i \colon V\to N $. Let us consider the set--
	$$  i\left(   \langle \xi, \eta_1, \ldots, \eta_l \rangle \mapsto e\left( \xi, \eta_1, \ldots, \eta_l \right) \right)\left( \kappa, \beta_1, \ldots, \beta_l \right)\subseteq i(P)\setminus \kappa $$
	it is $ \leq^* $-dense open in $ i(P)\setminus \kappa $, and thus meets a condition $ r\in H $. Since $ r\in N $, it can be represented using a sequence of generators $ \langle \beta'_1, \ldots, \beta'_s \rangle $, on which we can assume that it contains $ \langle \beta_1, \ldots, \beta_l \rangle $. Let-- 
	$$ \langle \xi, \eta'_1, \ldots, \eta'_s  \rangle \mapsto r\left(   \xi, \eta'_1, \ldots, \eta'_s  \right)\in P\setminus \xi $$
	be a function in $ V $, such that--
	$$ r=  i\left(    \langle \xi, \eta'_1, \ldots, \eta'_s  \rangle \mapsto r\left(   \xi, \eta'_1, \ldots, \eta'_s  \right)  \right) \left(   \kappa, \beta'_1, \ldots, \beta'_s  \right) $$
	Now, for every $ \langle \xi, \eta'_1, \ldots, \eta'_s, \nu_1, \ldots, \nu_k \rangle $, let us define the condition $ p\left(  \xi, \eta'_1, \ldots , \eta'_s, \nu_1, \ldots, \nu_k \right) \in P$. We do this recursively, and define, for every $ 1\leq i\leq k $, a condition $ p\left( \xi, \eta'_1, \ldots, \eta'_s, \nu_1, \ldots, \nu_i \right)\in P $. Simultaneously, we prove that--
			\begin{align*}
	\{ \xi<\kappa \colon & p\left( \xi, \theta_{f_{ \beta'_1 }(\xi)} , \ldots, \theta_{ f_{ \beta'_s }(\xi) } , \mu_{\alpha_1}(\xi), \ldots, \mu_{\alpha_{i}}(\xi) \right) \restriction_{ \mu_{\alpha_i}(\xi) }  \Vdash \\
	& p\left( \xi, \theta_{f_{ \beta'_1 }(\xi)} , \ldots, \theta_{ f_{ \beta'_s }(\xi) } , \mu_{\alpha_1}(\xi), \ldots, \mu_{\alpha_{i}}(\xi) \right)\setminus \mu_{\alpha_i}(\xi)\in  \\
	&e\left( \xi, \theta_{f_{ \beta'_1 }(\xi)} , \ldots, \theta_{ f_{ \beta'_s }(\xi) } , \mu_{\alpha_1}(\xi), \ldots, \mu_{\alpha_i}(\xi) \right)   \mbox{ and- }\\
	&p\left( \xi,\theta_{f_{ \beta'_1 }(\xi)} , \ldots, \theta_{ f_{ \beta'_s }(\xi) } ,  \mu_{\alpha_1}(\xi)  , \ldots, \mu_{\alpha_i}(\xi)\right)   \in G
	\}
\end{align*}
This will complete the proof of the lemma, and thus, the proof of theorem \ref{Theorem: Multivariable Fusion Replacement}. 
\begin{itemize}
		\item First, fix $ \xi, \eta_1, \ldots, \eta_s $, and let us define $ p\left( \xi, \eta_1,\ldots, \eta_s  \right) $.	If $p\restriction_{ \xi } \Vdash r\left( \xi, \eta_1, \ldots, \eta_s \right)\in e\left(\xi, \eta_1, \ldots, \eta_l \right)$, set $ p\left( \xi, \eta_1, \ldots, \eta_s \right) ={ p\restriction_{\xi}}^{\frown} r\left( \xi, \eta_1,\ldots, \eta_s \right) $. Else, let $ p\left( \xi, \eta_1,\ldots, \eta_s \right) $ be an arbitrary condition above $ p $. We argue that--
			\begin{align*}
			\{ \xi<\kappa \colon   &{p\restriction_{ \xi }} \Vdash r\left( \xi, \theta_{f_{ \beta'_1 }(\xi)} , \ldots, \theta_{ f_{ \beta'_s }(\xi) }  \right) \in e\left( \xi,\theta_{f_{ \beta'_1 }(\xi)} , \ldots, \theta_{ f_{ \beta'_s }(\xi) }   \right)  \mbox{ and }\\
			&p\left( \xi, \theta_{f_{ \beta'_1 }(\xi)} , \ldots, \theta_{ f_{ \beta'_s }(\xi) } \right)  \in G
			\}\in W
		\end{align*}
		Recall that $ r\in H $ was defined such that-- 
		$$p\Vdash  r\in i\left(  \langle \xi, \eta_1, \ldots, \eta_l \rangle \mapsto e\left( \xi, \eta_1, \ldots, \eta_l \right) \right)\left( \kappa, \beta_1, \ldots, \beta_l \right) $$ applying the embedding $ k \colon N\to M $ and reflecting down modulo $ W$ gives--
			\begin{align*}
			\{ \xi<\kappa \colon   {p\restriction_{ \xi }} \Vdash r\left( \xi, \theta_{f_{ \beta'_1 }(\xi)} , \ldots, \theta_{ f_{ \beta'_s }(\xi) }  \right)\in e\left( \xi, \theta_{f_{ \beta'_1 }(\xi)} , \ldots, \theta_{ f_{ \beta'_s }(\xi) }   \right) 			\}\in W
		\end{align*}
	Finally, $p\Vdash  r\in H $ and thus $p\Vdash  k(r)\in j_W(G) $, by lemma \ref{lem3-1}. Reflecting this down gives--
		\begin{align*}
		\{ \xi<\kappa \colon  p\left( \xi, \theta_{f_{ \beta'_1 }(\xi)} , \ldots, \theta_{ f_{ \beta'_s }(\xi) } \right)  \in G
		\}\in W
	\end{align*}
		\item Fix $ \xi, \eta'_1, \ldots, \eta'_s, \nu_1  $ and let us define $ p\left(  \xi, \eta'_1, \ldots, \eta'_s, \nu_1 \right) $. Denote $ g_1 = g_1\left( \xi, \eta'_1,\ldots, \eta'_s \right) $. 
		
		If $ p\left(  \xi, \eta'_1, \ldots, \eta'_s \right)\restriction_{ \xi } \Vdash  p\left(  \xi, \eta'_1, \ldots, \eta'_s \right)\setminus \xi \in e\left( \xi, \eta'_1, \ldots, \eta'_s \right)$, then $ p\left(  \xi, \eta'_1, \ldots, \eta'_s \right)\restriction_{ \xi } = p\restriction_{\xi}$ decides the statements--
		$$ F_{1}\left(  \xi, \eta_1, \ldots, \eta_l, \nu_1, \ldots, \nu_j \right) = \lusim{W}_{ g_{1}  }\cap V  \ , \ t^{q}_{ g_{1} } = t_{1}\left(  \xi, \eta_1, \ldots, \eta_l, \nu_1, \ldots, \nu_i \right)$$
		and, if it decides them positively, it forces that there exists a sequence $ \langle q(\nu) \colon \nu<g_1 \rangle $ witnessing this. Define-- 
		$$ p\left( \xi, \eta'_1, \ldots, \eta'_s, \nu_1 \right) =  {p\left( \xi, \eta'_1, \ldots, \eta'_s \right)\restriction_{ \nu_1 }}^{\frown} q(\nu_1) $$
		If $ p\left(  \xi, \eta'_1, \ldots, \eta'_s \right)\restriction_{ \xi } \nVdash  p\left(  \xi, \eta'_1, \ldots, \eta'_s \right)\setminus \xi \in e\left( \xi, \eta'_1, \ldots, \eta'_s \right)$ , or $ p\left(  \xi, \eta'_1, \ldots, \eta'_s \right)\restriction_{ \xi } \Vdash  p\left(  \xi, \eta'_1, \ldots, \eta'_s \right)\setminus \xi \in e\left( \xi, \eta'_1, \ldots, \eta'_s \right)$ but the statements--
		$$ F_{1}\left(  \xi, \eta_1, \ldots, \eta_l, \nu_1, \ldots, \nu_j \right) = \lusim{W}_{ g_{1}  }\cap V  \ , \ t^{q}_{ g_{1} } = t_{1}\left(  \xi, \eta_1, \ldots, \eta_l, \nu_1, \ldots, \nu_i \right)$$
		are decided negatively,  let $ p\left( \xi, \eta'_1,\ldots, \eta'_s ,\nu_1 \right) $ be an arbitrary condition above $  p\left( \xi, \eta'_1,\ldots, \eta'_s \right)  $.\\
		We argue that--
		\begin{align*}
			\{ \xi<\kappa \colon & p\left( \xi, \theta_{f_{ \beta'_1 }(\xi) }, \ldots, \theta_{ f_{ \beta'_s }(\xi)} , \mu_{\alpha_1}(\xi)\right) \restriction_{ \mu_{\alpha_1}(\xi) }  \Vdash \\
			& p\left( \xi, \theta_{f_{ \beta'_1 }(\xi)} , \ldots, \theta_{ f_{ \beta'_s }(\xi) } , \mu_{\alpha_1}(\xi)\right)\setminus \mu_{\alpha_1}(\xi)\in  \\
			&e\left( \xi, \theta_{f_{ \beta'_1 }(\xi)} , \ldots, \theta_{ f_{ \beta'_s }(\xi) } , \mu_{\alpha_1}(\xi) \right)   \mbox{ and- }\\
			&p\left( \xi, \theta_{f_{ \beta'_1 }(\xi)} , \ldots, \theta_{ f_{ \beta'_s }(\xi) } ,  \mu_{\alpha_1}(\xi) \right)   \in G
			\}
		\end{align*}
		First, by the previous point, 
		\begin{align*}
			  \{\xi<\kappa  \colon   & p\left(  \xi, \theta_{f_{ \beta'_1 }(\xi)} , \ldots, \theta_{ f_{ \beta'_s }(\xi) }  \right)\restriction_{ \xi } \Vdash  p\left(  \xi, \theta_{f_{ \beta'_1 }(\xi)} , \ldots, \theta_{ f_{ \beta'_s }(\xi) }  \right)\setminus \xi \in  \\
		&e\left( \xi, \theta_{f_{ \beta'_1 }(\xi)} , \ldots, \theta_{ f_{ \beta'_s }(\xi) }  \right)\}\in W
		\end{align*}
		By the properties of the set $ e\left(  \xi,  \theta_{f_{ \beta'_1 }(\xi)} , \ldots, \theta_{ f_{ \beta'_s }(\xi) }  \right) $, the condition-- 
		$$p\left(  \xi, \theta_{f_{ \beta'_1 }(\xi)} , \ldots, \theta_{ f_{ \beta'_s }(\xi) }  \right)\restriction_{\xi}$$ decides the statements--
		$$  F_{1}\left(  \xi , \theta_{f_{ \beta'_1 }(\xi)} , \ldots, \theta_{ f_{ \beta'_s }(\xi) }  \right) = \lusim{W}_{ g_1 }\cap V $$
		and--
		$$ t^{ p\left(  \xi, \theta_{f_{ \beta'_1 }(\xi)} , \ldots, \theta_{ f_{ \beta'_s }(\xi) }  \right) }_{ g_1} = t_1\left( \xi, \theta_{f_{ \beta'_1 }(\xi)} , \ldots, \theta_{ f_{ \beta'_s }(\xi) }  \right) $$
		\begin{claim} \label{Claim: The statements are decided positively}
			For a set of $ \xi $-s in $ W $, the above statements are decided in a positive way.\\
		Before the proof of the claim, let us proceed with our argument. By the claim and definition \ref{Def: the sets e},  
		$$ p\left(  \xi, \theta_{f_{ \beta'_1 }(\xi)} , \ldots, \theta_{ f_{ \beta'_s }(\xi) } , \mu_{\alpha_1}(\xi) \right) = {p\left(  \xi, \theta_{f_{ \beta'_1 }(\xi)} , \ldots, \theta_{ f_{ \beta'_s }(\xi) }  \right)\restriction_{ \mu_{\alpha_1}(\xi) } }^{\frown} q\left( \mu_{\alpha_1}(\xi) \right) $$
		and, by the properties of the set $ e\left(  \xi,  \theta_{f_{ \beta'_1 }(\xi)} , \ldots, \theta_{ f_{ \beta'_s }(\xi) }  \right) $, the condition-- 
		$$p\left(  \xi, \theta_{f_{ \beta'_1 }(\xi)} , \ldots, \theta_{ f_{ \beta'_s }(\xi) }  \right)$$ 
		forces that-- 
		$$p\left(  \xi, \theta_{f_{ \beta'_1 }(\xi)} , \ldots, \theta_{ f_{ \beta'_s }(\xi) } , \mu_{\alpha_1}(\xi) \right) = {p\left(  \xi, \theta_{f_{ \beta'_1 }(\xi)} , \ldots, \theta_{ f_{ \beta'_s }(\xi) }  \right)\restriction_{ \mu_{\alpha_1}(\xi) } }^{\frown} q\left( \mu_{\alpha_1}(\xi) \right) \in \lusim{G}$$
		and--
		\begin{align*}
		&p\left(  \xi, \theta_{f_{ \beta'_1 }(\xi)} , \ldots, \theta_{ f_{ \beta'_s }(\xi) } , \mu_{\alpha_1}(\xi) \right) \setminus  \mu_{\alpha_1}(\xi) =  q\left( \mu_{\alpha_1}(\xi) \right) \in \\
		& e\left(    \xi,  \theta_{f_{ \beta'_1 }(\xi)} , \ldots, \theta_{ f_{ \beta'_s }(\xi) } , \mu_{\alpha_1}(\xi) \right)
		\end{align*}
		Thus, for a set of $ \xi $-s in $ W $, 
		\begin{align*}
	\{\xi<\kappa  \colon   & p\left(  \xi, \theta_{f_{ \beta'_1 }(\xi)} , \ldots, \theta_{ f_{ \beta'_s }(\xi) }  \right)\restriction_{ \xi } \Vdash  p\left(  \xi, \theta_{f_{ \beta'_1 }(\xi)} , \ldots, \theta_{ f_{ \beta'_s }(\xi) }  \right)\setminus \xi \in  \\
	&e\left( \xi, \theta_{f_{ \beta'_1 }(\xi)} , \ldots, \theta_{ f_{ \beta'_s }(\xi) }  \right)\}\in W
		\end{align*}
		Which finishes the second step. Thus, it remains to prove claim \ref{Claim: The statements are decided positively}:
		\end{claim}
		\pr Let us prove first that--
		$$ \{ \xi<\kappa \colon p\left(  \xi, \theta_{f_{ \beta'_1 }(\xi)} , \ldots, \theta_{ f_{ \beta'_s }(\xi) }  \right)\restriction_{\xi} \Vdash F_{1}\left(  \xi , \theta_{f_{ \beta'_1 }(\xi)} , \ldots, \theta_{ f_{ \beta'_s }(\xi) }  \right) = \lusim{W}_{ g_1 }\cap V   \} $$
		Assume otherwise. Then in $ M\left[j_W(G)\right] $,
		\begin{align*}
		&j_W\left( \langle \xi, \eta_1, \ldots, \eta_s \rangle\mapsto F_1\left( \langle \xi, \eta_1, \ldots, \eta_s \rangle \right)  \right)\left( \kappa, j_{ 0, \alpha }\left( \beta'_1 \right), \ldots, j_{ 0, \alpha }\left( \beta'_s \right) \right) \neq \\
		&\left[  \xi \mapsto  W_{ g_1\left(  \xi, \theta_{f_{ \beta'_1 }(\xi)} , \ldots, \theta_{ f_{ \beta'_s }(\xi) }  \right) }\cap V \right]_W 
		\end{align*}
		but both sides are equal to  $ k_1\left( U_{ \mu_{\alpha_1} }  \right) $, contradicting property \ref{Property:   U mu alpha} of the embedding $ k_{\alpha_1} $. \\
		Now let us prove that--
		\begin{align*}
			 \{ \xi<\kappa \colon &p\left(  \xi, \theta_{f_{ \beta'_1 }(\xi)} , \ldots, \theta_{ f_{ \beta'_s }(\xi) }  \right)\restriction_{\xi} \Vdash \\
			 &t^{ p\left(  \xi, \theta_{f_{ \beta'_1 }(\xi)} , \ldots, \theta_{ f_{ \beta'_s }(\xi) }  \right) }_{ g_1} = t_1\left( \xi, \theta_{f_{ \beta'_1 }(\xi)} , \ldots, \theta_{ f_{ \beta'_s }(\xi) }  \right)  \} 
		 \end{align*}
		Assume otherwise. Then the condition $ s = j_W\left(  \xi \mapsto  p\left(   \xi, \theta_{f_{ \beta'_1 }(\xi)} , \ldots, \theta_{ f_{ \beta'_s }(\xi) }  \right) \right)\left( \kappa \right) $ forces that--
		\begin{align*}
		 t^{s}_{ k_{ \alpha_1 }\left( \mu_{\alpha_1} \right)  } \neq  k_{\alpha_1}\left( t_{ \alpha_1 } \right) = t_{\alpha_1}
		\end{align*}

		Note that $ s\in j_W(G)\restriction_{ k_{ \alpha_1 }\left( \mu_{\alpha_1} \right) }   $ and $ t_{\alpha_1} $ is the initial segment of the Prikry sequence of $ k_{\alpha_1}\left( \mu_{\alpha_1} \right) $ below $ \mu_{\alpha_1} $ in $ M\left[ j_W(G) \right] $. Thus, one of the sequences $t^{s}_{ k_{ \alpha_1 }\left( \mu_{\alpha_1} \right)  } $ and $ t_{\alpha_1} $ is a strict initial segment of the other. By the second requirement in definition \ref{Def: the sets e} $, s\restriction_{ k_{\alpha_1}\left( \mu_{\alpha_1} \right) }   $ decides which one is an initial segment of the other. Now this yields a contradiction:
		\begin{enumerate}
			\item If $ t_{\alpha_1} $ is a strict initial segment of $t^{s}_{ k_{ \alpha_1 }\left( \mu_{\alpha_1} \right)  }$: Recall that $ s = k_{\alpha_1}(s') $, where--
			$$ s' =   j_{\alpha_1}\left(   \langle \xi, \eta_1, \ldots, \eta_s \rangle \mapsto  p\left( \xi, \eta_1 ,\ldots, \eta_s \right)  \right)\left( \kappa, j_{ 0,\alpha_1 }(\beta'_1), \ldots, j_{0, \alpha_1}\left( \beta'_s \right) \right) $$
			Then $ s'\restriction_{ \mu_{\alpha_1} }  $ forces that  $ t_{\alpha_1} $ is a strict initial segment of $ t^{s'}_{\mu_{\alpha_1}} $. Work over $ M_{\alpha_1} $. Let $ \gamma< \mu_{\alpha_1} $ be an ordinal, forced by $ s'\restriction_{ \mu_{\alpha_1}}  $ to be a bound on the first ordinal in  $ t^{s'}_{\mu_{\alpha_1}}\setminus t_{ \alpha_1 } $ (such a bound exists since the forcing $ j_{\alpha_1}(P)\restriction_{ \mu_{\alpha_1} } $ is $ \mu_{\alpha_1} $-c.c. in $ M_{\alpha_1} $). Applying $ k_{\alpha_1} \colon M_{\alpha_1}\to M $, $ \gamma < \mu_{\alpha_1} $ is an upper bound on the first ordinal in $ t^{s}_{ k_{\alpha_1}\left( \mu_{\alpha_1} \right)  }\setminus t_{\alpha_1} $. However, in $ M\left[ j_W(G) \right] $, this element is $ \mu_{\alpha_1}  $ itself, which is strictly above $ \gamma $. A contradiction.
			\item Else, $t^{s}_{ k_{ \alpha_1 }\left( \mu_{\alpha_1} \right)  }$ is a strict initial segment of $t_{\alpha_1}$: Denote $ \gamma = \max\left(t_{\alpha_1}\right) $. Then, by definition \ref{Def: the sets e}, $ s $ forces that the initial segment of the Prikry sequence of $ k_{\alpha_1}\left( \mu_{\alpha_1} \right) $ is $ t^{s}_{  k_{\alpha_1}\left( \mu_{\alpha_1} \right) } $, followed by an element strictly above $ \gamma $; in particular, $ t_{\alpha_1} $ is not an initial segment of the Prikry sequence of $ k_{\alpha_1}\left( \mu_{\alpha_1} \right) $ in $ M\left[ j_W(G) \right] $, which is a contradiction. 
		\end{enumerate}
		$\square$ of claim \ref{Claim: The statements are decided positively}.
			
		\item Assume now that $ 1\leq i< k $ is arbitrary, and for every $ \xi, \eta'_1, \ldots, \eta'_s, \nu_1, \ldots, \nu_{i} $, a condition $p\left( \xi, \eta'_1, \ldots, \eta'_s, \nu_1, \ldots, \nu_{i} \right)$ is defined. Denote $ g_{i+1} = g_{i+1}\left( \xi, \eta'_1, \ldots, \eta'_s, \nu_1, \ldots, \nu_i \right) $. For every $ \nu_{i+1}< g_{i+1} $, let us define the condition $ p\left( \xi,\eta'_1, \ldots, \eta'_s, \nu_1, \ldots, \nu_{i}, \nu_{i+1} \right) $.
		If $ p\left(  \xi, \eta'_1, \ldots, \eta'_s, \nu_1, \ldots, \nu_{i} \right)\restriction_{ \nu_{i} } \Vdash  p\left(  \xi, \eta'_1, \ldots, \eta'_s, \nu_1, \ldots, \nu_{i} \right)\setminus \nu_{i} \in e\left( \xi, \eta'_1, \ldots, \eta'_s, \nu_1, \ldots, \nu_{i} \right)$ and $p\left(  \xi, \eta'_1, \ldots, \eta'_s, \nu_1, \ldots, \nu_{i} \right)\restriction_{ \nu_{i} }$ forces the statements--
						$$ F_{i+1}\left(  \xi, \eta_1, \ldots, \eta_l, \nu_1, \ldots, \nu_i \right) = \lusim{W}_{ g_{i+1}  }\cap V  \ , \ t^{q}_{ g_{i+1} } = t_{i+1}\left(  \xi, \eta_1, \ldots, \eta_l, \nu_1, \ldots, \nu_i \right)$$
		then $ p\left(  \xi, \eta'_1, \ldots, \eta'_s, \nu_1, \ldots, \nu_{i} \right)\restriction_{ \nu_{i} } $ forces that there exists a sequence $ \langle q(\nu) \colon \nu<g_{i+1} \rangle $ witnessing this. In this case, define-- 
		$$ p\left( \xi, \eta'_1, \ldots, \eta'_s, \nu_1, \ldots, \nu_{i}, \nu_{i+1} \right) =  {p\left( \xi, \eta'_1, \ldots, \eta'_s, \nu_1, \ldots, \nu_{i} \right)\restriction_{ \nu_{i+1} }}^{\frown} q(\nu_{i+1}) $$
		Else, let $ p\left( \xi, \eta'_1,\ldots, \eta'_s ,\nu_1, \ldots, \nu_{i}, \nu_{i+1} \right) $ be an arbitrary condition which extends the condition $  p\left( \xi, \eta'_1,\ldots, \eta'_s, \nu_1, \ldots, \nu_{i} \right)  $.
		
		Let us argue now that--
			\begin{align*}
			\{ \xi<\kappa \colon & p\left( \xi, \theta_{f_{ \beta'_1 }(\xi)} , \ldots, \theta_{ f_{ \beta'_s }(\xi) } , \mu_{\alpha_1}(\xi), \ldots, \mu_{\alpha_{i}}(\xi), \mu_{\alpha_{i+1}}(\xi) \right) \restriction_{ \mu_{\alpha_{i+1}}(\xi) }  \Vdash \\
			& p\left( \xi, \theta_{f_{ \beta'_1 }(\xi)} , \ldots, \theta_{ f_{ \beta'_s }(\xi) } , \mu_{\alpha_1}(\xi), \ldots, \mu_{\alpha_{i}}(\xi) , \mu_{\alpha_{i+1}}(\xi)\right)\setminus \mu_{\alpha_{i+1}}(\xi)\in  \\
			&e\left( \xi, \theta_{f_{ \beta'_1 }(\xi)} , \ldots, \theta_{ f_{ \beta'_s }(\xi) } , \mu_{\alpha_1}(\xi), \ldots, \mu_{\alpha_i}(\xi) , \mu_{\alpha_{i+1}}(\xi)\right)   \mbox{ and- }\\
			&p\left( \xi, \theta_{f_{ \beta'_1 }(\xi)} , \ldots, \theta_{ f_{ \beta'_s }(\xi) } ,  \mu_{\alpha_1}(\xi)  , \ldots, \mu_{\alpha_i}(\xi), \mu_{\alpha_{i+1}}(\xi)\right)   \in G
			\}
		\end{align*}
		
		We do this as in the previous point. First,
		\begin{align*}
		\{\xi<\kappa  \colon   & p\left(  \xi, \theta_{f_{ \beta'_1 }(\xi)} , \ldots, \theta_{ f_{ \beta'_s }(\xi) } , \mu_{\alpha_1}(\xi), \ldots , \mu_{\alpha_i}(\xi) \right)\restriction_{ \mu_{\alpha_i}(\xi) } \Vdash \\
		& p\left(  \xi, \theta_{f_{ \beta'_1 }(\xi)} , \ldots, \theta_{ f_{ \beta'_s }(\xi) } , \mu_{\alpha_1}(\xi), \ldots , \mu_{\alpha_i}(\xi) \right)\setminus \mu_{\alpha_i}(\xi) \in  \\
			&e\left( \xi, \theta_{f_{ \beta'_1 }(\xi)} , \ldots, \theta_{ f_{ \beta'_s }(\xi) } , \mu_{\alpha_1}(\xi), \ldots , \mu_{\alpha_i}(\xi) \right)\}\in W
		\end{align*}
		Thus, for a set of $ \xi $-s in $ W $, the condition--
		$$p\left(  \xi, \theta_{f_{ \beta'_1 }(\xi)} , \ldots, \theta_{ f_{ \beta'_s }(\xi) } , \mu_{\alpha_1}(\xi), \ldots , \mu_{\alpha_i}(\xi) \right)\restriction_{ \mu_{\alpha_i}(\xi) }$$
		decides the statements--
		\begin{align*}
		&F_{i+1}\left(  \xi , \theta_{f_{ \beta'_1 }(\xi)} , \ldots, \theta_{ f_{ \beta'_s }(\xi) } , \mu_{\alpha_1}(\xi), \ldots , \mu_{\alpha_i}(\xi) \right) =  \\
		&\lusim{W}_{ g_{i+1}\left(  \xi , \theta_{f_{ \beta'_1 }(\xi)} , \ldots, \theta_{ f_{ \beta'_s }(\xi) } , \mu_{\alpha_1}(\xi), \ldots , \mu_{\alpha_i}(\xi)   \right) }\cap V 
		\end{align*}
		and--
		$$ t^{ p\left(  \xi, \theta_{f_{ \beta'_1 }(\xi)} , \ldots, \theta_{ f_{ \beta'_s }(\xi) }   , \mu_{\alpha_1}(\xi), \ldots , \mu_{\alpha_i}(\xi) \right) }_{ g_{i+1}\left( \xi , \theta_{f_{ \beta'_1 }(\xi)} , \ldots, \theta_{ f_{ \beta'_s }(\xi) } , \mu_{\alpha_1}(\xi), \ldots , \mu_{\alpha_i}(\xi)  \right)  } = t_{i+1}\left( \xi, \theta_{f_{ \beta'_1 }(\xi)} , \ldots, \theta_{ f_{ \beta'_s }(\xi) } , \mu_{\alpha_1}(\xi), \ldots , \mu_{\alpha_i}(\xi) \right) $$
		arguing as in claim \ref{Claim: The statements are decided positively}, both statements are decided positively for a set of $ \xi $-s in $ W $. Thus,
		\begin{align*}
		& p\left(  \xi, \theta_{f_{ \beta'_1 }(\xi)} , \ldots, \theta_{ f_{ \beta'_s }(\xi) }   , \mu_{\alpha_1}(\xi), \ldots , \mu_{\alpha_i}(\xi) , \mu_{\alpha_{i+1}}(\xi)\right) = \\
		&{p\left(  \xi, \theta_{f_{ \beta'_1 }(\xi)} , \ldots, \theta_{ f_{ \beta'_s }(\xi) }   , \mu_{\alpha_1}(\xi), \ldots , \mu_{\alpha_i}(\xi) \right)\restriction_{ \mu_{\alpha_{i+1}  }(\xi) }    }^{\frown} q\left(  \mu_{\alpha_{i+1}}(\xi) \right)
		\end{align*}
		and the condition $q\left(  \mu_{\alpha_{i+1}}(\xi) \right)$ is forced, by--
		$$p\left(  \xi, \theta_{f_{ \beta'_1 }(\xi)} , \ldots, \theta_{ f_{ \beta'_s }(\xi) }   , \mu_{\alpha_1}(\xi), \ldots , \mu_{\alpha_i}(\xi) \right)$$
		to be in-- 
		$$ G\setminus  \mu_{\alpha_{i+1}}(\xi) \cap e\left(   \xi, \theta_{f_{ \beta'_1 }(\xi)} , \ldots, \theta_{ f_{ \beta'_s }(\xi) }   , \mu_{\alpha_1}(\xi), \ldots , \mu_{\alpha_i}(\xi) , \mu_{\alpha_{i+1}}(\xi)  \right) $$  
		Therefore,
		\begin{align*}
			\{ \xi<\kappa \colon & p\left( \xi, \theta_{f_{ \beta'_1 }(\xi)} , \ldots, \theta_{ f_{ \beta'_s }(\xi) } , \mu_{\alpha_1}(\xi), \ldots, \mu_{\alpha_{i}}(\xi), \mu_{\alpha_{i+1}}(\xi) \right) \restriction_{ \mu_{\alpha_{i+1}}(\xi) }  \Vdash \\
			& p\left( \xi, \theta_{f_{ \beta'_1 }(\xi)} , \ldots, \theta_{ f_{ \beta'_s }(\xi) } , \mu_{\alpha_1}(\xi), \ldots, \mu_{\alpha_{i}}(\xi) , \mu_{\alpha_{i+1}}(\xi)\right)\setminus \mu_{\alpha_{i+1}}(\xi)\in  \\
			&e\left( \xi, \theta_{f_{ \beta'_1 }(\xi)} , \ldots, \theta_{ f_{ \beta'_s }(\xi) } , \mu_{\alpha_1}(\xi), \ldots, \mu_{\alpha_i}(\xi) , \mu_{\alpha_{i+1}}(\xi)\right)   \mbox{ and- }\\
			&p\left( \xi, \theta_{f_{ \beta'_1 }(\xi)} , \ldots, \theta_{ f_{ \beta'_s }(\xi) } ,  \mu_{\alpha_1}(\xi)  , \ldots, \mu_{\alpha_i}(\xi), \mu_{\alpha_{i+1}}(\xi)\right)   \in G
			\}
		\end{align*} 
	as desired.
	\end{itemize}
	$\square$ of lemma \ref{Lemma: Final lemma in the proof of the theorem - adding more generators}.
	$ \square $ of theorem \ref{Theorem: Multivariable Fusion Replacement}.

\subsection{Properties of $ k_{\alpha} $} \label{Subsection: Properties of kalpha}

In this subsection we complete the proof of properties $ (A)-(D) $ of $ k_{\alpha} $. After that, we will prove in lemma \ref{Lemma:  k_kappa* is the identity} that $ k_{\kappa^*}\colon M_{\kappa^*}\to M $ is the identity, and conclude the proof of theorems  \ref{Theorem: Restrictions of ultrapwers with simply generated measures} and \ref{Theorem: Measure generated from i restrics to iteration of N}.

\begin{lemma} \label{Lemma: mu alpha is measurable in M alpha}
	$ \mu_{\alpha} = \mbox{crit}\left( k_{\alpha} \right) $ is measurable in $ M_{\alpha} $. Moreover, $\mu_{\alpha}$ is the least measurable above $ \mbox{sup}\{ \mu_{\beta} \colon \beta<\alpha \} $ which has cofinality above $ \kappa $ in $ V $.
\end{lemma}

\pr Write $ \mu = \left[f\right]_W $ and $ \mu = j_{\alpha}\left( h \right)\left(  \kappa,  j_{0,\alpha}(\beta_1), \ldots, j_{0,\alpha}\left( \beta_k \right), \mu_{\alpha_1}, \ldots, \mu_{\alpha_m}   \right)  $, 
for some $ f\in V\left[G\right] $, $ h\in V $, $ \beta_1, \ldots, \beta_l $ generators of $ i $  and $ \alpha_1< \ldots < \alpha_k < \alpha $.

Since $ \mu< k_{\alpha}\left( \mu \right) $, we can assume that for every $ \xi<\kappa $, 
$$  f\left( \xi \right) < h\left(  \xi, \theta_{f_{ \beta_1 }(\xi)} , \ldots, \theta_{ f_{ \beta_l }(\xi) } ,  \mu_{\alpha_1}(\xi), \ldots, \mu_{\alpha_k}(\xi) \right) $$
and let $ p\in G $ be a condition which forces this. Given $ \xi, \eta_1, \ldots ,\eta_l, \nu_1, \ldots, \nu_k $, consider the set--
\begin{align*}
	e\left( \xi, \eta_1, \ldots, \eta_l, \nu_1, \ldots, \nu_k \right) = \{&  r\in P\setminus \nu_k \colon \mbox{ for some bounded subset }  A\subseteq h\left( \xi,  \eta_1, \ldots, \eta_l, \nu_1, \ldots, \nu_k \right),\\
	&r\Vdash \lusim{f}(\xi) \in A  \} 
\end{align*}
Then $ e\left(  \xi, \eta_1, \ldots, \eta_l, \nu_1, \ldots, \nu_k \right) $ is $ \leq^* $-dense open above conditions which extend $ p $ and force that--
$$ \langle \theta_{f_{ \beta_1 }(\xi)} , \ldots, \theta_{ f_{ \beta_l }(\xi) } , \mu_{\alpha_1}(\xi), \ldots, \mu_{\alpha_k}(\xi) \rangle = \langle \eta_1, \ldots, \eta_l, \nu_1, \ldots, \nu_k \rangle $$
By Theorem \ref{Theorem: Multivariable Fusion Replacement}, the   sequence $ \langle \beta_1, \ldots, \beta_l \rangle $ can be extended to a sequence $ \langle \beta'_1, \ldots, \beta'_{s} \rangle $, and $ p $ can be extended to a system of conditions,
$$ \langle p\left( \xi, \eta_1, \ldots, \eta_{s} , \nu_1, \ldots, \nu_k\right) \colon  \xi, \eta_1, \ldots, \eta_{s}<\kappa, \nu_1<\ldots< \nu_k<\kappa \rangle $$
such that, for a set of $ \xi $-s in $ W $,
\begin{align*}
	&p\left( \xi, \theta_{f_{ \beta'_1 }(\xi)} , \ldots, \theta_{ f_{ \beta'_s }(\xi) } , \mu_{\alpha_{1}}(\xi), \ldots, \mu_{\alpha_k}(\xi) \right)\restriction_{ \mu_{ \alpha_k }(\xi) } \Vdash \\
	&p\left( \xi, \theta_{f_{ \beta'_1 }(\xi)} , \ldots, \theta_{ f_{ \beta'_s }(\xi) } , \mu_{\alpha_{1}}(\xi), \ldots, \mu_{\alpha_k}(\xi) \right)\setminus \mu_{\alpha_k}(\xi) \in \\
	&e\left( \xi, \theta_{f_{ \beta'_1 }(\xi)} , \ldots, \theta_{ f_{ \beta'_s }(\xi) } , \mu_{\alpha_{1}}(\xi), \ldots, \mu_{\alpha_k}(\xi) \right) 
\end{align*}
and--
$$p\left( \xi, \theta_{f_{ \beta'_1 }(\xi)} , \ldots, \theta_{ f_{ \beta'_s }(\xi) } , \mu_{\alpha_{1}}(\xi), \ldots, \mu_{\alpha_k}(\xi) \right) \in G $$

Assume now that $ \langle \xi, \eta_1, \ldots, \eta_{s}, \nu_1, \ldots, \nu_k \rangle $ are given, such that--
\begin{align*}
	p\left( \xi, \eta_1, \ldots, \eta_{s} , \nu_1, \ldots, \nu_k \right)\restriction_{ \nu_k} \Vdash & p\left( \xi, \eta_1, \ldots, \eta_{s} , \nu_1, \ldots, \nu_k \right) \setminus \nu_k \in \\
	&e\left(  \xi, \eta_1, \ldots, \eta_{s}, \nu_1, \ldots, \nu_k \right) 
\end{align*}
Let $ \lusim{A} $ be a $ P_{ \nu_k } $-name, forced by $ p\left( \xi, \eta_1, \ldots, \eta_s, \nu_1, \ldots, \nu_k \right)\restriction_{ \nu_k} $ to be a witness to the fact that $p\left( \xi, \vec{\eta}, \vec{\nu}\right)\setminus { \nu_k} \in e\left( \xi, \vec{\eta},\vec{\nu} \right)$. Namely it is a bounded subset of $ h\left( \xi, \vec{\eta}, \vec{\nu} \right) $, and $ p\left( \xi, \vec{\eta}, \vec{\nu} \right)\setminus \nu_k \Vdash \lusim{f}(\xi) \in \lusim{A} $.

Let $ A\left(  \xi, \vec{\eta}, \vec{\nu} \right) $ be the set of ordinals $ \gamma< h\left( \xi, \vec{\eta}, \vec{\nu} \right) $ such that, some $ r\geq p\left( \xi, \vec{\eta}, \vec{\nu} \right)\restriction_{ \nu_k } $ forces that $ \gamma \in \lusim{A} $. Since $ \nu_k < h\left( \xi, \vec{\eta}, \vec{\nu} \right) $, $ A\left( \xi, \vec{\eta}, \vec{\nu} \right) $ is a bounded subset of $ h\left( \xi, \vec{\eta}, \vec{\nu} \right)  $. The function $ \langle  \xi, \vec{\eta}, \vec{\nu} \rangle \mapsto A\left(  \xi, \vec{\eta}, \vec{\nu} \right)$ lies in $ V $. 

By the results of  theorem \ref{Theorem: Multivariable Fusion Replacement}, there exists a set of $ \xi $-s in $ W $ for which--
\begin{align*}
G \ni  &p\left(     \xi, \theta_{f_{ \beta'_1 }(\xi)} , \ldots, \theta_{ f_{ \beta'_s }(\xi) } , \mu_{\alpha_{1}}(\xi), \ldots, \mu_{\alpha_k}(\xi)     \right) \Vdash \\
&\lusim{f}(\xi)\in A\left(  \xi, \theta_{f_{ \beta'_1 }(\xi)} , \ldots, \theta_{ f_{ \beta'_s }(\xi) } , \mu_{\alpha_{1}}(\xi), \ldots, \mu_{\alpha_k}(\xi)  \right)
\end{align*}
Thus, in $ M\left[  j_W(G) \right] $, 
\begin{align*}
	  \left[  f \right]_W \in& \left[  \xi \mapsto A\left(  \xi,  \theta_{f_{ \beta'_1 }(\xi)} , \ldots, \theta_{ f_{ \beta'_s }(\xi) } , \mu_{\alpha_1}(\xi), \ldots, \mu_{\alpha_k}(\xi) \right) \right]_W =  \\
	  & k_{\alpha}\left(  j_\alpha\left(   \langle \xi, \vec{\eta}, \vec{\nu}  \rangle \mapsto A\left( \xi, \vec{\eta}, \vec{\nu} \right) \right)\left(  \kappa, j_{0,\alpha}\left( \beta'_1 \right), \ldots, j_{0,\alpha}\left( \beta'_s \right), \mu_{\alpha_1}, \ldots, \mu_{\alpha_k}  \right) \right)  \subseteq \mbox{Im}\left( k_{\alpha} \right)
\end{align*}
where the last inclusion follows since-- 
$$  j_\alpha\left(   \langle \xi, \vec{\eta}, \vec{\nu}  \rangle \mapsto A\left( \xi, \vec{\eta}, \vec{\nu} \right) \right)\left(  \kappa, j_{0,\alpha}\left( \beta'_1 \right), \ldots, j_{0,\alpha}\left( \beta'_s \right), \mu_{\alpha_1}, \ldots, \mu_{\alpha_k}  \right) $$ 
is a bounded subset of--
$$ \mu_{\alpha} =  j_\alpha\left(   \langle \xi, \vec{\eta}, \vec{\nu}  \rangle \mapsto h\left( \xi, \vec{\eta}, \vec{\nu} \right) \right)\left(  \kappa, j_{0,\alpha}\left( \beta'_1 \right), \ldots, j_{0,\alpha}\left( \beta'_s \right), \mu_{\alpha_1}, \ldots, \mu_{\alpha_k}  \right)$$
which is $ \mbox{crit}\left( k_{\alpha} \right) $. \\
Thus we proved that $ \mu_{\alpha} \in \mbox{Im}\left( k_{\alpha} \right) $, which is a contradiction.
$\square$

\begin{lemma} \label{Lemma: Properties at stage alpha: mu alpha appears in its k }
	$\mu_{\alpha}$ appears in the Prikry sequence added to $ k_{\alpha}\left( \mu_{\alpha} \right) $ in $ M\left[j_W(G)\right] $.	
\end{lemma}

\pr
In $ M\left[H\right] $, denote by $ t^* $ the initial segment of the Prikry sequence of $ k_{\alpha}\left( \mu_{\alpha} \right) $ which consists of all the ordinals below $ \mu_{\alpha} $. Denote by $ n^* $ the length of $ t^* $. Let $ \langle  \xi, \vec{\eta}, \vec{\nu}\rangle\mapsto t^*\left( \xi, \vec{\eta}, \vec{\nu} \right) $ be a function in $ V $ such that--
$$t^* = j_{\alpha}\left(   \langle  \xi, \vec{\eta}, \vec{\nu}\rangle\mapsto t^*\left(  \xi, \vec{\eta}, \vec{\nu} \right) \right)  \left( \kappa,  j_{0,\alpha}\left(\beta_1   \right)   , \ldots, j_{0,\alpha}\left(\beta_l   \right) , \mu_{\alpha_0}, \ldots, \mu_{\alpha_k} \right)  $$
(we assumed here that $ t^* $ can be represented using the same generators as $ \mu_{\alpha} $. If this is not the case, modify the set of generators).

We can assume that for every $ \langle \xi, \vec{\eta}, \vec{\nu}\rangle $, $ t^*\left( \xi, \vec{\eta}, \vec{\nu} \right) $ is a sequence of length $ n^* $.  Since $ k_{\alpha}\left( t^* \right) = t^* $, 
$$ \left[ \xi \mapsto t^*\left(  \xi,    \theta_{f_{ \beta_1 }(\xi)} , \ldots, \theta_{ f_{ \beta_l }(\xi) }   , \mu_{\alpha_1}(\xi) , \ldots, \mu_{ \alpha_k}(\xi) \right) \right]_W = t^* $$
In $ V\left[G\right] $, denote, for every $ \xi<\kappa $,
\begin{align*}
\mu_{\alpha}(\xi) =& \mbox{the } \left(n^*+1\right)\mbox{-th element in the Prikry sequence of }\\ 
&h\left( \vec{\xi},   \theta_{f_{ \beta_1 }(\xi)} , \ldots, \theta_{ f_{ \beta_l }(\xi) }   , \mu_{\alpha_1}(\xi) , \ldots, \mu_{ \alpha_k}(\xi) \right)  
\end{align*}
Clearly $ \left[\xi\mapsto \mu_{\alpha}(\xi)\right]_W \geq \mu_{\alpha} $.

 We argue that equality holds. We will prove that for every $ \eta< \left[  \xi \mapsto \mu_{\alpha}(\xi)\right]_W $, $ \eta < \mu_{\alpha} $. Assume that such $ \eta $ is given, and let $ f\in V\left[G\right] $ be a function such that $ \left[f\right]_W = \eta $. Then we can assume that for every $ \xi<\kappa $,
$$ f(\xi) < \mu_{\alpha}(\xi) $$
and let $ p\in G $ be a condition which forces this. 

For every $\xi, \vec{\eta}, \vec{\nu}$, consider the set--
\begin{align*}
	e\left(  \xi, \vec{\eta}, \vec{\nu} \right) = \{ &r\in P\setminus \nu_k \colon \exists \gamma<h\left(   \xi, \vec{\eta}, \vec{\nu} \right), \ r\Vdash \mbox{ if } t^*\left( \xi, \vec{\eta}, \vec{\nu} \right) \mbox{is an initial segment of the }\\
	&\mbox{Prikry sequence of } h\left(\xi, \vec{\eta}, \vec{\nu} \right), \mbox{ then } \lusim{f}(\xi) < \gamma \} 
\end{align*}
then $ e\left( \vec{\xi}, \vec{\nu}_1, \ldots, \vec{\nu}_k \right) $ is $ \leq^* $ dense open above conditions which force that--
$$ \langle \theta_{f_{ \beta_1 }(\xi)} , \ldots, \theta_{ f_{ \beta_l }(\xi) } , \mu_{\alpha_1}(\xi), \ldots, \mu_{\alpha_k}(\xi) \rangle = \langle \eta_1, \ldots, \eta_l, \nu_1, \ldots, \nu_k \rangle $$
This, since, given a name for an element $ \lusim{f}(\xi) $ which is forced to be strictly below $ \mu_{\alpha}(\xi) $, (which is the element which appears right after $ t^*\left( \xi, \vec{\eta}, \vec{\nu} \right) $ in the Prikry sequence of $ h\left( \xi, \vec{\eta}, \vec{\nu} \right) $), the element can be decided by taking a direct extension. \\
By Theorem \ref{Theorem: Multivariable Fusion Replacement}, the  sequence $ \langle \beta_1, \ldots, \beta_l \rangle $ can be extended to a sequence $ \langle \beta'_1, \ldots, \beta'_{s} \rangle $, and $ p $ can be extended to a system of conditions,
$$ \langle p\left( \xi, \eta_1, \ldots, \eta_{s} , \nu_1, \ldots, \nu_k\right) \colon  \xi, \eta_1, \ldots, \eta_{s}<\kappa, \nu_1<\ldots< \nu_k<\kappa \rangle $$
such that, for a set of $ \xi $-s in $ W $,
\begin{align*}
	&p\left( \xi, \theta_{f_{ \beta'_1 }(\xi)} , \ldots, \theta_{ f_{ \beta'_s }(\xi) } , \mu_{\alpha_{1}}(\xi), \ldots, \mu_{\alpha_k}(\xi) \right)\restriction_{ \mu_{ \alpha_k }(\xi) } \Vdash \\
	&p\left( \xi, \theta_{f_{ \beta'_1 }(\xi)} , \ldots, \theta_{ f_{ \beta'_s }(\xi) } , \mu_{\alpha_{1}}(\xi), \ldots, \mu_{\alpha_k}(\xi) \right)\setminus \mu_{\alpha_k}(\xi) \in \\
	&e\left( \xi, \theta_{f_{ \beta'_1 }(\xi)} , \ldots, \theta_{ f_{ \beta'_s }(\xi) } , \mu_{\alpha_{1}}(\xi), \ldots, \mu_{\alpha_k}(\xi) \right) 
\end{align*}
and--
$$p\left( \xi, \theta_{f_{ \beta'_1 }(\xi)} , \ldots, \theta_{ f_{ \beta'_s }(\xi) } , \mu_{\alpha_{1}}(\xi), \ldots, \mu_{\alpha_k}(\xi) \right) \in G $$
Assume now that $ \langle \xi, \vec{\eta}, \vec{\nu} \rangle =  \langle \xi, \eta_1, \ldots, \eta_{s}, \nu_1, \ldots, \nu_k \rangle $ are given, such that--
\begin{align*}
	p\left( \xi, \vec{\eta}, \vec{\nu}\right)\restriction_{ \nu_k} \Vdash & p\left( \xi, \vec{\eta}, \vec{\nu} \right) \setminus \nu_k \in  e\left(  \xi, \vec{\eta}, \vec{\nu} \right) 
\end{align*}
Let $ \lusim{\gamma} $ be a $ P_{ \nu_k } $-name, forced by $ p\left( \xi, \eta_1, \ldots, \eta_s, \nu_1, \ldots, \nu_k \right)\restriction_{ \nu_k} $ to an ordinal below $ h\left( \xi, \vec{\eta}, \vec{\nu} \right) $, such that $ p\left( \xi, \vec{\eta}, \vec{\nu} \right)\setminus \nu_k \Vdash \lusim{f}(\xi) < \lusim{\gamma} $.
Let $ \gamma\left(  \xi, \vec{\eta}, \vec{\nu} \right) $ be the supremum of the set of ordinals $ \tau< h\left( \xi, \vec{\eta}, \vec{\nu} \right) $ such that, some $ r\geq p\left( \xi, \vec{\eta}, \vec{\nu} \right)\restriction_{ \nu_k } $ forces that $\lusim{\gamma} = \tau  $. Since $ \nu_k < h\left( \xi, \vec{\eta}, \vec{\nu} \right) $, $ \gamma\left( \xi, \vec{\eta}, \vec{\nu} \right)  <  h\left( \xi, \vec{\eta}, \vec{\nu} \right)  $. The function $ \langle  \xi, \vec{\eta}, \vec{\nu} \rangle \mapsto \gamma\left(  \xi, \vec{\eta}, \vec{\nu} \right)$ lies in $ V $.

By the results of  theorem \ref{Theorem: Multivariable Fusion Replacement}, there exists a set of $ \xi $-s in $ W $ for which--
\begin{align*}
	G \ni  &p\left(     \xi, \theta_{f_{ \beta'_1 }(\xi)} , \ldots, \theta_{ f_{ \beta'_s }(\xi) } , \mu_{\alpha_{1}}(\xi), \ldots, \mu_{\alpha_k}(\xi)     \right) \Vdash \\
	&\mbox{if } t^*\left(  \xi, \theta_{f_{ \beta'_1 }(\xi)} , \ldots, \theta_{ f_{ \beta'_s }(\xi) } , \mu_{\alpha_{1}}(\xi), \ldots, \mu_{\alpha_k}(\xi)    \right) \\
	&\mbox{is an initial segment of the Prikry sequence of } \\
	&h\left( \xi, \theta_{f_{ \beta'_1 }(\xi)} , \ldots, \theta_{ f_{ \beta'_s }(\xi) } , \mu_{\alpha_{1}}(\xi), \ldots, \mu_{\alpha_k}(\xi)    \right), \mbox{ then }\\
	&\lusim{f}(\xi)< \gamma\left(  \xi, \theta_{f_{ \beta'_1 }(\xi)} , \ldots, \theta_{ f_{ \beta'_s }(\xi) } , \mu_{\alpha_{1}}(\xi), \ldots, \mu_{\alpha_k}(\xi)  \right)
\end{align*}
Thus, in $ M\left[  j_W(G) \right] $, where indeed $ t^*  $ is an initial segment of the Prikry sequence of $ k_{\alpha}\left( \mu_{\alpha} \right) $,
\begin{align*}
	\left[  f \right]_W \in& \left[  \xi \mapsto \gamma\left(  \xi,  \theta_{f_{ \beta'_1 }(\xi)} , \ldots, \theta_{ f_{ \beta'_s }(\xi) } , \mu_{\alpha_1}(\xi), \ldots, \mu_{\alpha_k}(\xi) \right) \right]_W =  \\
	& k_{\alpha}\left(  j_\alpha\left(   \langle \xi, \vec{\eta}, \vec{\nu}  \rangle \mapsto \gamma\left( \xi, \vec{\eta}, \vec{\nu} \right) \right)\left(  \kappa, j_{0,\alpha}\left( \beta'_1 \right), \ldots, j_{0,\alpha}\left( \beta'_s \right), \mu_{\alpha_1}, \ldots, \mu_{\alpha_k}  \right) \right) < \mu_{\alpha}
\end{align*}
as desired.
$\square$

\begin{lemma}
Let $ U_{ \mu_{\alpha} }  = \{  X\subseteq \mu_{ \alpha} \colon \mu_{\alpha} \in  k_{\alpha}(X)  \}\cap M_{\alpha}$.  Then $ U_{ \mu_{\alpha} } \in M_{\alpha} $. Furthermore, 
	$  k_{\alpha}\left( U_{ \mu_{\alpha} } \right) = j_{W}\left(  \delta\mapsto U_{\delta} \right)\left( k_{\alpha}\left( \mu_{\alpha} \right) \right) $, 
	where, for every $ \delta\in \Delta $, $ U_{\delta} = W_{\delta}\cap V $, for $ W_\delta $ which is the measure used in the Prikry forcing at stage $ \delta $ in the iteration $ P $.
\end{lemma}

\pr We first prove that $ j_{W}\left(  \delta\mapsto U_{\delta} \right)\left( k_{\alpha}\left( \mu_{\alpha} \right) \right) \in \mbox{Im}\left( k_{\alpha} \right) $. Then, we will prove that the measure $ F\in M_{\alpha} $ for which $ j_{W}\left(  \delta\mapsto U_{\delta} \right)\left( k_{\alpha}\left( \mu_{\alpha} \right) \right) = k_{\alpha}(F) $ equals to $ U_{\mu_{\alpha}} $. 

In order to prove that $ j_{W}\left(  \delta\mapsto U_{\delta} \right)\left( k_{\alpha}\left( \mu_{\alpha} \right) \right)\in \mbox{Im}(k_{\alpha}) $, we prove that there exists a family $ \mathcal{F}\in M_{\alpha} $ of measures on $ \mu_{\alpha} $, with $ \left| \mathcal{F} \right|< \mu_{\alpha} $, such that $ j_{W}\left(  \delta\mapsto U_{\delta} \right)\left( k_{\alpha}\left( \mu_{\alpha} \right) \right)\in k_{\alpha}(F) = k_{\alpha}''\mathcal{F} $.

Fix, in $ V $, an enumeration $ W $ of all the normal measures on measurable cardinals below $ \kappa $. For every $ \langle \xi, \vec{\eta}, \vec{\nu} \rangle $, let $ \gamma\left( \xi, \vec{\eta}, \vec{\nu} \right) $ be the index of $ U_{ h\left( \xi, \vec{\eta}, \vec{\nu} \right) } $ in this enumeration. Note that each measure $  U_{ h\left( \xi, \vec{\eta}, \vec{\nu} \right) } $  belongs to $ V $, but the sequence $ \langle  U_{ h\left( \xi, \vec{\eta}, \vec{\nu} \right) }   \colon  \xi, \vec{\eta}, \vec{\nu}< \kappa \rangle $ might be external to $ V $. So the function $ \langle \xi, \vec{\eta}, \vec{\nu} \rangle \mapsto \gamma\left( \xi, \vec{\eta}, \vec{\nu} \right) $ doesn't necessarily belong to $ V $. 

Fix $ \langle \xi, \vec{\eta}, \vec{\nu} \rangle $ and consider the set--
	\begin{align*}
	e\left( \xi, \vec{\eta}, \vec{\nu} \right) = &\{ r\in P\setminus \nu_k \colon \mbox{there exists a set of ordinals } A \mbox{ of cardinality strictly smaller than }\\
	&h\left(\xi, \vec{\eta}, \vec{\nu} \right), \mbox{ such that }  r\restriction_{  h\left( \xi, \vec{\eta}, \vec{\nu}\right) } \Vdash \gamma\left( \xi, \vec{\eta}, \vec{\nu} \right) \in A \}
\end{align*}
Then $ e\left( \xi, \vec{\eta}, \vec{\nu} \right)\subseteq P\setminus \nu_k $ is $ \leq^* $-dense open, since $ P\restriction_{  h\left( \xi, \vec{\eta}, \vec{\nu} \right)  } $ is  $ h\left( \xi, \vec{\eta}, \nu \right) $-c.c.. 

Now apply theorem \ref{Theorem: Multivariable Fusion Replacement} and argue as in the previous lemma: There exists (in $ V $) a mapping $ \langle \xi, \vec{\eta}, \vec{\nu} \rangle\mapsto A\left( \xi, \vec{\eta}, \vec{\nu} \right) $ such that, in $ M\left[ j_W(G) \right] $,
\begin{align*}
&\left[  \xi \mapsto \gamma\left(   \xi,  \theta_{f_{ \beta'_1 }(\xi)} , \ldots, \theta_{ f_{ \beta'_s }(\xi) } , \mu_{\alpha_1}(\xi), \ldots, \mu_{\alpha_k}(\xi) \right)   \right]_W   \in  \\
&\left[  \xi \mapsto A\left(   \xi,  \theta_{f_{ \beta'_1 }(\xi)} , \ldots, \theta_{ f_{ \beta'_s }(\xi) } , \mu_{\alpha_1}(\xi), \ldots, \mu_{\alpha_k}(\xi) \right)   \right]_W = \\
& k_{\alpha}''\left(  j_{\alpha}\left(  \langle \xi, \vec{\eta}, \vec{\nu} \rangle\mapsto A\left(  \xi, \vec{\eta}, \vec{\nu} \right) \right)\left(   \kappa, j_{0,\alpha}\left( \beta'_1 \right), \ldots, j_{0,\alpha}\left( \beta'_s \right), \mu_{\alpha_1}, \ldots, \mu_{\alpha_k} \right)  \right)
\end{align*}
In $ M_{\alpha} $, let $ \mathcal{F} $ be the set of measures on $ \mu_{\alpha} $  which are indexed in the enumeration $ j_{\alpha}(W) $ by an index in the set $ \mathcal{A} =  j_{\alpha}\left(  \langle \xi, \vec{\eta}, \vec{\nu} \rangle\mapsto A\left(  \xi, \vec{\eta}, \vec{\nu} \right) \right)\left(   \kappa, j_{0,\alpha}\left( \beta'_1 \right), \ldots, j_{0,\alpha}\left( \beta'_s \right), \mu_{\alpha_1}, \ldots, \mu_{\alpha_k} \right)  $. Note that $ \left| \mathcal{A} \right| < \mu_{\alpha} $ and thus $ \left| \mathcal{F} \right| < \mu_{\alpha} $. Then $ j_W\left( \delta\mapsto U_{\delta} \right)\left( k_{\alpha}\left( \mu_{\alpha} \right) \right) $ is enumerated by the ordinal--
$$ \left[  \xi \mapsto \gamma\left(   \xi,  \theta_{f_{ \beta'_1 }(\xi)} , \ldots, \theta_{ f_{ \beta'_s }(\xi) } , \mu_{\alpha_1}(\xi), \ldots, \mu_{\alpha_k}(\xi) \right)   \right]_W \in k_{\alpha}''\mathcal{A}   $$
and thus $  j_W\left( \delta\mapsto U_{\delta} \right)\left( k_{\alpha}\left( \mu_{\alpha} \right) \right)\in k_{\alpha}''\mathcal{F} $ , as desired.

Let $ F\in M_{\alpha} $ be a measure on $ \mu_{\alpha} $ such that-- 
$$ j_W\left( \delta\mapsto U_{\delta} \right)\left( k_{\alpha}\left( \mu_{\alpha} \right) \right) = k_{\alpha}(F) $$
Let us argue that $ F = U_{ \mu_{\alpha} } $. It suffices to prove that $ F\subseteq U_{ \mu_{\alpha} } $. Fix a set $ X\in F $. Assume that-- 
$$ X  = j_{\alpha}\left( \langle \xi, \vec{\eta}, \vec{\nu} \rangle \mapsto X\left( \xi, \vec{\eta}, \vec{\nu} \right) \right)\left(   \kappa, j_{0,\alpha}\left( \beta_1 \right), \ldots, j_{0,\alpha}\left( \beta_l \right), \mu_{\alpha_1}, \ldots, \mu_{\alpha_k}   \right)$$
(We assumed again that $ X $ can be represented using the same generators as $ \mu_{\alpha} $. If this is not the case, modify the set of generators of $ \mu_{\alpha} $ ). Then $ k_{\alpha}(X) \in j_W\left( \delta\mapsto U_\delta \right)\left(  k_{\alpha}\left( \mu_{\alpha} \right)  \right) $.

As in the previous lemma, let $ n^* $ be the length of $ t^* $, the initial segment of the Prikry sequence of $ k_{\alpha}\left( \mu_{\alpha} \right) $ below $ \mu_{\alpha} $.
For every $\langle \xi, \vec{\eta}, \vec{\nu}\rangle $, let--
\begin{align*}
	e\left(\xi, \vec{\eta}, \vec{\nu}\right) = \{ & r\in P\setminus \nu_k \colon r\restriction_{ h\left( \xi, \vec{\eta}, \vec{\nu} \right) }\parallel X\left( \xi, \vec{\eta}, \vec{\nu}\right)\in U_{ h\left(\xi, \vec{\eta}, \vec{\nu} \right) } ,   \\
	& \mbox{if it decides positively, then } r\restriction_{ h\left(\xi, \vec{\eta}, \vec{\nu} \right) } \Vdash \lusim{A}^{r}_{ h\left( \xi, \vec{\eta}, \vec{\nu}\right) } \subseteq \\
	& X\left( \xi, \vec{\eta}, \vec{\nu} \right) ; \mbox{ else, } r\restriction_{ h\left( \xi, \vec{\eta}, \vec{\nu}\right) } \Vdash \lusim{A}^{r}_{ h\left( \xi, \vec{\eta}, \vec{\nu} \right) } \mbox{ is disjoint} \\
	&\mbox{from } X\left( \xi, \vec{\eta}, \vec{\nu} \right). \mbox{ Moreover, } r\restriction_{ h\left( \xi, \vec{\eta}, \vec{\nu}\right) } \parallel \mbox{lh}\left( t^{r}_{ h\left( \xi, \vec{\eta}, \vec{\nu} \right) } \right) > n^*, \\
	&\mbox{and if it decides positively, then there exists a bounded subset } \\
	& A\left( \xi, \vec{\eta}, \vec{\nu} \right)\subseteq h\left(\xi, \vec{\eta}, \vec{\nu} \right) \mbox{ for which } r\restriction_{\xi, \vec{\eta}, \vec{\nu} } \Vdash \mbox{ the } \left(n^*+1\right)\mbox{-th}\\
	&\mbox{element of }t^{r}_{ h\left( \xi, \vec{\eta}, \vec{\nu} \right) } \mbox{ belongs to } A\left( \xi, \vec{\eta}, \vec{\nu} \right)  \} 
\end{align*}
By theorem \ref{Theorem: Multivariable Fusion Replacement}, there exists a larger set of generators $ \beta'_1, \ldots, \beta'_s $ and, for every $ \langle \xi, \vec{\eta}, \vec{\nu} \rangle $, a condition $ p\left( \langle \xi, \vec{\eta}, \vec{\nu} \rangle \right) $, such that, for a set of $ \xi $-s in $ W $,
\begin{align*}
	&p\left( \xi, \theta_{f_{ \beta'_1 }(\xi)} , \ldots, \theta_{ f_{ \beta'_s }(\xi) } , \mu_{\alpha_{1}}(\xi), \ldots, \mu_{\alpha_k}(\xi) \right)\restriction_{ \mu_{ \alpha_k }(\xi) } \Vdash \\
	&p\left( \xi, \theta_{f_{ \beta'_1 }(\xi)} , \ldots, \theta_{ f_{ \beta'_s }(\xi) } , \mu_{\alpha_{1}}(\xi), \ldots, \mu_{\alpha_k}(\xi) \right)\setminus \mu_{\alpha_k}(\xi) \in \\
	&e\left( \xi, \theta_{f_{ \beta'_1 }(\xi)} , \ldots, \theta_{ f_{ \beta'_s }(\xi) } , \mu_{\alpha_{1}}(\xi), \ldots, \mu_{\alpha_k}(\xi) \right) 
\end{align*}
and--
$$p\left( \xi, \theta_{f_{ \beta'_1 }(\xi)} , \ldots, \theta_{ f_{ \beta'_s }(\xi) } , \mu_{\alpha_{1}}(\xi), \ldots, \mu_{\alpha_k}(\xi) \right) \in G $$
Let us argue first that for a set of $ \xi $-s in $ W $, 
$$ p\left( \xi, \theta_{f_{ \beta'_1 }(\xi)} , \ldots, \theta_{ f_{ \beta'_s }(\xi) } , \mu_{\alpha_{1}}(\xi), \ldots, \mu_{\alpha_k}(\xi) \right)\restriction_{ \mu_{ \alpha_k }(\xi) } $$ decides that--  
$$\mbox{lh}\left( t^{p\left( \xi, \theta_{f_{ \beta'_1 }(\xi)} , \ldots, \theta_{ f_{ \beta'_s }(\xi) } , \mu_{\alpha_{1}}(\xi), \ldots, \mu_{\alpha_k}(\xi) \right)}_{ h\left(\xi, \theta_{f_{ \beta'_1 }(\xi)} , \ldots, \theta_{ f_{ \beta'_s }(\xi) } , \mu_{\alpha_{1}}(\xi), \ldots, \mu_{\alpha_k}(\xi) \right) } \right) \leq n^*$$
Indeed, assume otherwise. Let $ A^*\left( \xi, \vec{\eta}, \vec{\nu} \right) $ be the bounded subset of $ h\left( \xi, \vec{\eta}, \vec{\nu} \right) $  which consists of all the ordinals, which are forced by some extension of $ p\left( \xi, \vec{\eta}, \vec{\nu} \right)\restriction_{  \nu_k} $ to be in $ A\left( \xi, \vec{\eta}, \vec{\nu} \right) $  (whenever $ p\left( \xi, \vec{\eta}, \vec{\nu} \right) $ forces that the length of $ t^{p \left( \xi, \vec{\eta}, \vec{\nu} \right) }_{ h\left( \xi, \vec{\eta}, \vec{\nu} \right) }   $ is greater than $ n^* $). Then, in $ M\left[ j_W(G) \right] $, 
$$ \mu_{\alpha} \in k_{\alpha}\left(  j_{\alpha}\left( \langle \xi, \vec{\eta}, \vec{\nu} \rangle \mapsto A^*\left( \xi, \vec{\eta}, \vec{\nu} \right)   \right)\left(  \kappa, j_{0,\alpha}\left( \beta'_1 \right), \ldots, j_{0,\alpha}\left( \beta'_s \right), \mu_{\alpha_1}, \ldots, \mu_{\alpha_k} \right)   \right) $$
But this is a contradiction, since $ j_{\alpha}\left( \langle \xi, \vec{\eta}, \vec{\nu} \rangle \mapsto A^*\left( \xi, \vec{\eta}, \vec{\nu} \right)   \right)\left(  \kappa, j_{0,\alpha}\left( \beta'_1 \right), \ldots, j_{0,\alpha}\left( \beta'_s \right), \mu_{\alpha_1}, \ldots, \mu_{\alpha_k} \right)   $ is a bounded subset of $ \mu_{\alpha} $. 

Therefore, we can assume that-- 
$$ p\left( \xi, \theta_{f_{ \beta'_1 }(\xi)} , \ldots, \theta_{ f_{ \beta'_s }(\xi) } , \mu_{\alpha_{1}}(\xi), \ldots, \mu_{\alpha_k}(\xi) \right)\restriction_{ \mu_{ \alpha_k }(\xi) } $$ forces that--  
$$\mbox{lh}\left( t^{p\left( \xi, \theta_{f_{ \beta'_1 }(\xi)} , \ldots, \theta_{ f_{ \beta'_s }(\xi) } , \mu_{\alpha_{1}}(\xi), \ldots, \mu_{\alpha_k}(\xi) \right)}_{ h\left(\xi, \theta_{f_{ \beta'_1 }(\xi)} , \ldots, \theta_{ f_{ \beta'_s }(\xi) } , \mu_{\alpha_{1}}(\xi), \ldots, \mu_{\alpha_k}(\xi) \right) } \right) \leq n^*$$
Denote now $ p^* =  \left[  \xi \mapsto p\left( \xi, \theta_{f_{ \beta'_1 }(\xi)} , \ldots, \theta_{ f_{ \beta'_s }(\xi) } , \mu_{\alpha_{1}}(\xi), \ldots, \mu_{\alpha_k}(\xi) \right) \right]_W $. Then $ p^*\restriction_{ k_{\alpha}(\mu_{\alpha}) } $ forces that $ \mu_{\alpha} \in \lusim{A}^{p^*}_{ k_{\alpha}\left( \mu_{\alpha} \right) } $. By the definition of the sets $ e\left( \xi, \vec{\eta}, \vec{\nu} \right) $, the set $ \lusim{A}^{p^*}_{ k_{\alpha}\left( \mu_{\alpha} \right) } $ is forced to be  either disjoint or contained in $ k_{\alpha}(X) $. Since $ k_{\alpha}(X)\in j_W\left(  \delta\mapsto U_{ \delta} \right)\left( k_{\alpha}\left( \mu_{\alpha} \right) \right) $, it cannot be disjoint (again, by the definition of $ e\left( \xi, \vec{\eta}, \vec{\nu} \right) $). Therefore $ \mu_{\alpha} \in k_{\alpha}(X) $ and thus $ X\in U_{\mu_{\alpha} } $, as desired.
$\square$.

Finally, let us argue that $ j_{ \kappa^*} = j_W\restriction_{V} $. Recall that $ \kappa^* = i(\kappa) $, and note that $ \kappa^* = \mbox{sup}\{ \mu_{\alpha} \colon \alpha<\kappa^* \} $.

\begin{lemma} \label{Lemma:  k_kappa* is the identity}
	$ M = M_{\kappa^*} $, $ j_W(\kappa) = i(\kappa)  $ and $ j_{ \kappa^*} = j_W\restriction_{V} $. 
\end{lemma}

\begin{remark}
	In particular, if $ i = j_U $ (namely $ W $ is simply generated) then $ j_W(\kappa) = j_{U}(\kappa) $. On the other hand, possibly $ j_U(\kappa) < i(\kappa) $, and then $ j_W(\kappa) > j_U(\kappa) $. 
\end{remark}

\pr
Define, similarly to $ k_{\alpha} \colon M_{\alpha}\to M $, the embedding $ k_{\kappa^*} \colon M_{\kappa^*} \to M $ as follows: 
\begin{align*}
&k_{\kappa^*}\left(   j_{\kappa^*}\left( f \right)\left(  \kappa,  j_{0,\kappa^*}(\beta_1), \ldots, j_{0,\kappa^*}\left( \beta_l \right), \mu_{\alpha_1}, \ldots, \mu_{\alpha_k}   \right)  \right)  =   \\
&j_W\left( f \right)\left(  \kappa, \theta_{\left[  f_{\beta_1}(\xi) \right]_W} , \ldots,   \theta_{\left[ f_{ \beta_l }(\xi)\right]_W} , \mu_{\alpha_1}, \ldots, \mu_{\alpha_m}  \right) 
\end{align*}
for every $ f\in V $, $ \beta_1, \ldots, \beta_l $ generators of $ i $  and $ \alpha_1< \ldots < \alpha_m < \kappa^* $. Clearly $ \mbox{crit}\left( k_{ \kappa^*} \right) \geq \kappa^* $. It suffices to prove that $ k_{\kappa^* } $ is the identity function.

Let $ \tau $ be an ordinal, and let $ f\in V\left[G\right] $ be a function such that $ \left[f\right]_W = \tau $. By the $ \kappa $-c.c. of $ P_{\kappa} $, there exists $ F\in V $ such that for every $ \xi<\kappa $,
$ f(\xi)\in F(\xi)$ and $ \left| F(\xi) \right| < \kappa $. Therefore, in $ M\left[j_W(G)\right] $,
$$ \tau = \left[ f \right]_W \in \left[ F \right]_W = k_{\kappa^*}\left(  j_{\kappa^*} (F)(\kappa) \right)$$
But--
$$ \left|   j_{\kappa^*} (F)(\kappa) \right| < j_{\kappa^*}\left( \kappa \right) = \kappa^* \leq \mbox{crit}\left( k_{\kappa^*} \right) $$
so $ \eta\in \mbox{Im}\left( k_{\kappa^*} \right) $ as desired.
$\square$

\section{ Further directions and open problems}

It is likely that results of Section 4 can be extended to wider context of Prikry type forcing notions. 
The first candidates are one element Prikry forcings and Prikry forcings with non-normal ultrafilters.
For the former it seems that the present arguments can be applied without much changes. 
The latter looks to require more work since $[id]$ is not $\kappa$ anymore and additional generators may appear.
Another example is Extender based Prikry forcings. Here some new ideas seems to be needed due to the Cohen parts of the forcings.

Let us state some open questions.

\begin{question}
 Are there other ways to generate normal ultrafilters $W$ in $V[G]$ beyond those given in \ref{Theorem: Conditions imposed on i}?
\end{question}

Let $ W\in V\left[G\right] $ be a normal measure on $ \kappa $. Assuming $\neg o^{\P}$ (or even no inner model  with a Woodin cardinal) and exploring closure of the ultrapower, it seems possible to argue that $N$ of the type of \ref{Theorem: Conditions imposed on i} should exist (see subsection \ref{Subsection: on existence of N}).  So we may extend the above question and ask, whether $ W $ must be generated from the embedding $ i\colon V\to  N $ as in section \ref{Section: The general Framework}.

\begin{question}
	What are the possibilities for non-normal $\kappa-$complete ultrafilters in $V[G]$?
\end{question}

Recall that, given $ i\colon V\to N $ and a measure $ W $ generated from it as in theorem \ref{Theorem: Conditions imposed on i}, the assumption that $ \leq^*_{Q_{\alpha}} = \leq_{Q_{\alpha} } $ holds for a final segments of $ \alpha\in \Delta $ suffices for $ N=M $ (where $ M $ is the ground model of $ \mbox{Ult}\left( V\left[G\right], W \right) $).

\begin{question}
 Suppose that for unboundedly many $\alpha<\kappa$,  $\leq_\alpha\not=\leq_\alpha^*$. Is then $M\neq N$?	
\end{question}

\begin{question}
 What are the exact conditions on $Q_\alpha$'s that insure $M=N$?
\end{question}

In section \ref{Section: on jwkappa above jukappa} we studied sufficient and necessary conditions for having $ j_W(\kappa) > j_U(\kappa) $. In proposition \ref{proposition: set unbounded in jWkappa}, we proved, under the assumption that $ \kappa $ is a limit of cardinals $ \alpha<\kappa $ which are all $ \kappa $-strong, that there are measures $ U\in \mathcal{K} $ and $ W\in \mathcal{K}\left[G\right] $ on $ \kappa $ extending $ U $,   such that $ j_W(\kappa) > j_U(\kappa) $.

\begin{question}
Is the assumption that $ \kappa $ is a limit of $ \kappa $-strong cardinals really necessary?
\end{question}

\begin{question}
In theorems \ref{Theorem: Measure generated from i restrics to iteration of N} and \ref{Theorem: Restrictions of ultrapwers with simply generated measures}, can we omit the assumption that the normal measures used in the iteration $ P = P_{\kappa} $ below $ \kappa $ are simply generated? How is the structure of $ j_W\restriction_{V} $ influenced from such a change?
\end{question}

\end{document}